\title{A monopole invariant for families of contact structures}
\documentclass[12pt , a4paper]{article}  

\usepackage{authblk} 

\author{Juan Mu\~{n}oz-Ech\'{a}niz}
\date{ }
\affil{Department of Mathematics, Columbia University}    
\usepackage[nochapters]{classicthesis} 
\usepackage{xfrac}
\usepackage{amsmath}
\usepackage[
backend=biber,
style=alphabetic,
sorting=nyt, doi=false, isbn=false, url=false, eprint=false
]{biblatex}
\usepackage{amsfonts}
\usepackage{amsthm}
\usepackage{upgreek}
\usepackage{tikz-cd}
\usepackage[margin=2.5cm]{geometry}
\usepackage{enumerate}
\usepackage{mathabx} 
\usepackage{tikz-cd}
\usepackage{graphicx}
\usepackage{stmaryrd}
\usepackage{dirtytalk}

 \usepackage{titlesec}

\titleformat{\subsubsection}[runin]
       {\normalfont\bfseries}
       {\thesubsubsection}
       {0.5em}
       {}
       [.]

\titleformat{\subsection}[runin]
       {\normalfont\bfseries}
       {\thesubsection}
       {0.5em}
       {}
       [.]
       
\newcommand{\inner}[2]{\langle #1 , #2 \rangle }
\newcommand{\R}{\mathbb{R}}
\newcommand{\C}{\mathbb{C}}
\newcommand{\Z}{\mathbb{Z}}
\newcommand{\Q}{\mathbb{Q}}

\newcommand{\ogr}[1]{\widetilde{\mathrm{Gr}}_{2}( #1)}

\newcommand{\s}{\mathfrak{s}}
\newcommand{\afr}{\mathfrak{a}}
\newcommand{\bfr}{\mathfrak{b}}
\newcommand{\cfr}{\mathfrak{c}}

\makeatletter
\newcommand{\dirlim@}[2]{%
  \vtop{\m@th\ialign{##\cr
    \hfil$#1\operator@font lim$\hfil\cr
    \noalign{\nointerlineskip\kern1.5\ex@}#2\cr
    \noalign{\nointerlineskip\kern-\ex@}\cr}}%
}
\newcommand{\dirlim}{%
  \mathop{\mathpalette\dirlim@{\rightarrowfill@\textstyle}}\nmlimits@
}

\newcommand{\HMto}{\widecheck{HM}}
\newcommand{\HMtilde}{\widetilde{HM}}
\newcommand{\fc}{\mathbf{Fc}}
\newcommand{\fctilde}{\widetilde{\mathbf{Fc}}}

\makeatother

\theoremstyle{plain}
\newtheorem{Theorem}{Theorem}[section]
\newtheorem*{thm*}{Theorem}

\newtheorem{Corollary}[Theorem]{Corollary}
\newtheorem{Lemma}[Theorem]{Lemma}
\newtheorem*{Lemma*}{Lemma}
\newtheorem*{Proposition*}{Proposition}
\newtheorem{Proposition}[Theorem]{Proposition}
\newtheorem{Slogan}[Theorem]{Slogan}

\newtheorem{Question}[Theorem]{Question}

\theoremstyle{definition}
\newtheorem{Definition}[Theorem]{Definition}
\newtheorem{Example}[Theorem]{Example}
\newtheorem{Remark}[Theorem]{Remark}

\DeclareMathAlphabet{\pazocal}{OMS}{zplm}{m}{n}

\begin{document}

\maketitle

\abstract{We use monopole Floer homology to study the topology of the space of contact structures on a $3$-manifold. Our main tool is a generalisation of the Kronheimer--Mrowka--Ozsváth--Szabó contact invariant to an invariant for families of contact structures, and we establish foundational results that describe the interaction between this invariant and the module structure in monopole Floer homology. We apply these results in several examples of contact manifolds, such as the links of non-rational surface singularities, to deduce several applications. Namely, we are able to obstruct the existence of sections of a natural fibration over the $2$-sphere whose total space is the space of contact structures on the $3$-manifold, and from this we are able to detect the existence of exotic loops of contact structures on contact 3-manifolds with convex sphere boundary.
}

\section{Introduction}

Let $Y$ always denote a smooth, closed, connected and oriented $3$-manifold. A (co-oriented and positive) \textit{contact structure} on $Y$ is a co-oriented $2$-plane field on $Y$ which is ``maximally non-integrable" in the following sense: for any $1$-form $\alpha$ with $\xi = \mathrm{ker} \alpha$ we have that $\alpha \wedge d \alpha$ is a positive volume form on $Y$. A $1$-form $\alpha$ such that $\xi = \mathrm{ker} \alpha$ as a co-oriented distribution is a called a \textit{contact form} for $\xi$.

A fundamental problem in low-dimensional topology is to classify contact structures on a given $3$-manifold $Y$. Much progress has been achieved in this direction and two tools are particularly useful here: the monopole or Heegard Floer homology groups \cite{KM,OS} together with their associated contact invariant \cite{monocont,monolens,OScontact} and Giroux's theory of convex surfaces \cite{convexite}. Beyond the classification problem, the topology of contact structures on $3$-manifolds as they vary in \textit{families} remains quite mysterious. Other than some small examples \cite{GGP04,Bou06,DG10,GK14,GM17,presas,EM} little is known is general. 

The contact invariant is an element $\mathbf{c}(\xi )$ in the monopole Floer homology group $\widecheck{HM}_\ast (- Y )$ canonically attached to a contact structure $\xi $ on $Y$ (well-defined up to $\pm$ sign ambiguity). It first appeared in \cite{monolens}, where it was used to establish that rational homology $3$-spheres carrying taut foliations aren't $L$-spaces, but the essential part of the construction of $\mathbf{c}(\xi_0 )$ was carried out in \cite{monocont}. The invariant $\mathbf{c}(\xi )$ enjoys several desirable properties, such as vanishing for overtwisted contact structures \cite{MR,OScontact,mariano}, non-vanishing for strongly symplectically fillable contact structures \cite{HFgenusbounds,mariano} or naturality under strong symplectic cobordisms \cite{MR,mariano}. For the reader more familiar with the other incarnations of the monopole Floer homology groups, namely the Heegaard Floer homology or Embedded Contact homology groups, we note that under the isomorphisms between the three theories (see \cite{tauECH},\cite{KLT}, \cite{CGH} and subsequent papers) the contact invariants defined in each theory are known to agree (see Theorem 1.1 in \cite{taubes5} and Theorem 6.2.4 in \cite{CGH1}). For convenience, we include here a dictionary between some of the monopole and Heegard-Floer groups: $ \widecheck{HM} = HF^{+}$, $\widehat{HM} = HF^{-}$, $\widetilde{HM} = \widehat{HF}$.

With a view towards the general problem of understanding the topology of the space of contact structures on a $3$-manifold, in this article we construct a \textit{families} version of the contact invariant. Before describing its formal properties (Theorems \ref{mainthm},\ref{formula},\ref{thmtriangles1}) we present some new applications of this invariant. 

\subsection{Non-existence of winding families of contact structures}

Fix a contact structure $\xi $ on $Y$. We denote $\mathcal{C}(Y, \xi )$ the path-component of $\xi$ in the space of contact structures on $Y$ (here we equip the space of contact structures with Whitney's $C^\infty$ topology). By Gray's stability Theorem \cite{geigesbook} any two contact structures $\xi_0 , \xi_1$ in $\mathcal{C}(Y, \xi )$ are \textit{isotopic}, that is, there exists a path of diffeomorphisms $f_t$ of $Y$ starting at $f_0 = \mathrm{id}$ such that $(f_1)_\ast \xi_0 = \xi_1 $. After a choice of point $p \in Y$ we can consider the \textit{evaluation map}
\begin{align}
 ev : \mathcal{C}(Y, \xi ) \rightarrow S^2  \label{ev}
\end{align}
which sends a contact structure $\xi^\prime$ to its plane $\xi^{\prime}(p)$ at the point $p$, with $S^2$ regarded as the space of co-oriented $2$-planes in $T_p Y \cong \mathbb{R}^3$. One can show that the map $ev$ is a \textit{fibration}. If $B \subset (Y , \xi )$ is a Darboux ball centered at $p$, the fibre of $ev$ is homotopy equivalent to the subspace $\mathcal{C}(Y , \xi , B ) \subset \mathcal{C}(Y , \xi )$ consisting of contact structures which agree with $\xi$ over $B$. In other words, $\mathcal{C}(Y, \xi , B )$ consists of those contact structures which agree with the standard model $dz - ydx$ over the ball $B$. 

\begin{Question}\label{section?}
When does the evaluation map $ev : \mathcal{C}(Y, \xi ) \rightarrow S^2$ admit a section ? 
\end{Question}

This is equivalent to finding a family $T \in \pi_2 \mathcal{C}(Y, \xi )$ such that the degree of the mapping $S^2 \xrightarrow{T} \mathcal{C}(Y, \xi ) \xrightarrow{ev} S^2$ is one (a \textit{winding} $S^2$-family of contact structures).

When $(Y, \xi )$ is either the tight $S^1 \times S^2$, the tight $S^3$, or a quotient of the latter by any finite subgroup $\Gamma \subset \mathrm{SU}(2)$ then the evaluation map does admit sections (see \S\ref{examples}). In general, the "flexible" version of Question \ref{section?} has a simple answer: on the space of co-oriented $2$-plane fields $\Xi (Y,\xi )$ in the connected component of a given one $\xi$, the evaluation map $ev : \Xi (Y, \xi ) \rightarrow S^2$ admits a section whenever the Euler class of $\xi$ vanishes (see Example \ref{flexible}). In contrast with these results, as an application of our invariant for families of contact structures we show that the location of the contact invariant $\mathbf{c}(\xi ; \mathbb{Q} ) \in \widecheck{HM}_\ast (- Y  ; \mathbb{Q})$ (here we use coefficients in $\mathbb{Q}$) obstructs the existence of sections of the evaluation map:

\begin{Theorem}\label{propsection}
If $\mathbf{c}(\xi ; \mathbb{Q} )$ does \textbf{not} lie in the image of the $U$ map 
\[U : \HMto_\ast (-Y ; \mathbb{Q} ) \rightarrow \HMto_\ast (-Y ; \mathbb{Q} )
\]
then the evaluation map $ev : \mathcal{C}(Y , \xi ) \rightarrow S^2$ does \textbf{not} admit sections, even after pullback by any homotopically non-trivial map $S^2 \rightarrow S^2$.
\end{Theorem}

The criterion $\mathbf{c}(\xi ; \mathbb{Q}) \notin \mathrm{Im} U$ applies in large classes of examples, including the Milnor fillable contact structures on the Seifert-fibered integral homology spheres (all except for $S^3$ and the Poincaré sphere $\Sigma (2,3,5)$). See [\cite{geigesspheres},Theorem 1.10] for a related result.

\subsection{Non-trivial loops of contact structures and contact mapping classes on the complement of a ball}

From Theorem \ref{propsection} one can easily deduce the existence of exotic phenomena in $3$-dimensional contact topology, see \cite{GGP04,Bou06} for previously known examples. Below we consider the compact contact $3$-manifold $(Y \setminus B , \xi )$ with convex sphere boundary obtained by removing a Darboux chart $B$ from $(Y, \xi )$.

\begin{Corollary}\label{loop}
If $\mathbf{c}(\xi ; \mathbb{Q} ) \notin \mathrm{Im}U$, then the Darboux ball complement $Y\setminus B $ admits a homotopically \textbf{non-trivial} loop of contact structures based at $\xi$ and agreeing with $\xi$ over $\partial (Y \setminus B)$. Furthermore, this loop of contact structures has infinite order. If, in addition, the Euler class of $\xi$ vanishes, then this loop is homotopically trivial as a loop of co-oriented $2$-plane fields on $Y \setminus B$ agreeing with $\xi$ over $\partial (Y \setminus B )$. 
\end{Corollary}

\begin{proof}
When the evaluation map admits no sections, then the connecting map in the long exact sequence in homotopy of the fibration $ev$ gives an injective map $\delta : \mathbb{Z} = \pi_2 S^2 \hookrightarrow \pi_1 \mathcal{C}(Y, \xi , B )$, and hence the loop $\delta (1)$ has infinite order. When the Euler class of $\xi$ vanishes then the loop $\delta (1)$ becomes trivial in $\pi_1$ of the corresponding space of co-oriented $2$-plane fields $\Xi (Y, \xi , B ) = \{ \xi^\prime \in \Xi (Y, \xi ) \, | \, \xi^\prime = \xi \text{ over } B \}$ by Example \ref{flexible} below.
\end{proof}


\subsection{The families contact invariant}

We now discuss the main technical result of this article, which describes formal properties of an invariant for families of contact structures generalising $\mathbf{c}(\xi )$. 

Throughout we fix a coefficient ring $R$ which we assume is commutative and unital. Associated to a $3$-manifold with a spin-c structure $\s$ one can associate various flavors of the monopole Floer groups (see \cite{KM,monolens}); the ones relevant to us will be $\HMto_\ast (Y, \s  ; R )$, $\HMtilde_\ast (Y, \s  ; R)$. A contact structure $\xi$ (or even just a co-oriented plane field) induces a spin-c structure denoted $\s_\xi$. A common feature of all the flavors is a canonical grading by the set of homotopy classes of co-oriented $2$-plane fields $\xi$ inducing the spin-c structure $\mathfrak{s}$, which we denote $\pi_0 \Xi (Y , \mathfrak{s} )$ and which carries a natural $\mathbb{Z}$-action. This induces a relative $\mathbb{Z} / \mathrm{div}(c_1 (\mathfrak{s}) )\mathbb{Z}$-grading where $\mathrm{div}( c_1 (\mathfrak{s} )) \in \mathbb{Z}_{\geq 0}$ is the divisibility of the first Chern class $c_1 (\s )$. 


The most basic version of our families invariant is a map of $R$-modules
\begin{align}
\fc : H_{\ast}(\mathcal{C}(Y , \xi ); \Lambda_R ) \rightarrow \widecheck{HM}_{\ast}(-Y , \s_{\xi}; R) \label{Uintro}
\end{align}
where $H_{\ast} (\mathcal{C}(Y , \xi ) ; \Lambda_R  )$ is the singular homology group of $\mathcal{C}(Y , \xi )$ with coefficients in a certain \textit{local system} $\Lambda_R$ of free $R$-modules of rank $1$ over the space $\mathcal{C}(Y , \xi )$. We have $\Lambda_R = \Lambda_\Z \otimes_{\Z} R$ where $\Lambda_\Z  $ is the local system of $\mathbb{Z}$-modules associated to the determinant line bundle of certain family of Fredholm operators parametrised by $\mathcal{C}(Y, \xi )$ (see Definition \ref{orientationlocalsystem}). In particular, if the characteristic of $R$ is two, then the local system $\Lambda_R$ is trivial. 


If we choose one of the two generators of the $\mathbb{Z}$-module $\Lambda_{\mathbb{Z}} (\xi )$ given by the fiber of $\Lambda_{\mathbb{Z}}$ over the point $\xi$, then this fixes the sign of the usual contact invariant $\mathbf{c}(\xi ; R ) \in \widecheck{HM}_\ast (-Y  , \s_{\xi } ; R) $. In addition, it also picks out a preferred generator, denoted by $1_R$, for the $R$-module $H_0 ( \mathcal{C}(Y , \xi) ; \Lambda_R )$. The element $1_{\Z}$ is either non-torsion or has order two, according as to whether the local system $\Lambda_\Z$ is trivial over $\mathcal{C}(Y , \xi )$ or not, respectively.

The Floer group $\widecheck{HM}_{\ast}(-Y , \s )$ carries the structure of a module over the graded $R$-algebra
\begin{align}
\mathbb{A}(R) = R [ U] \otimes_\Z \Lambda^{\ast} \big( H_1 (Y ; \Z ) / \mathrm{torsion} \big) \label{algebra}
\end{align}
where $U$ decreases grading by $2$ and $H_1 (Y ; \Z ) / \mathrm{torsion}$ decreases grading by $1$. In analogy with this, we will see that $H_{\ast} ( \mathcal{C}(Y , \xi ) ; \Lambda_R)$ can also be endowed with a natural module structure over the graded algebra $\mathbb{A}(R)$. In particular, the action of $U $ on $H_\ast (\mathcal{C}(Y , \xi ) ; \Lambda_R )$ is defined in terms of the evaluation map (\ref{ev}) and the cap product operation 
\begin{align*}
U : H_{\ast} (\mathcal{C}(Y , \xi ); \Lambda_R ) \rightarrow H_{\ast-2}(\mathcal{C}(Y , \xi ) ; \Lambda_R ) \quad T \mapsto T \cap  ev^\ast ([S^2]^\vee ) .
\end{align*}
We refer to Definition \ref{modulecontact} for the full action of $\mathbb{A}(R)$ on $H_\ast (\mathcal{C}(Y , \xi ) ; \Lambda_R )$. 

\begin{Remark}It follows that $U^2 = 0$ on $H_{\ast}(\mathcal{C}(Y , \xi ) ; \Lambda_R)$, so the latter is really a module over the graded $R$-algebra $$ R [U] / (U^2 ) \otimes_{\mathbb{Z}} \Lambda^\ast \big( H_1 (Y ; \mathbb{Z} )/ \mathrm{torsion} \big) .$$
\end{Remark}

The main technical result of the article is the following

\begin{Theorem}\label{mainthm}
There exists a "families contact invariant" given by a collection of $R$-module maps
\begin{align}
\fc : H_{j} (\mathcal{C}(Y , \xi ); \Lambda_R) \rightarrow \widecheck{HM}_{[\xi]+j}(-Y ,  \s_{\xi} ; R ) \quad , \quad j\geq 0 \label{invtmap}
\end{align}
which are natural with respect to orientation preserving diffeomorphisms and satisfy the following properties:
\begin{enumerate}[(a)]
\item[(A)] The $j = 0$ map recovers the usual contact invariant: $\fc (1_{R}) = \mathbf{c}(\xi ; R)$.

\item[(B)] $\fc$ is a map of graded $\mathbb{A}(R)$-modules: $\fc(a \cdot T ) = a \cdot \fc(T )$ for $a \in \mathbb{A}(R)$ and $T \in H_{\ast}(\mathcal{C}(Y , \xi ); \Lambda_R )$.
\end{enumerate}
\end{Theorem}

\begin{Remark}
\textit{Naturality.} The above assertion on naturality has the following meaning. Let $f$ be an orientation-preserving diffeomorphism of $Y$, and let $\xi_{1}$ be the contact structure obtained by pulling back another one $\xi_0$, $f^{\ast} \xi_0 = \xi_{1}$. By pulling back we have a homeomorphism $F = f^{\ast} : \mathcal{C}(Y,\xi_0 ) \xrightarrow{\cong} \mathcal{C}(Y , \xi_{1} )$. The assertion is that then there is a canonical isomorphism of local systems $\eta : \Lambda_{\Z } \xrightarrow{\cong} F^{\ast} \Lambda_{\Z}$ such that the following diagram (where the vertical arrows are isomorphisms) commutes

\begin{center}
\begin{tikzcd}
H_{\ast} ( \mathcal{C}(Y , \xi_0 ); \Lambda_{R} ) \arrow{d}{(F , \eta )_{\ast} } \arrow{r}{\fc_{0}}  & \widecheck{HM}_{\ast}(-Y , \s_{\xi_0}; R ) \arrow{d}{f_{\ast}}\\
H_{\ast} ( \mathcal{C}(Y , \xi_{1} ) ; \Lambda_{R} ) \arrow{r}{\fc_1} & \widecheck{HM}_{\ast}(-Y , \s_{\xi_{1}} ; R ).
\end{tikzcd}
\end{center}
\end{Remark}

\begin{Remark} \textit{Criterion for triviality of $\Lambda_{R}$}. It is unclear to the author whether $\Lambda_{R}$ can be non-trivial. However, a simple criterion is available:

\begin{Corollary}\label{triviallocalsystem}
Suppose the contact invariant $\mathbf{c}(\xi ; \Z) \in \widecheck{HM}(-Y , \s_{\xi} ; \Z)$ is not $2$-torsion, i.e. $2 \mathbf{c}(\xi ; \Z) \neq 0$. Then $\Lambda_{\Z}$ is trivial. 
\end{Corollary}
\begin{proof}
By Theorem \ref{mainthm}(A) it follows that $\fc(1_{ \mathbb{Z}})$ is not $2$-torsion, and hence that $H_0 (\mathcal{C}(Y, \xi ); \Lambda_\Z)$ is isomorphic to $\Z$ rather than $\Z/2\Z$. Hence $\Lambda_{\Z}$ is trivial.
\end{proof}

This criterion applies in many cases of interest. For instance, whenever the contact structure admits a \textit{strong symplectic filling}, in which case one has $\mathbf{c}(\xi ; \mathbb{Q} ) \neq 0$ already \cite{HFgenusbounds}.

\end{Remark}
\begin{Remark} \textit{Sign-ambiguity}.
Even if the local system $\Lambda_{\Z}$ over $\mathcal{C}(Y , \xi )$ is trivial, there is no canonical choice of generator of the $\mathbb{Z}$-module $\Lambda_{\mathbb{Z}}(\xi)$ for a given contact structure $\xi$. In fact, Lin--Ruberman--Saveliev \cite{lin-ruberman-saveliev} show that there is no way of fixing the sign so that the usual contact invariant $\mathbf{c}(\xi)$ becomes natural with respect to orientation-preserving diffeomorphisms of $Y$. Indeed, they show that the unique tight contact structure on $Y = -\Sigma(2,3,7)$ admits a contactomorphism $f$ which reverses the sign of $\mathbf{c}(\xi ; \Z)$ (i.e. $f_{\ast} \mathbf{c}(\xi ; \Z) = - \mathbf{c}(\xi ; \Z )$). We also note that the local system $\Lambda_{\Z}$ is trivial in this example, because this contact structure has a strong symplectic filling.
\end{Remark}

Going back to our motivating Question \ref{section?}, observe that the existence of a section of $ev$ is equivalent to the surjectivity of the \textit{degree map} $\mathrm{deg}:= ev_\ast : \pi_2 \mathcal{C}(Y , \xi ) \rightarrow \pi_2 S^2 = \Z$. This map can also be defined homologically, and gives a map $\mathrm{deg}( - ; R ) : H_2 (\mathcal{C}(Y , \xi ) ; R ) \rightarrow H_2 (S^2 ; R ) = R$. Theorem \ref{mainthm}(B) allows us to recover the degree of a family from the $U$ action on its families invariant:

\begin{Corollary}\label{formula}
If $\Lambda_\Z$ is trivial and a trivialisation is fixed, then for any $T \in H_2 (\mathcal{C}(Y, \xi );R)$ we have $U \cdot \fc(T ) = \mathrm{deg}(T ) \cdot \mathbf{c}(\xi ;R )$.
\end{Corollary}

\begin{proof}[Proof of Theorem \ref{section?}] 
If $\mathbf{c}(\xi ; \Q ) \neq 0$ then $\Lambda_\Z$ is trivial by Corollary \ref{triviallocalsystem}. A section $s : S^2 \rightarrow \mathcal{C}(Y, \xi )$ of $ev$ would yield a family $T := s_\ast [S^2] \in H_2 ( \mathcal{C}(Y , \xi ) ; \mathbb{Q} )$ with $\mathrm{deg}(T ; \mathbb{Q}) = 1$. Then by Corollary \ref{formula} we have $\mathbf{c}(\xi ; \Q ) = U \cdot \fc (T ) \in \mathrm{Im} U $. The argument for the non-existence of sections of the pullback of $ev$ by a homotopically non-trivial map $S^2 \rightarrow S^2$ is similar.
\end{proof}

\subsection{Examples}\label{examples}

We now give examples of contact $3$-manifolds $(Y, \xi )$ such that the criterion $\mathbf{c}(\xi ; \mathbb{Q} ) \notin \mathrm{Im}U$ applies. All of the examples below are irreducible $3$-manifolds, and in many of them $\xi$ has vanishing Euler class.

\begin{Example}(Links of singularities)
The simplest example is the Brieskorn sphere $$\Sigma (p,q,r) = \big\{ (x,y,z) \in \mathbb{C}^3 \, | \, x^p + y^q + z^r = 0 \text{  and  } |x|^2 + |y|^2 + |z|^2 = \epsilon \big\} \, , $$ where $\epsilon \in \mathbb{R}_{> 0}$ is small and $p,q,r \geq 1$ are integers with 
\begin{align}
1/p + 1/q + 1/r < 1 \, .\label{range}
\end{align}
These $3$-manifolds are rational homology spheres with infinite fundamental group \cite{milnor}. They carry a natural contact structure given by the complex tangencies $\xi_{\mathrm{sing}} = T \Sigma (p,q,r) \cap i T \Sigma (p,q,r)$. 

More generally, we consider an isolated normal surface singularity germ $(X, o )$ and the contact $3$-manifold $(Y, \xi_{\mathrm{sing}} )$ given by the \textit{link} of the singularity. The contact structure $\xi_{\mathrm{sing}}$ is called the \textit{Milnor fillable} contact structure associated to the singularity. It is known that such $3$-manifolds $Y$ are irreducible \cite{neumann}. Provided that $Y$ is a rational homology sphere, then the following are equivalent statements, as shown by Bodnár--Plamenevskaya \cite{plamenevskaya} and Némethi \cite{nemethi}:
\begin{enumerate}
    \item[(a)] $\mathbf{c}(\xi_{\mathrm{sing}}  ) \notin \mathrm{Im} U$
    \item[(b)] $Y$ is not an $L$-space
    \item[(c)] $(X, o )$ is not a rational singularity.
\end{enumerate}
It is known that any \textit{rational} singularity $(X,o)$ given by a single equation, or more generally a complete intersection, must be one of the \textit{ADE singularities} \cite{nemethinotes}. The links of the latter singularities have finite fundamental group, and so it follows that the Brieskorn singularity $x^p + y^q + z^r = 0$ with $p,q,r \geq 2$ in the range given by (\ref{range}) are not rational. Hence, by the equivalence $(a) \iff (c)$ the contact $3$-manifolds $(\Sigma (p,q,r), \xi_{\mathrm{sing}})$ in the range (\ref{range}) satisfy $\mathbf{c} (\xi_{\mathrm{sing}} ) \notin \mathrm{Im}U$. Similarly, any Seifert fibered integral homology sphere, excluding $S^3$ or the Poincaré sphere $\Sigma (2,3,5)$, carries a Milnor fillable contact structure $\xi_{\mathrm{sing}}$ with $\mathbf{c} (\xi_{\mathrm{sing}} ) \notin \mathrm{Im}U$. Indeed, the non-trivial Seifert fibered integral homology spheres are given by the manifolds $ \Sigma ( p_1, p_2, \ldots , p_n )$, where $p_i \geq 2$ are pairwise coprime integers and $n \geq 3$, where $\Sigma (p_1, p_2 , \ldots , p_n )$ is the link of the isolated singularity described by the complete intersection $f_1 = \ldots = f_{n-2} = 0$ with $f_j = \sum a_{ij} x_{i}^{p_i}$ for sufficiently general coefficients $a_{ij} \in \C$ (these manifolds have infinite fundamental group with the exception of $\Sigma (2,3,5)$ \cite{NR}). Alternatively, this can be deduced from $(a) \iff (b)$, since it is known that no non-trivial Seifert fibered integral homology sphere is an $L$-space except the Poincaré sphere $\Sigma (2,3,5 )$ by \cite{eftekhary} (see also \cite{LS,MOY}).

To spell out one concrete example, for $\Sigma (2,3,7 )$ we have\footnote{The gradings in Floer homology for the examples in this section are taken appropriately \textbf{shifted} so that the contact invariant $\mathbf{c} (\xi )$ lies in degree $0$. Also, all identities involving contact invariants are understood to hold \textbf{up to signs}.}$$\HMto_\ast ( -\Sigma (2,3,7 ) , \s_{\xi_\mathrm{sing}} ; \mathbb{Z} ) \cong \mathbb{Z}_{(0)} \oplus \mathbb{Z}[U,U^{-1}]/U \mathbb{Z}[U ]$$ and the $U$ action is trivial on the $\mathbb{Z}_{(0)}$ summand. In this case $\mathbf{c} (\xi; \mathbb{Z} ) $ generates the summand $\mathbb{Z}_{(0)}$ and hence $\mathbf{c}(\xi ; \mathbb{Q}) \notin \mathrm{Im}U$.

\end{Example}

\begin{Example}
Several surgeries on the Figure Eight knot are hyperbolic (hence irreducible) and support contact structures with $\mathbf{c}(\xi ; \mathbb{Q}) \notin \mathrm{Im}U$. Contact structures on these manifolds have been classified by Conway and Min \cite{conway}. 
\end{Example}

\begin{Example}
 All but one of the $\frac{n(n-1)}{2}$ tight contact structures supported on $- \Sigma ( 2,3,6n-1)$ up to isotopy, classified by Ghiggini and Van Horn-Morris \cite{ghiggini-vanhorn-morris}, satisfy $\mathbf{c}(\xi ; \mathbb{Q}) \notin \mathrm{Im}U$.
\end{Example}

For comparison, we include also examples for which the evaluation map does admit sections, and hence with $\mathbf{c}(\xi) \in \mathrm{Im}U$.

\begin{Example}\textit{ADE singularities.}\label{linksADE}
Consider the flat hyperk\"{a}hler structure $(g , I_1 , I_2 , I_3 )$ on $\R^4$. The radial vector field $v = x \partial_x + y \partial_y + z \partial_z + w \partial_w $ in $\mathbb{R}^4$ is Liouville for all symplectic structures in the family $\omega_t = \sum_{i = 1}^{3} t_i g ( I_i \cdot , \cdot )$ parametrised by $t \in S^2$ (i.e. $\mathcal{L}_{v} \omega_t = \omega_t $) and $v$ is transverse to $S^3 \subset \mathbb{R}^4$. Thus there is a family of contact forms $\alpha_t$ on $S^3$ given by $\alpha_t= \iota_{v} \omega_t$ which provides a section of $ev$ on tight $S^3$. Since this family of contact structures is $\mathrm{SU}(2)$-invariant, we have also constructed a section of $ev$ on the quotients of tight $S^3$ by a finite subgroup $\Gamma \subset \mathrm{SU}(2)$. The contact manifolds $S^3 / \Gamma$ are precisely the the links of a special class of \textit{rational} singularities: the \textit{ADE singularities}. These singularity links include, among others, the lens spaces $L(p,p-1)$ or the Poincaré sphere $\Sigma (2,3,5 )$. Let $\xi$ be any contact structure in the $S^2$-family $ \xi_t = \mathrm{ker} \alpha_t$. We have $\HMto_\ast (- S^3 / \Gamma , \s_{\xi} ; \mathbb{Z}) \cong \mathbb{Z}[U , U^{-1}]/U \mathbb{Z}[U]$ and $\mathbf{c}(\xi ) = 1$. If $T$ denotes the $S^2$-family of contact structures given by the $\xi_t$ then from Theorem \ref{formula} we deduce $\fc (T) = U^{-1}$, and $U \cdot \fc (T) = \mathbf{c}(\xi )$. 
\end{Example}

\begin{Example} \textit{Tight $S^1 \times S^2$}. Let $\xi_{\mathrm{tight}}$ be a tight contact structure on $S^1 \times S^2$. We consider two families of contact structures in $\mathcal{C}(Y, \xi_{\mathrm{tight}} )$. First, consider the family $T \in H_2 ( \mathcal{C}(S^1 \times S^2 , \xi_{\mathrm{tight}} ); \mathbb{Z} )$ of contact structures $\xi_t$ parametrised by $t \in S^2$ given by the kernels of $\alpha_t = \sum_{i = 1}^{3} t_i \alpha_i $ where
\begin{align*}  
    \alpha_1 = z d \theta + xdy - y dx \,\, , \,\,
    \alpha_2 = xd \theta + ydz - zdy \,\, , \,\,
    \alpha_3 = y d \theta + zdx - x dz .
\end{align*}
It is a simple exercise to check that this family provides a section for the evaluation map. 

As a $\mathbb{Z}[U]$ module we have $$\HMto_{\ast} (-S^1 \times S^2 , \s_{\xi_{\mathrm{tight}}} ; \mathbb{Z}) \cong \mathbb{Z}[U , U^{-1}]/U \mathbb{Z}[U] \otimes_{\mathbb{Z}} H_\ast (S^1 ; \mathbb{Z}) $$ 
We have $\mathbf{c} (\xi_{\mathrm{tight}}) = 1$, 
$\fc (T ) = U^{-1}$ and $U \cdot \fc (T) = \mathbf{c} (\xi_\mathrm{tight} )$ by Theorem \ref{formula}. 

\end{Example}

\begin{Example}(\textit{Evaluation of $2$-plane fields})\label{flexible}
We can compare Question \ref{section?} with its "flexible" analogue: consider now the evaluation map $\Xi (Y, \xi ) \rightarrow S^2$ on the space co-oriented $2$-plane fields, which is also a fibration. This evaluation map admits a section as long as the Euler class of $\xi $ vanishes. Indeed, in this case we may identify $\Xi (Y, \xi )$ with the space $\mathrm{Map}_0 (Y , S^2 )$ of null-homotopic smooth maps $Y \rightarrow S^2$. The evaluation mapping becomes identified with the obvious evaluation mapping on this latter space, this admits a section given by the constant maps $Y \rightarrow S^2$. 

The homotopy type of $\Xi (Y , \xi )$ is often well-understood. For instance, whenever $Y$ is an integral homology sphere then $\Xi (Y, \xi )\simeq \mathrm{Map}_0 (S^3, S^2 )$ \cite{hansen}. 
\end{Example}

\subsection{The $U$-map and families of contact structures}

We now describe a refinement of Theorem \ref{mainthm} in the case of the action of $U \in \mathbb{A} (R)$. For the remainder of this section we assume that the local system $\Lambda_{\Z}$ over $\mathcal{C}(Y , \xi )$ is \textit{trivial} (recall once more the criterion which ensures this, Corollary \ref{triviallocalsystem}) and fix a trivialization (i.e. a choice of generator of the $\mathbb{Z}$-module $\Lambda_{\mathbb{Z}}(\xi )$ ) so that the families invariant gives a map
$$
\fc : H_\ast ( \mathcal{C}(Y , \xi ) ) \rightarrow \HMto_\ast (- Y , \mathfrak{s}_\xi ).
$$

From now on we omit the ring $R$ from our notation if it is not essential. 

Going beyond Question \ref{section?}, one could ask how the homotopy type of the space $\mathcal{C}(Y , \xi )$ differs from that of $\mathcal{C}(Y , \xi , B )$. Often the latter has "simpler" topology. For example, for the tight contact structure $\xi$ on $S^3$ one has $\mathcal{C}(S^3 , \xi ) \simeq \mathrm{U}(2)$ whereas $\mathcal{C}(S^3 ,\xi, B ) \simeq \{ \ast \}$ \cite{EM}. At the homological level, the passage from $\mathcal{C}(Y, \xi , B )$ to $\mathcal{C}(Y, \xi )$ amounts to understanding how cycles in the total space of the fibration $ev$ intersect with the fibres, and this is encoded into the \textit{Wang exact sequence} for the fibration (\ref{ev}) (easily assembled from the Serre spectral sequence)
\begin{center}
\begin{tikzpicture}[baseline= (a).base]
\node[scale=1 , trim left=-7cm] (a) at (0,0){
\begin{tikzcd}
\cdots \arrow{r} & H_{\ast}(\mathcal{C}(Y , \xi )) \arrow{r}{U_B}  & H_{\ast-2}(\mathcal{C}(Y , \xi , B )) \arrow{r}{\chi}  &  H_{\ast-1}(\mathcal{C}(Y , \xi , B)) \arrow{r}{\iota_\ast} & \cdots 
 \end{tikzcd}
};
\end{tikzpicture}
\end{center}
In geometric terms, the map $U_B$ acts on a generic cycle in $\mathcal{C}(Y , \xi )$ by taking its intersection with the fibre of (\ref{ev}), and $\iota_{\ast}$ is the inclusion of the fibre. The map $\chi$ is the differential in the $E^2$ page of the spectral sequence. The map $H_{\ast}(\mathcal{C}(Y , \xi )) \xrightarrow{U} H_{\ast-2}(\mathcal{C}(Y , \xi ))$ defined earlier can be recovered from the diagram above as the composition $U = \iota_{\ast} \circ U_p $. 

On the Seiberg--Witten gauge-theory side one can find a structure analogous to the evaluation map $ev : \mathcal{C}(Y , \xi ) \rightarrow S^2$. The space of \textit{irreducible} configurations modulo gauge transformations $\mathcal{B}^{\ast}( Y ,  \s_{\xi} )$ also carries a partially-defined evaluation map 
\begin{align}
\mathcal{B}^{\ast} (Y ,  \s_{\xi} ) \dashrightarrow \mathbb{P}(S_p ) \cong \C P^1= S^2 \label{evgauge}
\end{align}
which assigns to the class of a configuration $(B , \Psi )$ the complex line in the spinor bundle fibre $S_p \approx \C^2$ spanned by $\Psi$ at the point $p$. The relevance of this evaluation map is its close relation with $U$-action on the Floer theory. Indeed, the action of $U$ on the Floer homology is defined as a sort of cap product with the first Chern class of a canonical complex line bundle $\mathcal{U} \rightarrow \mathcal{B}^{\ast} (Y ,  \s_{\xi} )$, with (\ref{evgauge}) arising as the map to $\mathbb{C}P^1$ determined by a certain ``pencil'' of hyperplanes in the class of the line bundle $\mathcal{U}$. Thus, resembling the contact case, this operation corresponds geometrically to taking intersections of moduli of Floer trajectories with the fibres of (\ref{evgauge}). Similarly to the Wang long exact sequence, on the Floer theory we have a \textit{Gysin exact sequence}
\[
\begin{tikzcd}
\cdots \arrow{r}{p} & \HMto_{\ast}(-Y , \mathfrak{s}) \arrow{r}{U} & \HMto_{\ast-2} (-Y, \mathfrak{s}) \arrow{r}{j} & \HMtilde_{\ast-1}(-Y , \mathfrak{s} ) \arrow{r}{p} & \cdots  
\end{tikzcd} 
\]

involving the "tilde" group $\HMtilde_\ast (Y, \mathfrak{s} )$. Whereas the group $\HMto_\ast (Y, \mathfrak{s} )$ arises "formally" as the $S^1$-equivariant Morse homology of the Chern-Simons functional (and thus carries the module structure over the polynomial algebra $R[U] \cong H_{S^1}^\ast (\mathrm{point} ; R )$, with $U$ decreasing grading by $2$), the "tilde" flavor should be regarded as the usual (non-equivariant) Morse homology, and thus carries the structure of an $H_\ast (S^1 ) = R[\chi]/(\chi^2 ) $-module, with $\chi$ raising grading by one. From the Gysin sequence, the map $\chi$ can be recovered as $\chi = j p$.

The connection between the two evaluations (\ref{ev}) and (\ref{evgauge}) is seen by certain map which assigns canonical irreducible configurations to contact structures
\[
\begin{tikzcd}
f : \mathcal{C}(Y , \xi )  \arrow{r} &  \mathcal{B}^\ast (Y , \s_{\xi } ).
\end{tikzcd}
\]
Under the familiar identification $S^2 = \mathbb{C}P^1$ coming from spin geometry, the map $f$ intertwines our two evaluation maps (\ref{ev}) and (\ref{evgauge}). On a heuristic level, one should regard the families contact invariant $\fc$ as the "map induced by $f$ in homology" (one should be able to formalise this by working at the level of spectra, but we don't pursue this direction in this article). At this point, Theorem \ref{mainthm}(B) and the following refinement should be regarded as algebraic manifestations of the phenomenon just described:

\begin{Theorem} \label{thmtriangles1}
Associated to any closed contact $3$-manifold $(Y , \xi )$ with trivial local system $\Lambda_{\Z}$  there is a natural diagram which is commutative (up to signs)
\begin{center}
\begin{tikzpicture}[baseline= (a).base]
\node[scale=1, trim left=-7cm] (a) at (0,0){
\begin{tikzcd}
 \llap{ }\arrow{r}{p}  & \widecheck{HM}_{\ast} ( -Y , \s_{\xi }) \arrow{r}{U} &  \widecheck{HM}_{\ast- 2} ( -Y , \s_{\xi} ) \arrow{r}{j} & \widetilde{HM}_{\ast-1}(-Y , \s_{\xi } ) \arrow{r}{p} & \llap{ }\\ 
\llap{} \arrow{r}{\iota_\ast} & H_{\ast}(\mathcal{C}(Y , \xi )) \arrow{r}{U_B} \arrow{u}{\fc} & H_{\ast-2}(\mathcal{C}(Y , \xi , B )) \arrow{r}{\chi} \arrow{u}{\fc \cdot \iota_\ast} &  H_{\ast-1}(\mathcal{C}(Y , \xi , B )) \arrow{r}{\iota_{\ast}} \arrow{u}{\fctilde} & \llap{ } 
 \end{tikzcd}
};
\end{tikzpicture}
\end{center}
where the top row is the Gysin exact sequence, the bottom row is the Wang exact sequence of the fibration (\ref{ev}) and $\fctilde$ is another "families contact invariant" that we construct in \S \ref{trianglessection}.
\end{Theorem}

As a particular case, Theorem \ref{thmtriangles1} recovers a property about the contact invariant $\mathbf{c}(\xi )$ which is well-known from the Heegaard--Floer point of view: that $U \cdot \mathbf{c}(\xi ) = 0 $ and we have a canonical element $\widetilde{\mathbf{c}}(\xi ):= \fctilde (1) \in \HMtilde_{[\xi]} (-Y , \mathfrak{s}_{\xi} )$ such that $p \widetilde{\mathbf{c}}(\xi ) = \mathbf{c}(\xi )$. Conjecturally, the invariant $\widetilde{\mathbf{c}}(\xi )$ corresponds to the Heegaard--Floer contact invariant that takes values in $\widehat{HF}(-Y, \mathfrak{s}_\xi )$ defined in \cite{OScontact}.

The group $\HMtilde_\ast (-Y , \mathfrak{s}_{\xi })$ is an $R[\chi]/(\chi^2 )$ module, and so is $H_\ast ( \mathcal{C}(Y, \xi , B ) )$. One deduces from Theorem \ref{thmtriangles1} the

\begin{Corollary}
The $\fctilde$ is a map of $R [\chi ] / (\chi^2 )$ modules:
$$ \fctilde \cdot \chi = \chi \cdot \fctilde .$$
\end{Corollary}

In particular, we deduce from this and the diagram that
\begin{center}
    $\mathbf{c}(\xi) \in \mathrm{Im}U$ \textit{  if and only if  } $\chi \widetilde{\mathbf{c} }(\xi ) = 0 .$
\end{center}

\begin{Example}
Continuing with the examples from \S \ref{examples}, we have that as a $\mathbb{Z}[\chi]/(\chi^2)$ module
$$
\HMtilde_\ast (- \Sigma (2,3,7) , \mathfrak{s}_{\xi_\mathrm{sing}} ; \Z ) \cong \mathbb{Z}_{(0)} \oplus \mathbb{Z}[\chi]/(\chi^2)
$$
where $\chi$ acts trivially on the $\mathbb{Z}_{(0)}$ summand. Here $\widetilde{\mathbf{c}}(\xi_\mathrm{sing} ) $ generates the bottom grading of $\mathbb{Z}[\chi]/(\chi^2)$ and $\fctilde (\mathcal{O}_{\xi_{\mathrm{sing}}} )  = \chi \cdot \widetilde{\mathbf{c}}(\xi_\mathrm{sing})$ generates the top grading, where $\mathcal{O}_{\xi_{\mathrm{sing}}} := \chi \cdot [\xi_{\mathrm{sing}}] \in H_1 (\mathcal{C}( \Sigma (2,3,7) , \xi_{\mathrm{sing}} , B ) )$. 

\end{Example}

\subsection{Sketch of the construction of (\ref{invtmap})}\label{sketch}

We first review the construction of the usual contact invariant $\mathbf{c}(\xi )$ in monopole Floer homology (see \cite{monolens}). For this one studies the Seiberg-Witten equations on certain non-compact Riemannian manifold $Z^+$ which is formed by taking the symplectization $( K = [1 , + \infty)_s \times Y , \omega = d (\frac{s^2}{2} \alpha ) )$ of a contact form $\alpha$ for $\xi$ and gluing a cylindrical end $Z = (- \infty , 0] \times \R$ to $K$. The natural Riemannian metric to consider over the end $K$ is \textit{conical}, whereas over $Z$ we consider a cylindrical metric. One then considers a Seiberg-Witten moduli space $M([\afr] , \xi )$ over $Z^+$ in which solutions to the equations must approach a translation-invariant solution $[\afr]$ over $Z$ and over $K$ approach a \textit{canonical} solution $(A_{\xi }, \Phi_{\xi} )$ to the equations over $K$ which is provided by the symplectic structure $\omega$. The module $\widecheck{HM}_\ast (-Y, \s_{\xi} )$ is generated by translation-invariant solutions $[\afr]$, and enumerating points in rigid moduli spaces $M ([\afr] , \xi )$ produces a Floer cycle $\psi (\xi )$ such that $\mathbf{c}(\xi ) := [\psi (\xi )]$. The assignment $\xi \mapsto \mathbf{c}(\xi )$ can be regarded as a sort of topological field theory restriction map for the $4$-dimensional manifold $K$ with boundary $\partial K= -Y$, even if $K$ happens to be non-compact.

In order to define an invariant of a family of contact structures $\xi_t$, the simplest situation is that when the space parametrising the contact structures is contractible, let's say it is the $j$-simplex $\Delta$. One then considers a similar set of equations on $Z^+$ parametrised by $t \in \Delta$. The essential dependence of the equations on $t \in \Delta$ is through the ``boundary condition'' over the symplectic end $K$ given by the canonical solutions $(A_{\xi_t} , \Phi_{\xi_t} )$. In turn, the metrics and spin-c structures are constant (do not depend on $t \in \Delta$) over the cylinder $Z$. 
This enables us to set up moduli spaces over simplices $M([\afr] , \Delta ) \rightarrow \Delta$ much as before, and we extract an invariant of our $\Delta$-family in the standard fashion: by enumerating points in rigid moduli spaces $M([\afr] , \Delta )$ as $[\afr]$ varies, the counts assemble into a Floer chain $\psi (\Delta )$. The latter is not a cycle anymore, but the assignment $\Delta \mapsto \psi (\Delta )$ is a chain map and this yields (\ref{invtmap}) when passing to homology. 


\subsection{Outline}
The structure of the article is as follows. In \S \ref{background} we provide some necessary background on families of spin-c structures. In \S \ref{invtsection} we present the construction of the families contact invariant (\ref{invtmap}), from which Theorem \ref{mainthm}(A) follows. In \S \ref{modulesection} we define and interpret geometrically the module structures that Theorem \ref{mainthm}(B) refers to. The technical core of the article is \S \ref{proofmainsection}, where a neck-stretching argument is carried to prove Theorem \ref{mainthm} (B). In \S \ref{trianglessection} we prove Theorem \ref{thmtriangles1}. Throughout the article, various results regarding transversality with evaluation constraints, compactness and orientations of Seiberg-Witten moduli spaces are assumed. We discuss them in more detail in \S\ref{appendix}, which has the nature of an appendix.

\tableofcontents

\subsubsection*{Acknowledgements} The author is grateful to his advisor Francesco Lin for suggesting this problem and for the many useful discussions. This work also benefited from conversations with Francisco Presas, Eduardo Fern\'{a}ndez and Luya Wang. The author thanks Hokuto Konno and the anonymous referee for making corrections and suggestions that have greatly improved the quality of the exposition.
\section{Families of spin-c structures and irreducible configurations}\label{background}

In this section we discuss preliminary material regarding spin-c structures as they vary in families. 

\subsection{Basic notions about spin-c structures}

We let $M$ be an oriented manifold of dimension $n = 2m$ or $n = 2m+1$. Our case of interest is $n = 4$ or $n = 3$. 

\begin{Definition}
A \textit{spin-c structure} on $M$ is a triple $\s = (g , S, \rho )$ consisting of the following data:
\begin{enumerate}
    \item[(a)] a Riemannian metric $g$ on $M$
    \item[(b)] a unitary vector bundle $S \rightarrow M$ of rank $2^m$
    \item[(c)] a vector bundle map $\rho : T^\ast M \rightarrow \mathrm{Hom} (S , S )$ which is skew-adjoint $\rho (v )^{\ast} = - \rho (v )$ and satisfies the \textit{Clifford identity} $\rho(v )^2 = - |v |_{g}^{2} \cdot \mathrm{id}_{S}$ , for all $v \in T^{\ast}M$.
\end{enumerate} 
The bundle $S$ is referred to as the \textit{spinor bundle} of $\s$ and its sections are \textit{spinors}; the map $\rho$ is the \textit{Clifford multiplication} of $\s$.
\end{Definition}

The Clifford multiplication $\rho$ naturally extends to a vector bundle map from the complexified exterior algebra $\rho: \Lambda^{\bullet} T^{\ast} M \otimes \C \rightarrow \mathrm{Hom}(S, S)$ by linearity and the following rule: for a one-form $\alpha$ and a form $\beta$ $$\rho( \alpha \wedge \beta ) =\frac{1}{2} \big( \rho( \alpha ) \rho(\beta ) + (-1)^{ \mathrm{deg} \beta } \rho (\beta ) \rho (\alpha )\big).$$
From the canonical volume element $\omega$ determined from the metric $g$ we form the complex volume element $\omega_{\C} = i^{\lfloor \frac{n+1}{2} \rfloor}\omega \in \Gamma ( M , \Lambda^n T^\ast M \otimes \C )$. One sees that $\rho ( \omega_{\C} )^2 = 1$. In the case $n = 2m$ the bundle $S$ decomposes $S = S^+ \oplus S^- $ as the sum of the $\pm 1$-eigensubbundles of $\rho ( \omega_{\C} )$. Each $S^{\pm }$ has rank $2^{m-1}$ and these are referred to as \textit{positive} or \textit{negative} spinor bundles. In the case $n = 2m+1$ we \textit{require} in the definition of a spin-c structure that $\rho(\omega_{\C})$ acts on $S$ by $-1$.

If $X$ is an oriented manifold of dimension $2m$ with $\partial X = Y$ and we are given a spin-c structure $\s_X = (g_X , S_X , \rho_X )$ on $X$, we can \textit{restrict} it to $Y$ and obtain a spin-c structure $\s_X|_{Y} = (g_X|_Y , S_X^{+}|_Y , \rho_Y )$. Here $\rho_Y$ is defined by $\rho_Y (v) = \rho_X (n)^{-1} \rho_X (v)$, where $n$ stands for the unit outward normal to $Y$.

We now describe some further differential geometric notions associated with a spin-c structure:

\begin{Definition}
A unitary connection $A$ on the unitary bundle $S \rightarrow M$ is a \textit{spin-c connection} if the Clifford action $\rho : T^{\ast} M \rightarrow \mathrm{Hom}(S , S )$ is parallel with respect to the connection on $TM \otimes \mathrm{Hom}(S,S)$ induced by $A$ and the Levi-Civita connection of $g$.
\end{Definition}

There is a one-to-one correspondence between spin-c connections on $S$ and unitary connections on the associated line bundle $L = \mathrm{det} S^{+}$ if $n = 2m$ (and $L = \mathrm{det} S$ if $n = 2m+1$). The connection on $L$ induced by $A$ is denoted by $\hat{A}$, and the correspondence is just $A \mapsto \hat{A}$. Thus, the space of spin-c connections is an affine space over $\Omega^1 ( M ; i \R )$. 
\begin{Definition}
The \textit{Dirac operator} coupled with a spin-c connection $A$ is the differential operator $$D_{A} : \Gamma ( X , S ) \rightarrow (X , S )$$ defined by $D_{A} \Phi = \rho ( \nabla_{A} \Phi ),$ where the latter expression denotes the contraction of the $T^{\ast} X$ and the $S$ component of $\nabla_{A}$ using the Clifford action $\rho$. 
\end{Definition}

The differential operator $D_{A}$ is elliptic and self-adjoint. In the case $\mathrm{dim} M =n =  2m$, the Dirac operator decomposes $D_{A} = D_{A}^{+} \oplus D_{A}^{-}$ as a sum of two elliptic differential operators $D^{\pm}_{A} : \Gamma ( X , S^{\pm} ) \rightarrow \Gamma (X, S^{\mp} )$. 

\subsection{Changing the metric of a spin-c structure}

Given a spin-c structure $\s_0  = (g_0 , S_0 , \rho_0 )$ and a different Riemannian metric $g_1$ on $M$, there is a natural device for producing a new spin-c structure $(g_1 , S_1 , \rho_1 )$. We describe this now following \cite{BG}. 

Consider first a real finite-dimensional vector space $V$ equipped with two inner products $g_0, g_1$. Then there is a \textit{canonical linear isometry} $b_{g_1, g_0} : (V , g_0) \xrightarrow{\cong} (V , g_1)$, characterised by the property that it is positive and symmetric with respect to $g_0$. It is constructed as follows. Write $g_1 = g_0 (H \cdot , \cdot )$ for a (unique)  symmetric positive endomorphism $H$ of $(V , g_0)$. Then $b_{g_1, g_0} = H^{-1/2}$ is the required isometry. Finally, given two Riemannian metrics $g_0, g_1$ on a manifold $M$, the previous construction applies fibrewise to produce an isometry $b_{g_1 , g_0 } : (TM , g_0 ) \xrightarrow{\cong} (TM , g_1 )$.

\begin{Remark}The canonical isometry satisfies $b_{g_1 , g_0 }^{-1} = b_{g_0 , g_1 }$. Unfortunately, in general it is \textit{not} functorial: $b_{g_2, g_1 } \circ b_{g_1 , g_0 } \neq b_{g_2 , g_0}$ (see \cite{BG}).
\end{Remark}

This construction allows us to change the metric in a spin-c structure $\s_0 = (g_0 , S_0 , \rho_0 )$. Given another Riemannian metric $g_1$, we define $S_1 = S_0$ and $\rho_1 : T^\ast X \rightarrow \mathrm{Hom}(S_0 , S_0 )$ as $\rho_1 (v ) = \rho_0 (b_{g_1 , g_0 }^{\ast} v )$. This yields a new spin-c structure $\s_1 = (g_1, S_1, \rho_1 )$.

\begin{Definition}
Given two spin-c structures $\s_i = (g_i , S_i , \rho_i )$ ($i = 0,1$) on $M$, an \textit{isomorphism} between them consists of an isomorphism of unitary vector bundles $h : S_0 \xrightarrow{\cong} S$ such that $\rho_1 (v) = h \circ \rho_0 ( b_{g_1 , g_0}^{\ast} v ) \circ h^{-1}$ for all $v \in T^{\ast} X$.
\end{Definition}

It can be shown using Schur's Lemma that set of isomorphism classes\footnote{It is clear that "isomorphism" gives a reflexive and symmetric relation on the set of spin-c structures $(g,S,\rho )$; it can be shown using Schur's Lemma that it is also \textit{transitive}, even if $b_{g_2, g_1} \circ b_{g_1, g_0} = b_{g_2 , g_1}$ does not hold in general.} of spin-c structures on $M$ is a \textit{torsor} over the cohomology group $H^2 (M ; \Z )$. Given a unitary line bundle $Q$ over $M$, the action of $c_1 (Q) \in H^2 (M ; \Z )$ on the isomorphism class of the spin-c structure $\s = (g, S , \rho )$ is defined by $$ c_1 (Q) \cdot [\s ] = [ (g , S \otimes Q , \rho \otimes \mathrm{id}_Q ) ] .$$

\subsection{Irreducible configurations}

The space of configurations $(A, \Phi )$, where $A$ is a spin-c connection on $S$, and $\Phi \in \Gamma ( M , S^{+} )$ (resp. $\Phi \in \Gamma ( M , S )$ ) when $n = 2m$ (resp. $2m+1$) is denoted by $\mathcal{C}(M, \s)$. We equip $\mathcal{C}(M, \s )$ with the $C^\infty$ topology.  We denote by $\mathcal{C}^{\ast}(M, \s ) \subset \mathcal{C}(M, \s )$ the open subset of \textit{irreducible configurations}, namely those such that $\Phi$ is not identically vanishing on $M$. Configurations $( A , 0 )$ are called \textit{reducible}.

The automorphism group $\mathcal{G}$ of a spin-c structure $\s = (g, S, \rho )$ is referred to as the \textit{gauge group}. It can be shown using Schur's Lemma that $\mathcal{G}$ agrees with space of smooth mappings $\mathcal{G} = \mathrm{Map}(M , \mathrm{U}(1) ).$ We make $\mathcal{G}$ into a topological group by equipping it with the $C^\infty$ topology. There is a continuous $\mathcal{G}$-action on $\mathcal{C}(M, \s )$: given $v \in \mathcal{G}$ and configuration $(A, \Phi )$ we set
\begin{align*}
v \cdot (A, \Phi ) = (A - v^{-1} dv , v \Phi ) .
\end{align*}
The $\mathcal{G}$-action is free on $\mathcal{C}^{\ast}(X, \s )$, whereas it has stabiliser$\cong \mathrm{U}(1)$ at the reducible configurations.

\begin{Definition}
The \textit{configuration space} modulo gauge is the quotient space $\mathcal{B}(X , \s ) = \mathcal{C}( X , \s ) / \mathcal{G}$. The subspace $\mathcal{C}^{\ast}(X, \s )/ \mathcal{G} \subset \mathcal{B}(X , \s )$ is denoted $\mathcal{B}^{\ast}(X, \s )$.
\end{Definition}

The space $\mathcal{B}(X , \s )$ is Hausdorff. If an isomorphism $h : \s_0 \xrightarrow{\cong} \s_1$ of two spin-c structures on $X$ is given, there is an induced homeomorphism $\mathcal{B}(X , \s_0 ) \xrightarrow{\cong} \mathcal{B}(X , \s_1 )$ given by $(A, \Phi ) \mapsto (A_h, h ( \Phi ))$ where $A_h$ is the unique spin-c connection (for $\s_1 )$ such that $\widehat{(h( A ))} = \widehat{(A_h )}$.

\subsection{Families of spin-c structures and irreducible configurations}  \label{familiesirred}

We now consider continuously-varying families of spin-c structures $\s_t = (g_t , S_t , \rho_t )$ on a \textit{fixed} oriented smooth manifold $M$ parametrised by a "nice" connected topological space $T$. Note that the isomorphism class of the spin-c structures on $M$ given by $[\s_t ]$ is independent of $t \in T$. We denote such a $T$-family by $\underline{\s} =  ( \s_t  )_{t \in T}$. By a $T$-family of irreducible configurations on $M$ we mean a $T$-family of spin-c structures $\underline{\s}$ together with a continuosly varying family of (smooth) irreducible configurations $(A_t , \Phi_t ) \in \mathcal{C}^\ast (M, \s_t )$. Similarly, we adopt the notation $( \underline{A}, \underline{\Phi} )$ for such a family.

\begin{Remark} We will need to work with smoothly-varying families later on (with $T$ a smooth manifold); however, the discussion that follows applies equally well with only minor modifications. 
\end{Remark}

There is an obvious notion of "isomorphism" for two families of spin-c structures (resp. irreducible configurations) parametrised by the same space $T$: a continuosly varying $T$-family of isomorphisms of spin-c structures (resp. carrying the irreducible configurations onto one another). Much as before, the set of isomorphism classes of $T$-families of spin-c structures on $M$ is a torsor over the cohomology group $H^2 ( M \times T ; \Z )$. When it comes to families of irreducible configurations, the relevant "moduli functor" is represented by the irreducible configuration space:

\begin{Lemma}\label{correspondence}
There is a one-to-one correspondence 
between
\begin{enumerate}
    \item[(i)] the set of isomorphism classes of $T$-families of irreducible configurations $(\underline{A}, \underline{\Phi} )$ on $M$ with underlying isomorphism class of spin-c structure on $M$ represented by $\s_M$, and
    \item[(ii)] the set of continuous maps $\mathrm{Map}\big( T , \mathcal{B}^{\ast}(M, \s_M ) \big)$.
\end{enumerate}
\end{Lemma}

The main point is that $\mathcal{B}^\ast (M, \s_M )$ parametrises a \textit{universal} family $(\underline{A_\infty} , \underline{\Phi_\infty } )$ of irreducible configurations on $M$. This is constructed as follows. Say $\s_M = (g, S , \rho )$. The pullback of $S$ over the product $M \times \mathcal{C}^{\ast}(M, \s_M)$ is a $\mathcal{G}$-equivariant unitary vector bundle: the action of $v \in \mathcal{G}$ on the fibres of $S$ over $\{ m\} \times\mathcal{C}^{\ast}(M , \s_M)$ is given by multiplication by $v(m) \in \mathrm{U}(1)$; and the action on the base is the natural action on the second factor. The $\mathcal{G}$-action on the base is free, and passing to the quotient we obtain a unitary vector bundle $\underline{S_{\infty}}$ over $ M \times \mathcal{B}^{\ast}(M, \s_M )$ with a $\mathcal{B}^\ast (M, \s_M )$-family of Clifford multiplications. This yields a family of spin-c structures on $M$ parametrised by $\mathcal{B}^{\ast}(M , \s_M)$. Furthermore the tautological family of irreducible configurations on $M$ parametrised by $\mathcal{C}^{\ast}(M, \s_M)$ descends to a corresponding family of irreducible configurations $(\underline{A_\infty}, \underline{\Phi_\infty} )$ parametrised by $\mathcal{B}^\ast (M, \s_M )$.

Conversely, given a family of irreducible configurations $(\underline{A} , \underline{\Phi} )$ we construct an associated \textit{classifying map}
\begin{align}
f_{\underline{A} , \underline{\Phi}} : T \rightarrow \mathcal{B}^{\ast}(M, \s_M )\label{classifyingmap}
\end{align}
as follows. For a given $t \in T$ we choose an isomorphism $ \s_t \xrightarrow{\cong} \s_M$. Using this, we carry the irreducible configuration $(A_t , \Phi_t ) \in \mathcal{C}^{\ast}(M, \s_t )$ to an irreducible configuration $(A_{t}^\prime, \Phi_{t}^\prime ) \in\mathcal{C}^{\ast}(M, \s_M )$. This is done as follows. If $h: S_t \cong S$ is the underlying unitary isomorphism of spinor bundles, then we define $\Phi_{t}^\prime := h (\Phi_t )$; then $A_{t}^\prime$ is defined as the unique spin-c connection whose induced connection $\hat{A}_{t}^{\prime}$ on $\Lambda^2 S^+$ agrees with the connection $\hat{A}_{t}$ on $\Lambda^{2} S^{+}_t$ under the isomorphism $\Lambda^{2}S^{+}_t \cong \Lambda^{2} S^{+}$ induced by $h$. Note that choosing a different isomorphism $\s_t \cong \s_M$ only results in a gauge-equivalent irreducible configuration in $\mathcal{C}^{\ast}(M, \s_M)$. We then set $f_{\underline{A}, \underline{\Phi} }(t) = [( A_{t}^\prime , \Phi_{t}^{\prime})]$, which is easily verified to give a continuous map as we vary $t \in T$.

\begin{proof}[Proof of Lemma \ref{correspondence}] 
To go from (i) to (ii) we send a $T$-family of irreducible configurations $(\underline{A}, \underline{\Phi} )$ to its classifying map $f_{(\underline{A}, \underline{\Phi} )}$. For the other direction, if $f : T \rightarrow \mathcal{B}^{\ast}(M, \s_M )$ is a continuous map then the universal family of irreducible configurations $(\underline{A_{\infty} } , \underline{\Phi_{\infty}})$ parametrised by $\mathcal{B}^{\ast}(M, \s_M)$ can be pulled back along $f$ to produce a $T$-family of irreducible configurations on $M$.

The two assignments described above are inverse to each other. Indeed, given a family of irreducible configurations $(\underline{A}, \underline{\Phi})$ there is a \textit{unique} isomorphism of $(\underline{A}, \underline{\Phi} )$ with the pullback of the unversal family $(\underline{A_\infty} , \underline{\Phi}_\infty )$ by the map $f_{\underline{A}, \underline{\Phi}}$. The uniqueness follows again from the fact that $\mathcal{G}$ acts freely on irreducible configurations.
\end{proof}

The elementary correspondence from Lemma \ref{correspondence} implies the following "slogan" which plays a role in the upcoming construction of the families contact invariant:

\begin{Slogan}
One can trade a (possibly non-trivial) $T$-family of spin-c structures on $M$ carrying a $T$-family of irreducible configurations for a \textbf{constant} family of spin-c structures on $M$ together with a $T$-family of irreducible configurations which are only \textbf{well-defined up to $\mathcal{G}$-action}.
\end{Slogan}

In concrete terms, what this means is the following. Fix an open cover $T = \bigcup_{i \in I} U_i$ by contractible open sets. Then there is a correspondence between:
\begin{enumerate}
    \item[(i)] the set of isomorphism classes of $T$-families of irreducible configurations $(\underline{A}, \underline{\Phi} )$ on $M$ with underlying isomorphism class of spin-c structure on $M$ represented by $\s_M$, and   
    \item[(ii)] isomorphism classes of $I$-tuples of continuous maps $\big( (A_i , \Phi_i ) : U_i \rightarrow \mathcal{C}^\ast (M , \s_M )\big)_{i \in I} $ such that for each overlap $U_i \cap U_j$ there exists a (unique) continuous map $v_{ji} : U_i \cap U_j \rightarrow \mathcal{G}$ such that $v_{ji}(t)\cdot ( A_i (t) , \Phi_i (t) ) = (A_j (t), \Phi_j (t) )$.
    \end{enumerate}

Let us mention at this point that the role played by certain families of irreducible configurations (coming from families of contact structures) is going to be to provide natural "boundary conditions" for the Seiberg--Witten equations over an end of a non-compact $4$-manifold. It will be necessary later to "trivialise" the family of spin-c structures, and the $\mathcal{G}$-ambiguity of the resulting family of irreducible configurations will pose no issue due to the $\mathcal{G}$-invariance of the Seiberg--Witten equation.

\subsection{Families of irreducible configurations from symplectic and contact structures}

\subsubsection{Symplectic $4$-manifolds}\label{spincsympl}

Let $(X , \omega)$ be a symplectic $4$-manifold, oriented by the volume form $\omega^2$. We make the auxiliary choice of an $\omega$-compatible almost-complex structure $J$. This means that the tensor $g = \omega ( \cdot , J \cdot )$ defines a Riemannian metric. We refer to such a triple $(\omega , J , g )$ as an \textit{almost-K\"{a}hler} structure on $X$. 

\begin{Definition} The \textit{canonical spin-c structure} $\s_{\omega , J , g } = ( g , S_{\omega , J , g} , \rho_{\omega , J , g} )$ determined from the almost-K\"{a}hler structure $(\omega , J , g)$ is given by the following data:
\begin{itemize}
\item $S_{\omega , J , g}^+ = \C \oplus \Lambda^{0,2}_{J}T^{\ast} X$ and $S_{\omega , J , g}^- = \Lambda^{0,1}_{J}T^{\ast} X$,  equipped with the hermitian metrics naturally induced from $g$. 

\item the Clifford multiplication by $\eta \in T^{\ast} X$ has the component $\rho_{\omega , J , g}^+ (\eta): S_{\omega , J , g}^+ \rightarrow S_{\omega , J , g}^-$ defined by
\begin{align*}
\rho_{\omega , J , g}^+ (\eta ) (\alpha , \beta )  = \sqrt{2} ( \eta^{0,1} \wedge \alpha - \iota_{\eta_{0,1}} \beta ) .
\end{align*}
\end{itemize}
\end{Definition}

Above, $\iota_{X}$ stands for contraction by $X$ on the first component, and $\eta_{0,1}$ is the $(0,1)$-part of the metric dual tangent vector of $\eta$. The remaining component of the Clifford action, $\rho_{\omega , J , g}^- : S_{\omega , J , g}^- \rightarrow S_{\omega , J , g}^+$ can be recovered from the above, using the fact that $\rho_{\omega , J , g}$ should be skew-adjoint.


A computation shows that the Clifford action of the symplectic form $ \rho_{\omega , J , g}(\omega) : S_{\omega , J , g}^+ \rightarrow S_{\omega , J , g}^+$ is given by $-2i$ on $\C$ and $+ 2i$ on $\Lambda^{0,2}_{J}$. Observe that there is a canonical section $\Phi_{\omega , J , g}$ of $S_{\omega , J , g}^+$ given by constant $1$ on the $\C$ component. 

\begin{Lemma}\cite{taubsymp} \label{canonicalconnection}
There exists a unique spin-c connection $A_{\omega , J , g}$ on $S_{\omega , J , g}$ such that $$D_{A_{\omega , J , g}}^{+} \Phi_{\omega , J , g} = 0.$$
\end{Lemma}
\begin{Remark}
Alternatively, $A_{\omega , J , g}$ is uniquely characterised by the property that the covariant derivative $\nabla_{A_{\omega , J , g}} \Phi_{\omega , J , g}$ is a 1-form with values in the subbundle $\Lambda^{0,2}_{J}T^{\ast} X$.
\end{Remark}

\begin{Definition}
The \textit{canonical configuration} associated to $(\omega, J , g)$ is the pair $(A_{\omega , J , g} ,\Phi_{\omega , J , g} ) \in \mathcal{C}^{\ast}(X , \s_{\omega , J , g} )$.
\end{Definition}

Thus, the space of almost-Kähler structures on $X$ parametrises a family of irreducible configurations on $X$.

\begin{Remark}It is a fundamental Theorem of Taubes \cite{taubsymp} that the Seiberg--Witten invariant $\mathrm{SW}(\s_{\omega, J , g} ) \in \Z$ of the canonical spin-c structure of a closed symplectic $4$-manifold with $b^{+}(X) > 1$ is non-vanishing. Taubes' proof shows that, for a suitable large perturbation of the Seiberg--Witten equations, the canonical configuration becomes the only solution to the equations, modulo $\mathcal{G}$-action.
\end{Remark}


\subsubsection{Contact $3$-manifolds}\label{spinccont}

Let $(Y , \xi )$ be a contact $3$-manifold. We now choose the auxiliary data of a complex structure $j$ on the contact distribution (inducing the positive orientation) and a (positive) contact form $\alpha$. We will refer to such a triple as \textit{contact metric structure}. Indeed, given $(\xi , \alpha , j )$ there exists a \textit{unique} Riemannian metric $g_{\xi , \alpha ,j}$ on $Y$ characterised by
\begin{itemize}
\item $ |\alpha|_{g_{\xi, \alpha , j}} = 1$
\item $d \alpha  = 2 \ast \alpha$ where $\ast$ is the Hodge star operator of $g_{\xi , \alpha , j}$
\item $j $ is an isometry of $(\xi , g_{\xi , \alpha , j} )$.
\end{itemize}

Observe that the \textit{Reeb vector field} $R$ (determined uniquely by the requirement that $\alpha (R ) = 1$ and $d\alpha (R , \cdot ) = 0$) is $g_{\xi , \alpha , j}$-orthogonal to the contact plane $\xi$. It is convenient to regard $j$ as an endomorphism of $TY$ by setting $j (R)= 0$. Then, we can write down explicitly the Riemannian metric $g_{\alpha , j}$ as
\begin{align}
g_{\xi , \alpha , j} = \alpha \otimes \alpha + \frac{1}{2} d \alpha (\cdot , j \cdot ). \label{contactmetric}
\end{align}
\begin{Definition} The \textit{canonical spin-c structure} $\s_{\xi , \alpha , j} = ( g_{\xi , \alpha , j} , S_{\xi , \alpha , j} , \rho_{\xi , \alpha , j} )$ determined from the contact structure $\xi$ and the auxiliary data $\alpha , j$ is given by the following data
\begin{itemize}
\item $S_{\xi , \alpha , j } = \C \oplus \langle \alpha \rangle^{\perp}$ where the second factor is the $g_{\xi , \alpha , j}$-orthogonal complement of $\alpha$ inside of $T^{\ast} Y$ 
\item $\rho_{\xi , \alpha , j} (\eta ) (x , y )  = (i \eta (R ) x , - i \eta (R ) y ) - \sqrt{2} (\eta^{0,1}x - \iota_{\eta_{0,1}} y )$ for $\eta \in T^\ast Y$.
\end{itemize}
\end{Definition}
For the above, note that we can decompose $\eta \in T^{\ast} Y$ as $\eta = \eta (R) \alpha + \eta^{1,0} + \eta^{0,1}$ where $\eta^{p,q} \in \langle \alpha \rangle^{\perp} \otimes_{\R} \C$ stands for the $(p,q)$ component of the projection to $\langle \alpha \rangle^{\perp}$ of $\eta$, using the complex structure $j$ on $\langle \alpha \rangle^{\perp}$.

The $3$-dimensional contact analogue of Taubes' Theorem about closed symplectic $4$-manifolds now states that for a contact structure $\xi$ which admits a strong symplectic filling, the contact invariant $\mathbf{c}(\xi ) \in \widecheck{HM}^{\ast}(-Y , -\s_{\xi, \alpha , j} )$, is non-vanishing. A monopole Floer proof of this result has recently been given by Echeverria \cite{mariano}.

Given a contact form $\alpha$ for $(Y, \xi)$, by its \textit{symplectization} we will mean the symplectic manifold (with concave boundary) $(K , \omega )$ where $K = [1, + \infty ) \times Y$ and
\begin{align}
\omega = d (\frac{s^2}{2} \alpha ) = s ds \wedge \alpha + \frac{s^2}{2} d \alpha .\label{symplectization}
\end{align}
If we start with a triple $(\xi , \alpha , j)$ on $Y$, we obtain an almost-K\"{a}hler structure $(\omega , J , g )$ on $K$ by having $J$ agree with $j$ on $\xi = \ker \alpha$ and setting 
\begin{align}
J ( \partial / \partial s ) = \frac{1}{s} R \label{almostcomplex} 
\end{align} where $R$ is the Reeb vector field of $\alpha$. It follows that the Riemannian metric $g = \omega( \cdot , J \cdot)$ over $K = [1, + \infty) \times Y$ is the \textit{cone metric over} $(Y , g_{\xi , \alpha , j } )$, namely 
\begin{align*}
g = ds^2 + s^2 g_{\xi , \alpha , j} . 
\end{align*}

\begin{Lemma}[\cite{marianothesis}, Lemma 35]\label{restriction} There is a canonical identification of spin-c structures on $-Y = \partial K$ $$(\mathfrak{s}_{\omega , J , g })|_{ - Y} \cong  -\mathfrak{s}_{\xi , \alpha, j }.$$
\end{Lemma}

Above we denote by $- \mathfrak{s}_{\xi , \alpha , j}$ the induced spin-c structure on $-Y$ (obtained by adding a negative sign to the Clifford multiplication). 

\begin{Definition}
The \textit{canonical configuration} associated to $(\xi, \alpha , j )$ is the pair $(A_{\xi, \alpha, j }, \Phi_{\xi, \alpha , j} ) \in \mathcal{C}^{\ast}(Y, \s_{\xi, \alpha  , j} )$ obtained by restriction onto $Y$ of the canonical configurations $(A_{\omega , J , g} ,\Phi_{\omega , J , g} ) \in \mathcal{C}^{\ast}(K , \s_{\omega , J , g} )$ associated to the almost-K\"ahler structure $(\omega , J , g )$ on $K= [1, + \infty) \times Y$.
\end{Definition}

Thus, the space of contact metric structures parametrises a family of irreducible configurations on $Y$ and $K = [1, + \infty) \times Y$.

\section{Construction of the families contact invariant}\label{invtsection}

In this section we construct the families contact invariant (\ref{invtmap}). There is a Poincar\'{e} duality for the Floer groups [\cite{KM}, \S 3], $\widecheck{HM}_{\ast} (- Y ,\s_{\xi_0} ) \cong \widehat{HM}^{\ast}( Y , \s_{\xi_0} )$, and the map (\ref{invtmap}) most naturally arises as a map into the latter group: the \textit{from} version of the monopole Floer cohomology groups. We give a rough outline of this construction before going into the details.

First, we equip contact structures with auxiliary structures. Let $\mathcal{CM}(Y,\xi_0 )$ be the space of \textit{contact metric structures} (see \S \ref{spinccont}) $(\xi, \alpha , j )$ on $Y$ such that the contact structures $\xi$ and $\xi_0$ are isotopic. The forgetful projection induces a weak homotopy equivalence $\mathcal{CM}(Y , \xi_0 ) \simeq \mathcal{C}(Y , \xi_0 )$. We will also find it convenient to work within the realm of Banach spaces. A way to do this is by considering triples $(\xi , \alpha , j )$ where $\alpha$ (and hence $\xi$) is only assumed to be of class $C^{l}$, and the complex structure $j$ is of class $C^{l-1}$, for a suitable positive integer $l$. The metric $g_{\xi , \alpha , j}$ (see (\ref{contactmetric})) determined from the triple $(\xi, \alpha , j)$ is therefore of class $C^{l-1}$. The space of such triples is a Banach manifold homotopy equivalent to the space of $C^\infty$ triples. From now on, we will reserve the notation $\mathcal{CM}( Y , \xi_0 )$ for this more convenient Banach manifold version only.

Associated to each triple in $\mathcal{CM}(Y, \xi_0 )$ and each element of a certain Banach space of perturbations $\mathcal{P}$ we consider the Seiberg--Witten monopole equations over a certain non-compact $4$-manifold $Z^+$, with suitable asymptotics over its ends to canonical configurations determined by the contact structures together with critical points of the Chern-Simons functional. This leads to a \textit{universal} moduli space of solutions, which is a Banach manifold equipped with a Fredholm map
\begin{align*}
\mathfrak{M}(Z^+ ) \xrightarrow{\pi} \mathcal{CM}(Y , \xi_0 ) \times \mathcal{P} .
\end{align*}
 
 The moduli space decomposes according to critical points of the Chern-Simons-Dirac functional $$ \mathfrak{M}(Z^+ ) = \bigcup_{[\afr]} \mathfrak{M}([\afr], Z^+).$$

Given a generic cycle $T$ in $\mathcal{CM}(K , \xi_0 ) \times \mathcal{P}$ transverse to the Fredholm map $\pi$, we count \textit{isolated} points in $\mathfrak{M}(Z^+)$ which lie over $T$, and this leads to integers $\# ( \mathfrak{M}([\afr] , Z^+ ) \cdot T ) \in \Z$. Indexing the counts by the critical points $[\afr]$ we obtain a cocycle in the Floer cochain complex$$\psi ( T ) = \sum_{[\afr]} \# \big( \mathfrak{M}([\afr] , Z^+ ) \cdot T \big) \cdot [\afr]  \in \widehat{C}^{\ast}(Y , \s_{\xi_0}; R ).$$ This yields the homomorphism (\ref{invtmap}). In fact, we will be able to define the homomorphism at the chain level.



\subsection{Differential-geometric aspects}

\subsubsection{The symplectic end and the cylindrical end}\label{AKend}

We start by discussing the various metric structures that come into the construction.

\begin{Remark} For ease in notation we will denote elements of $\mathcal{CM}(Y,\xi_0 )$ by the symbol $t$. When we need to make reference to it, the contact metric structure on $Y$ associated to $t$ is denoted $(\xi_t , \alpha_t , j_t )$. From now on, we also fix a $C^{\infty}$ base triple $(\xi_0 , \alpha_0 , j_0 ) \in \mathcal{CM}(Y , \xi_0 )$.
\end{Remark}

 Let $Z^+$ be the non-compact $4$-manifold $\R \times Y$ with the product orientation. Let $K = [1, + \infty) \times Y \subset Z^+$ and $Z = (- \infty , 0]  \times Y $. For each $t \in \mathcal{CM}(Y , \xi_0 )$ we have an almost-K\"{a}hler structure $( \omega_t , J_t , g_t )$ over $K$ obtained from (\ref{symplectization}-\ref{almostcomplex}), where $\omega_t$ is $C^{l-1}$, $J_t$ is $C^{l}$ and $g_t$ is $C^{l-1}$. Recall that $g_t$ is the cone metric over $(Y , g_{\xi_t, \alpha_t , j_t})$, namely $g_t = ds^2+ s^2 g_{\xi_t, \alpha_t , j_t}.$ 

We now extend the metric $g_t$ from $K$ to the whole of $Z^+$. Over $Z = ( - \infty , 0 ] \times Y$ the metric $g_t$ agrees with the cylindrical product metric $ds^2 + g_{\xi_0 , \alpha_0 , j_0}.$ We fix the behaviour of the metric $g$ over the region $[0, 1] \times Y $ as follows. Choose a smooth function $\kappa  : [0 , 1] \rightarrow \R_{\geq 0}$ such that $\kappa \equiv 1$ on a neighbourhood of $[0 , 1/2]$ and $\kappa \equiv 0$ on a neighbourhood of $1$. Then the metric $g_t$ over the region $[0 , 1] \times Y$ is defined as $
ds^2 + \kappa (t) g_{\xi_0 , \alpha_0 , j_0} + (1-\kappa(s))s^2 g_{\xi_t , \alpha_t , j_t}.$

We will refer to $K$ as the \textit{conical} or \textit{symplectic end} of $Z^+$, and to $Z$ as the \textit{cylindrical end}. Observe that in this construction the family of metrics $g_t$ restricted over the cylindrical end $Z$ is independent of $t$ (it only depends on the fixed base triple $(\xi_0 , \alpha_0 , j_0 )$). We will denote by $(g_0 , \omega_0 , J_0 )$ the corresponding structures determined by the base triple $(\xi_0 , \alpha_0 , j_0 )$.

\subsubsection{Families of spin-c structures and canonical configurations}\label{confscon}

We move on now to discuss families of spin-c structures and irreducible configurations on the non-compact manifold $Z^+$. The latter will provide us with the right boundary conditions for the Seiberg--Witten equations over the symplectic end later on.

We consider the "trivial" family of spin-c structures on $Z^+$ parametrised by $\mathcal{CM}(Y, \xi_0 )$. Namely, we start with the spin-c structure $\s_{\omega_0 , J_0 , g_0}$ on $K$ determined by the almost-Kähler structure $(\omega_0 , J_0 , g_0 )$ (see \S \ref{spincsympl}). By Lemma \ref{restriction} we have $\s_{\omega_0 , J_0 , g_0}|_{-Y} = -\mathfrak{s}_{\xi_0 , \alpha_0 , j_0 }$. Hence, we may extend the spin-c structure $\s_{\omega_0 , J_0 , g_0} $ from $K$ over to the whole of $Z^+$ in a translation-invariant manner over $Z^+ \setminus K$. We denote this spin-c structure on $Z^+$ by $\mathfrak{r}_{0} = (g_0 , S , \rho )$. Finally, by \textit{changing metrics} we obtain a family of spin-c structures on $Z^+$ parametrised by $\mathcal{CM}(Y, \xi_0 )$: for $t \in \mathcal{CM}(Y, \xi_0 )$ we set 
$ \mathfrak{r}_{t } = (g_t , S, \rho_t )$ with $\rho_t = \rho \circ b_{g_t , g_0}^\ast $. In particular, observe that the spinor bundle over $Z^+$ of $\mathfrak{r}_t$ is always $S$, independently of $t$, and that the metrics and spin-c structures on the cylinder end $Z$ obtained by restriction of $\mathfrak{r}_t$ are independent of $t$.

We now discuss irreducible configurations on $Z^+$ parametrised by $\mathcal{CM}(Y, \xi_0 )$. The Clifford action of the symplectic structures gives a trace-less skew-adjoint map $\rho_t (\omega_t ) :  S^+|_K \rightarrow S^+|_K  $ such that $\rho (\omega_t )^2 = -4 \cdot \mathrm{id}_{S^+|_{K}}$. This induces a decomposition 
\begin{align*}
    S^+|_K = E_{- } (t) \oplus E_{+} (t )
\end{align*} 
into $\mp 2i$ eigensubbundles (each with with rank $1$). Because $\mathfrak{r}_t$ is (non-canonically!) isomorphic over $K$ to the spin-c structure induced from $(\omega_t , J_t , g_t )$ it follows that the $- 2i$ eigensubbundle admits a trivialisation $E_- (t) \approx \mathbb{C}$ for each $t$. We warn the reader that, as a bundle over $K \times \mathcal{CM}(Y, \xi_0 )$, the bundle $E_-$ need \textit{not} admit a trivialisation.

Let $U \subset \mathcal{CM}(Y, \xi_0 )$ be an open contractible subset. Then we may choose a unitary trivialisation of $E_-$ over $K \times U$ and obtain a $U$-family of nowhere-vanishing spinors $\Phi_t \in \Gamma (K , E_- )$ with pointwise unit length. As in \S \ref{spincsympl} there is a unique $U$-family of spin-c connections $A_t$ over $K$ such that $D_{A_t}^+ \Phi_t = 0$. We refer to the $U$-family of irreducible configurations $(A_t, \Phi_t )$ as the \textit{canonical configurations} over the symplectic end (associated to a given trivialisation of $E_-$ over $K \times U$).

\begin{Remark}
By putting together the canonical configurations $(A_t , \Phi_t )$ over all $U \subset \mathcal{CM}(Y,\xi_0 )$ (and applying a further change of metrics back to $g_0$) we obtain a continuous map $$ f : \mathcal{CM}(Y, \xi_0 ) \rightarrow \mathcal{B}^\ast (K , \mathfrak{s}_{\omega_0 , J_0 , g_0 } ) .$$ The interpretation of this map is clear: the space $\mathcal{CM}(Y, \xi_0 )$ parametrises a family of almost-Kähler structures on $K$, which in turn parametrises a family of irreducible configurations over $K$ as in \S \ref{spincsympl}. The classifying map for this family (in the sense of Lemma \ref{correspondence}) is the map $f$. At a heuristic level, the families contact invariant that we construct in this section should be regarded as a sort of "pushforward map in homology" induced by $f$.
\end{Remark}

\subsubsection{Translation invariance of the canonical configurations}\label{translationinvariance}
 
We now discuss how the canonical configurations over the symplectic end can be made \textit{translation invariant} in a suitable sense, which will become convenient occasionally.

Recall that $g_t = ds^2 + s^2 g_{\xi_t , \alpha_t , j_t }$ is a conical metric over $K = [1,+\infty ) \times Y$. Denote by $\overline{g}_t = ds^2 + g_{\xi_t , \alpha_t , j_t }$ the corresponding cylindrical metric. The \textit{rescaling operator} $\mathcal{R}_t : (TK , \overline{g}_t )  \rightarrow (TK , g_t )$ gives an isometry between the two metrics:
\begin{align*}
    \mathcal{R}_t (\partial_s ) &= \partial_s\\
    \mathcal{R}_t ( v ) &= \frac{1}{s} v \,\, , \,\, v \in TY .
\end{align*}
The almost complex structure $J_t$ is carried to a translation-invariant almost complex structure $\overline{J}_t = \mathcal{R}_{t}^{-1} \circ J_t \circ \mathcal{R}_t$, and we have corresponding unitary vector bundle isometries
\begin{align*}
    \mathcal{R}_{t}^\ast : ( \Lambda_{J_t}^{p,q} T^\ast K , h_t ) \rightarrow ( \Lambda_{\overline{J}_t}^{p,q} T^\ast K , \overline{h}_t )
\end{align*}
where $h_t$ and $\overline{h}_t$ stand for the hermitian metrics determined by the pairs $(g_t , J_t )$ and $(\overline{g}_t , \overline{J}_t )$.

At this point, we recall that the spinor bundle $S = S^+ \oplus S^-$ underlying our family of spin-c structures $\mathfrak{r}_t$ is independent of $t$ and has the following form over $K$ (see \S \ref{confscon} and \S \ref{spincsympl}):
$$ S^+ = \mathbb{C} \oplus \Lambda^{0,2}_{J_0}T^\ast K \,\, , \,\, S^- = \Lambda^{1,1}_{J_0} T^\ast K .
$$
The "rescaled" unitary bundle $\overline{S^+} := \mathbb{C} \oplus \Lambda^{0,2}_{\overline{J}_0}T^\ast K$ (and likewise for $\overline{S^-} $) is identified with the pullback to $K = [1, +\infty )\times Y$ of a bundle over $Y$, and hence one can speak of translation-invariant sections or connections on this bundle. We have:

\begin{Lemma}\label{spinorasymp}
Let $(A_{t} , \Phi_{t})$ be a $U$-family of canonical configurations on $K $, for a contractible open $U \subset \mathcal{CM}(Y, \xi_0 )$. After applying a $U$-family of smooth gauge transformations $g_t : K \rightarrow \mathrm{U}(1)$ to $(A_t , \Phi_t )$ we may assume that:
\begin{enumerate}[(i)]
\item the sections $ \overline{\Phi}_t := \mathcal{R}_{0}^{\ast} \Phi_{t} \in  \Gamma ( K ,\overline{S^+} )$ are translation-invariant, and
\item the connections $\overline{A}_t := \mathcal{R}^{\ast}_{0} A_t$ on $\overline{S} = \overline{S^+}\oplus \overline{S^-}$ are translation-invariant.
\end{enumerate}
\end{Lemma}

\begin{Definition}\label{transinvform}
A family of canonical configurations $(A_{t} , \Phi_{t} )$ over $K$ parametrised by $U$ is in \textit{translation-invariant form} if it satisfies (i)-(ii) above.
\end{Definition}

\begin{proof}[Proof of Lemma \ref{spinorasymp}]
The family of spin-c structures $\mathfrak{r}_t= (g_t , S, \rho_t )$ is isomorphic, via the rescaling operator $\mathcal{R}^{\ast}$, to the family of spin-c structures $\overline{\mathfrak{r}}_t = (\overline{g}_t , \overline{S} , \overline{\rho}_t )$ where 
\begin{align}
\overline{\rho}_t (v ) = \mathcal{R}_{0}^{\ast} \rho_t ((\mathcal{R}_{t}^{\ast})^{-1} v) (\mathcal{R}_{0}^{\ast})^{-1} , \,\, v \in T^\ast K . \label{transcliff}
\end{align}
Now, we have a translation-invariant non-degenerate $2$-form $\overline{\omega}_t := \mathcal{R}^{\ast} \omega_t =  ds \wedge \alpha_t + \frac{1}{2} d \alpha_t$, so the $-2i$-eigensubbundle $\overline{E}_- \subset \overline{S^+}$ corresponding to the action of $\overline{\rho}_u (\overline{\omega}_u )$ on $\overline{S^+}$ is also translation-invariant. Thus, in view of (\ref{transcliff}) the assertion $(i)$ is clear: one chooses a trivialisation of $\overline{E}_-$ over $K \times U$ by a translation-invariant unit section $\overline{e} \in \Gamma (K \times U , \overline{E}_- )$ to obtain a $U$-family of sections $ (\mathcal{R}_{t}^{\ast})^{-1} \overline{e}(\cdot,t ) \in \Gamma (K , E_{-}(t) )$ which agree with $\Phi_t$ up to gauge transformations. Next, we verify that $(ii)$ holds assuming that $\Phi_{t}$ satisfies $(i)$. Write the covariant derivative with respect to $(\mathcal{R}_{0})^\ast A_t$ as $\frac{d}{ds} + \nabla_{B_t (s)} + c_t (s) ds$, where for each $t$ we have paths of connections $s \mapsto B_{t}(s)$ and $i \mathbb{R}$-valued functions $s\mapsto c_t (s)$ on $Y$ (parametrised by $s \in [1, + \infty )$). Recall that $A_t$ can be characterised by the property that $\nabla_{A_t}\Phi_t $ is orthogonal to $\Phi_t$ (see Remark after Lemma \ref{canonicalconnection}) we obtain $$ \langle \nabla_{B_t (s)} \overline{\Phi}_{t} , \overline{\Phi}_t \rangle_{\overline{h_0}} + c_t (s) ds \equiv  0 .$$ It follows that both terms above vanish. The vanishing of $\langle \nabla_{B_t (s)} \overline{\Phi}_{t} , \overline{\Phi}_t \rangle_{\overline{h_0}}$ also gives us that $B_t (s)$ is independent of $s$ since $\overline{\Phi}_t$ is translation-invariant. Thus, (ii) follows.
\end{proof}

\subsubsection{A basic example}

It may be instructive to review the above constructions in a particularly simple case. 

Consider the flat hyperk\"{a}hler structure $(g_0 , J_1 , J_2 , J_3 )$ on $\R^4$. The radial vector field $v = x \partial_x + y \partial_y + z \partial_z + w \partial_w $ in $\mathbb{R}^4$ is Liouville for all symplectic structures in the family $\omega_t = \sum_{i = 1}^{3} t_i g_0 ( J_i \cdot , \cdot )$ parametrised by $t \in S^2 \subset \mathbb{R}^3$ (i.e. $\mathcal{L}_{v} \omega_u = \omega_u $) and $v$ is transverse to $S^3 \subset \mathbb{R}^4$. Thus there is an $S^2$-family of contact forms $\alpha_t$ on $S^3$ given by $\alpha_t = \iota_{V} \omega_t$. This family of contact structures is $\mathrm{SU}(2)$-invariant, and thus will descend to a family of contact structures on the quotients $S^3 / \Gamma$ by a finite subgroup $\Gamma \subset \mathrm{SU}(2)$. The manifolds $S^3 / \Gamma$ are precisely the links of the ADE singularities (which include e.g. the lens spaces $L(p,p-1)$ or the Poincaré sphere $\Sigma (2,3,5 )$ ). 

Let $Y = S^3 / \Gamma$. A complex structure $j_t$ on the contact distribution $\xi_t = \mathrm{ker} \alpha_t$ is obtained by restricting $J_t = \sum_i t_i J_i$. We thus have a family $(\xi_t , \alpha_t , j_t ) \in \mathcal{CM}(S^3 , \xi_0 ) $, and we take the base triple $(\xi_0 , \alpha_0 , j_0 ):= (\xi_{t}, \alpha_t , j_t )|_{t = (1,0,0)} $. The associated family of almost-Kähler structures on $K = [1,+\infty) \times Y$ agrees with the flat hyperkähler structure under the identification $$K \cong (\mathbb{R}^4 \setminus \{ x^2+y^2+z^2+w^2 < 1\} )/ \Gamma \, \,, \,\, (s,(x,y,z,w)) \mapsto (sx,sy,sz,sw).$$

We next calculate canonical configurations. The positive spinor bundle is $S^+ = \mathbb{C} \oplus \Lambda_{J_1}^{0,2} T^\ast K \cong \mathbb{C}^2 $, and we have trivialised $\Lambda_{J_1}^{0,2} T^\ast K$ using $\frac{1}{2}d\overline{z}^1  \wedge d \overline{ z}^2$. Likewise $S^- \cong \mathbb{C}^2$. The Clifford multiplications $\rho_t = \rho$ are independent of $t$. The symplectic forms $\omega_i = g_0 (J_i , )$ have the following Clifford actions on $S^+ = \mathbb{C}^2$:
\begin{align*}
\rho(\omega_1) = \begin{pmatrix} -2i & 0 \\ 0 & 2i \end{pmatrix} \quad \rho(\omega_2) = \begin{pmatrix} 0 & -2 \\ 2 & 0 \end{pmatrix} \quad \rho(\omega_3) = \begin{pmatrix} 0 & 2i \\ 2i & 0 \end{pmatrix}.
\end{align*}
and thus
\begin{align*}
    \rho (\omega_t ) = \begin{pmatrix} - 2i t_1 & -2 (t_2 - i t_3 )\\ 2(t_2 + i t_3) & 2it_1 \end{pmatrix}.
\end{align*}
An $S^2$-family of sections of $-2i$-eigenspace of $\rho( \omega_t )$ is $( 1+t_1 ,  it_2-t_3).$ Each has a transverse zero at $t = (-1,0,0)$ and is non-vanishing elsewhere. This means that $E_-$ is not trivial over $K \times S^2$. Normalising we obtain a family of unit length sections of $E_- (t)$ over $K$ for $t \in U_1 := S^2 \setminus (-1,0,0)$:
$$\Phi_t = \big( \sqrt{ \frac{1+t_1}{2} }, \frac{it_2 - t_3}{\sqrt{2(1+t_1)}} \big) \in S^+ = \mathbb{C}^2 .$$ The corresponding family of spin-c connections $A_t$ is independent of $t$ and is given by the trivial connection on $S =S^+ \oplus S^- = \mathbb{C}^2 \oplus \mathbb{C}^2$. The family of canonical configurations carried by $U_1$ that we just constructed is also in translation-invariant form.


\subsection{Space of configurations}\label{confsection}

We now construct a suitable space of configurations $(A, \Phi )$ over the symplectic end $K$ which has the structure of a Banach manifold.

\subsubsection{Sobolev spaces on non-compact manifolds}\label{Sobolev}

To work in the convenient setting of Fredholm theory we make use of Sobolev spaces over the non-compact symplectic end. On a Riemannian manifold $(M,g)$ of \textit{bounded geometry}, the various possible definitions of Sobolev spaces of sections will agree. We refer the reader to [\cite{Eichhorn}, Chapter 11], or to [\cite{marianothesis}, \S 3.2] for an exposition of these results. The cone over a closed Riemannian manifold, the case that concerns us, falls into this desirable category.

In the above setting, given an Euclidean vector bundle $E \rightarrow M$ with an orthogonal connection $A$, the space of Sobolev sections $L^{2}_{k , g , A}(M,E)$ can be defined as the space of measurable sections $s$ of $E$ with distributional derivatives up to order $k$ and such that $$||s||_{L^{2}_{k}}^{2} := \sum_{j \leq k} \int_{M} | \nabla_{A}^{(j)} s |_{h}^2 \mathrm{dvol}_g < + \infty .$$

In the above formula, $\nabla_{A}^{(j)}$ is the connection on $ E \otimes (T^\ast M )^{\otimes (j-1 )}$ induced from $A$ and the Levi-Civita connection of $g$, and the symbol $| \cdot |_h$ denotes the metric induced from $h$ and $g_t$ on the bundle $E \otimes (T^\ast K)^{\otimes j}$. The vector space $L_{k,g,A}^{2}(M,E)$ equipped with the $L^{2}_{k}$ inner product becomes a Hilbert space. When $g$, $E$, or $A$ are understood we might drop them from the notation. 

From now on, we will fix an integer $k \geq 4$, which ensures that $L^{2}_{k}$ configurations over a $4$-manifold of bounded geometry are in $C^0$ by the Sobolev embedding theorem. We recall that we have been working thus far with the space $\mathcal{CM}(Y , \xi_0 )$ of triples $(\xi, \alpha , j)$, where the regularity of $\xi$ and $\alpha$ is $C^l$, and $j$ is $C^{l-1}$. Hence $g_{\xi , \alpha , j}$ is $C^{l-1}$. We fix the the integer $l$ so that $l-k-2 \geq 2$, because we will later need that $C^{l-k-2} \subset C^{2}$.


\subsubsection{Boundary conditions over the symplectic end} \label{confspaces}
We now set up the relevant configuration spaces over the symplectic end, with asymptotics to the canonical configurations provided by the contact geometry. Because canonical configurations only exist over sufficiently small neighbourhoods $U \subset \mathcal{CM}(Y , \xi_0 )$, our construction of configuration spaces will involve taking suitable limits over such neighbourhoods. 

In what follows, it is convenient to consider the slightly larger region containing the symplectic end: $K^{\prime} = [0 , 1] \times Y \cup K  \subset Z^+ .$
Let $U \subset \mathcal{CM}(Y , \xi_0 )$ be an open contractible subset, carrying a family of canonical configurations $\gamma := ((A_{t} , \Phi_{t}))_{t \in U}$ defined over $K$ which are in translation-invariant form (Definition \ref{transinvform}).

\begin{Definition}\label{conf} 
For $(U, \gamma)$ as above, the \textit{configuration space} for $(U , \gamma )$, denoted $\mathcal{C}_{k}(K^{\prime} , \gamma)_U$, is the space of triples $(t, A , \Phi )$, where $t \in U$, $A$ is a locally $L^{2}_{k}$ spin-c connection on the spinor bundle $S$ (for the spin-c structure $\mathfrak{r}_t$) defined over $K^{\prime}$ and $\Phi$ is a locally $L^{2}_{k}$ section of $S^+$ over $K^\prime$, subject to the following asymptotics:
\begin{align}
\Phi - \Phi_{t} \in L^{2}_{k , g_t, A_{t}}(K , S^{+}) \label{spinorasy}\\
A - A_{t} \in L^{2}_{k, g_t}(K , T^{\ast} K \otimes i \R). \label{connasym}
\end{align}
\end{Definition}

The relevant \textit{gauge group} in this setting is the group $\mathcal{G}_{k+1}(K^{\prime} )$ of locally $L^{2}_{k+1}$ maps $v : K^{\prime} \rightarrow \mathrm{U}(1)$ which approach the identity, i.e. 
\begin{align*}
1-v \in L^{2}_{k+1, g_t }(K).
\end{align*}
Again, the Sobolev space above does not depend on $t$. Observe that configurations in $\mathcal{C}_{k}(K^{\prime} , \gamma )_U$ are necessarily \textit{irreducible} (i.e. $\Phi$ doesn't vanish everywhere on $K$) due to the asymptotic condition (\ref{spinorasy}). Hence $\mathcal{G}_{k+1}(K^{\prime})$ acts freely on $\mathcal{C}_{k}(K^{\prime} , \gamma_0 )_U$.

Since for any two conical metrics $g_0 , g_1$ over $K$ the difference $g_1 g_{0}^{-1}$ is bounded over $K$ and the configurations $(A_t , \Phi_t )$ were chosen in translation-invariant form (Definition \ref{transinvform}), it follows that the Sobolev spaces $L^{2}_{k, g_t , A_{t}} (K , S^+ )$ and $L^{2}_{k , g_t }(K , T^\ast K \otimes i \mathbb{R} )$ are \textit{independent} of $t \in U$. The configuration space for $(U , \gamma ) $ then forms a trivial bundle of affine Hilbert spaces $$\mathcal{C}_{k}(K^{\prime}, \gamma )_U \rightarrow U .$$ We make $\mathcal{C}_{k}(K^{\prime}, \gamma )_U$ into a Banach manifold by identifying it with $L^{2}_k \times U$ via $(A, \Phi , t ) \mapsto (A- A_t , \Phi - \Phi_t , t )$. In this "chart", the $\mathcal{G}_{k+1}(K^\prime )$-action $\mathcal{G}_{k+1}(K^\prime) \times (L^{2}_k \times U ) \rightarrow (L^{2}_k \times U )$ acquires the rather odd-looking form: $v \cdot (a , \phi , t ) =   (a - v^{-1} dv , v \phi - (1-v )\Phi_t , t )$. This action is only of class $C^{l-k-2}$. The reason is that $\Phi_t$ depends on first derivatives of the metric $g_t$ (and also on $\alpha_t$ and $j_t$) which has regularity $C^{l-1}$; thus we may only differentiate $l-k-2 = (l-2) - k$ times the action $\mathcal{G}_{k+1}(K^\prime ) \times \mathcal{C}_k (K^\prime , \gamma )_U \rightarrow\mathcal{C}_k (K^\prime , \gamma )_U \rightarrow$ in order to land inside $L^{2}_k$.


Most naturally, though, the tangent space at a given configuration $(A, \Phi , t )$ is identified with
\begin{align}
    T_{(A, \Phi , t )}\mathcal{C}_k (K^\prime , \gamma )_U = \Big\{(a , \phi , t ) \, \, | \,\, \dot{t} \in T_t \mathcal{CM}(Y, \xi_0 ) \,, \, a - \frac{\partial }{\partial \dot{t}}A_t \in L^{2}_{k}(K )  \, , \, \phi - \frac{\partial }{\partial \dot{t}} \Phi_t \in L^{2}_{k} (K) \Big\}.\label{tspace}
\end{align}

We omit the proof of the next result, which is done by carrying out the standard construction of slices for the gauge action (see \cite{KM} or \cite{marianothesis}).
\begin{Lemma}\label{banachconfspace}
The gauge group $\mathcal{G}_{k+1}(K^{\prime} )$ is a Hilbert Lie group that acts freely in a $C^{l-k-2}$ fashion on the Banach manifold $\mathcal{C}_{k}(K^{\prime} , \gamma )_U$ by $$v \cdot (t, A , \Phi ) = (t , A - v^{-1}dv, v \Phi )$$ and the quotient $\mathcal{B}_{k}(K^{\prime} , \gamma )_U = \mathcal{C}_{k}(K^{\prime} , \gamma)_U / \mathcal{G}_{k+1}(K^{\prime} )$ is naturally a $C^{l-k-2}
$ Banach manifold. 
\end{Lemma}

Consider now a second open contractible subset $\tilde{U} \subset \mathcal{CM}(Y , \xi_0 )$ together with a $\tilde{U}$-family of canonical configurations $\tilde{\gamma} = (( \tilde{A}_{t} , \tilde{\Phi}_{t} ))_{t \in \tilde{U}}$ and with $U \subset \tilde{U}$. We also assume that the families of canonical configurations $\gamma$ and $\tilde{\gamma}$ carried by $U$ and $\tilde{U}$, respectively, are in \textit{translation-invariant form} (Definition \ref{transinvform}). Then we find a unique $U$-family of gauge transformations $v_t : K \rightarrow \mathrm{U}(1)$ ($t \in U$) such that $v_t \cdot ( \Phi_{t} , A_{t} ) = (\tilde{\Phi}_{t} , \tilde{A}_{t} )$. The translation-invariance of $\gamma$ and $\tilde{\gamma}$ implies that the gauge-transformations $v_t$ are translation-invariant over the symplectic end $K = [1, + \infty) \times Y$, namely $v_t (s,y ) = v_t (1,y )$. In view of this, we may extend the $v_t$ over to the larger region $K^{\prime}$ by translation. We warn the reader that the $v_t$ need not satisfy the asymptotics $1-v_t \in L^{2}_{k+1,g_t}(K)$. However, we do obtain an inclusion map
\begin{align}
\mathcal{C}_{k}(K^{\prime}, \gamma )_U & \rightarrow   \mathcal{C}_{k}(K^{\prime}, \tilde{\gamma} )_{\tilde{U}} \nonumber \\
(t , A , \Phi )  & \mapsto   (t , A- v_{t}^{-1} d v_t , v_{t} \Phi ). \label{incl}
\end{align}

\begin{Lemma}
The map (\ref{incl}) is a well-defined smooth $\mathcal{G}_{k+1 }(K^{\prime})$-equivariant map which is an open embedding.
\end{Lemma}

\begin{proof}
The only issue which requires checking is whether (\ref{incl}) is well-defined. That is, we must check that if $(t ,A, \Phi )$ is in $\mathcal{C}_{k }(K^{\prime}, \gamma )_U $ then $(t,\tilde{A} , \tilde{\Phi }) := v_t \cdot (t , A , \Phi ) = ( t , A - v_{t}^{-1} dv_{t} , v_t \Phi )$ satisfies the conditions of Definition \ref{conf}:
\begin{itemize}
\item $\tilde{\Phi} - \tilde{\Phi}_{t} = v_t (\Phi - \Phi_{t})$. Thus, $\tilde{\Phi} - \tilde{\Phi}_{t}$ is in $L^{2}(K)$, because $\Phi - \Phi_{t} \in L^{2}(K)$ and $v_t$ has unit length
\item $\nabla_{\tilde{A}_{t}}( \tilde{\Phi} - \tilde{\Phi}_{t})  = \nabla_{A_{t} - v_{t}^{-1}dv_t } (v_t ( \Phi - \Phi_{t} ) ) = v_t \nabla_{A_{t}} (\Phi - \Phi_{t})$. Since $|v_t| = 1$ and $\nabla_{A_{t}} (\Phi - \Phi_{t}) \in L^{2} (K)$ then $\nabla_{\tilde{A}_{t}}( \tilde{\Phi} - \tilde{\Phi}_{t})$ is also in $L^{2} (K)$. Similarly, $\nabla^{l}_{\tilde{A}_{t}}( \tilde{\Phi} - \tilde{\Phi}_{t}) \in L^{2}(K)$ for all $l \geq 1$
\item $\tilde{A} - \tilde{A}_{t} = A - A_{t}$ over $K$, and so $\tilde{A} - \tilde{A}_{t} \in L^{2}_{k}(K)$.
\end{itemize}
\end{proof}

Thus, we have a \textit{directed system} whose objects are the Banach manifolds $\mathcal{C}_k (K^{\prime}, \gamma )_U$, one for each tuples $(U, \gamma )$ consisting of an open contractible set $U \subset \mathcal{CM}(Y , \xi_0  )$ carrying the family of canonical configurations $\gamma$ in translation-invariant form. A unique morphism (\ref{incl}), which is an open embedding of Banach manifolds, is associated with any two pairs $(U, \gamma ), (\tilde{U}, \tilde{\gamma} )$ such that $U \subset \tilde{U}$.

\begin{Definition} \label{confglobal} We define the \textit{configuration space} $\mathcal{C}_{k }(K^{\prime} )$ as the direct limit of the above directed system 
\begin{align*}
\mathcal{C}_{k}(K^{\prime} ) =\dirlim_{(U , \gamma)} \mathcal{C}_{k}(K^{\prime}, \gamma)_{U}.
\end{align*}
$C_{k}(K^{\prime})$ is a Banach manifold. It is the total space of a bundle of affine Hilbert spaces 
$$
\mathcal{C}_{k}(K^{\prime}) \rightarrow \mathcal{CM}(Y , \xi_0 ) 
$$
equipped with a preferred connection i.e. a complementary (horizontal) subbundle to the vertical subbundle of $T \mathcal{C}_{k}(K^{\prime} )$. Over each $U \subset \mathcal{CM}(Y , \xi_0 )$ carrying a family of canonical configurations $\gamma$ this connection induces the trivial splitting of $T \mathcal{C}_{k}(K^{\prime}, \gamma )_U$ obtained from the fact that the Sobolev spaces $L^{2}_{k, g_t , A_t}(K^{\prime} )$ are independent of $t \in U$.

We also have the configuration space modulo gauge
\begin{align*}
\mathcal{B}_{k}(K^{\prime}) =\mathcal{C}_{k}(K^{\prime} )  / \mathcal{G}_{k+1 }(K^{\prime}) \cong \dirlim_{(U , \gamma )} \mathcal{B}_{k}(K^{\prime} , \gamma)_{U}.
\end{align*}
\end{Definition}
By Proposition \ref{banachconfspace}, $\mathcal{B}_{k} (K^{\prime} )$ is a $C^{l-k-2}$ Banach manifold, and it carries a natural projection to $\mathcal{CM}(Y , \xi_0 )$.

\subsubsection{Configuration space on $Y$}\label{confspaceY}
For future reference, we also introduce here the relevant configuration spaces for the $3$-manifold $Y$. We refer the reader to \cite{KM} for further details. Given a spin-c structure $\s = (g , S , \rho )$ on $Y$, we have the configuration space $\mathcal{C}_{k-1/2}(Y , \s )$ of pairs $(B , \Psi )$ consisting of a spin-c connection $B$ and a section $\Psi$ of $S$, both of regularity $L^{2}_{k - 1/2}$. Those pairs with $\Psi $ not identically vanishing are called irreducible, and the locus of such is denoted $\mathcal{C}^{\ast}_{k-1/2}(Y , \s ) \subset \mathcal{C}_{k-1/2}(Y , \s )$. The blown-up configuration space $\mathcal{C}^{\sigma}_{k-1/2}(Y , \s )$ consists of triples $(B , s , \Psi )$ where now $s \geq 0$ is a non-negative real number, and $||\Psi||_{L^2} = 1$. The respective quotients by the (free) action of the group of $L^{2}_{k+1/2}$ gauge transformations are denoted $\mathcal{B}_{k-1/2}^{\ast}(Y , \s )$ and $\mathcal{B}^{\sigma}_{k-1/2}(Y , \s )$. They are Hilbert manifolds in a natural way [\cite{KM}, \S 9.3] (provided $k \geq 3$) and $\mathcal{B}^{\sigma}_{k-1/2}(Y , \s )$ has boundary given by configurations $(B , 0 , \Psi )$ with $||\Psi||_{L^2} = 1$.

\subsection{Moduli space and perturbations} \label{moduliandpert}
 
We now construct the promised Seiberg--Witten moduli space $\mathfrak{M}([\afr] , Z^+)$, which will be a Banach manifold equipped with a Fredholm map $\mathfrak{M}([\afr], Z^+ ) \xrightarrow{\pi} \mathcal{CM}(Y , \xi_0 ) \times \mathcal{P}$. This moduli is constructed by gluing together a moduli space over $K^{\prime}$ with a moduli space over the cylindrical end $Z = (- \infty , 0] \times Y$.

\subsubsection{The moduli space over $K^{\prime}$}

The Seiberg--Witten equations define a $\mathcal{G}_{k+1}(K^{\prime})$-equivariant section $\mathrm{sw}$ of a vector bundle $\Upupsilon_{k-1} \rightarrow \mathcal{C}_{k }(K^{\prime})$, which we now describe. On configuration spaces over an open $U \subset \mathcal{CM}(Y , \xi_0)$ equipped with a family of canonical configurations, we have the Seiberg--Witten map 
\begin{align*}
\mathrm{sw}_{\gamma , U} : \, & \mathcal{C}_{k } ( K^\prime , \gamma )_U  \rightarrow \Upupsilon_{k-1 , \gamma ,  U}\\
& (t , A , \Phi )  \mapsto ( \frac{1}{2} \rho_t (F_{\hat{A}}^{+ , g_t }) - (\Phi \Phi^{\ast})_0 , D_{A , g_t}^{+} \Phi ). \nonumber
\end{align*}
\begin{Remark}
We explain the notation from the above formula. First $\Upupsilon_{k-1 , \gamma , U}$ is the bundle over $\mathcal{C}_{k } (K^{\prime} , \gamma)_U$ with fibre over the point $(t, A , \Phi)$ given by $L^{2}_{k-1, g_t , A_{t} }( K^{\prime} , i \mathfrak{su}(S^+) \oplus S^-)$. Then $\rho_{t} ( F_{\hat{A}}^{+ , g_t} )$ is the self-adjoint endomorphism $S^+$ arising from the Clifford action of the self-dual component of the curvature $F_{\hat{A}}^{+ , g_t}$ of the $U(1)$ connection $\hat{A}$ on $\Lambda^2 S^+$. The quadratic term $(\Phi \Phi^{\ast})_{0}$ is the endomorphism which acts on a given spinor $\phi \in S^+$ by $$\phi \mapsto \inner{\Phi }{\phi} \Phi - \frac{1}{2} |\Phi|^2 \phi .$$
\end{Remark}

As before, given two open contractible subsets $U \subset \tilde{U}$ carrying canonical configurations, there is also a transition map
\begin{align*}
\Upupsilon_{k-1 , \gamma , U}&  \xrightarrow{\iota} \Upupsilon_{k-1 , \tilde{\gamma} , \tilde{U}}\nonumber \\
( ( \sigma , \Psi ) ,  (t , A , \Phi )) & \mapsto ( (\sigma , v_t \Psi ) , \iota(t , A , \Phi ))
\end{align*}
compatible with projections to the base, which thus yields a limiting bundle $\Upupsilon_{k-1} \rightarrow \mathcal{C}_{k }(K^{\prime} )$. The Seiberg--Witten maps fit in to give a commutative diagram
\[
  \begin{tikzcd}
    \mathcal{C}_{k}(K^{\prime}, \gamma )_U \arrow{r}{\iota} \arrow{d}{\mathrm{sw}_{\gamma ,  U}} & \mathcal{C}_{k}(K^{\prime} , \tilde{\gamma})_{\tilde{U}} \arrow{d}{\mathrm{sw}_{\tilde{\gamma} ,  \tilde{U}}} \\
      \Upupsilon_{k-1, \gamma , U } \arrow{r}{\iota} & \tilde{\Upupsilon}_{k-1, \tilde{\gamma} , \tilde{U}} .
  \end{tikzcd}
\] 
which provides a well-defined section $\mathrm{sw} $ of the bundle $\Upupsilon_{k-1} \rightarrow \mathcal{C}_{k }(K^{\prime})$ that we call the \textit{Seiberg--Witten map}.

In [\cite{KM}, \S 11.6], a Banach space $\mathcal{P}$ of \textit{tame} perturbations of the Chern-Simons-Dirac functional on a $3$-manifold $Y$ with a spin-c structure is constructed to achieve transversality for moduli spaces of gradient trajectories.  In our context, a suitable perturbation scheme, following the approaches of  \cite{KM} , \cite{monocont} and \cite{mariano}, is introduced as follows. Let $\mathcal{P}$ be such a Banach space of tame perturbations of the Chern-Simons-Dirac functional of $(Y , g_{\xi_0 , \alpha_0 , j_0 })$. We define a $\mathcal{G}_{k+1 }(K^{\prime} )$-equivariant section $\mu_{\gamma , U} : \mathcal{C}_{k }( K^{\prime} ,\gamma )_U  \times \mathcal{P} \rightarrow \Upupsilon_{k-1,\gamma, U}$, of the form 
\begin{align}
\mu_{\gamma , U} (t , A , \Phi , \mathfrak{p})= \varphi^1 \hat{\mathfrak{q}} ( A, \Phi ) + \varphi^2 \hat{\mathfrak{p}} ( A , \Phi )+ \varphi^3 \hat{\mathfrak{p}}_{K, t} .\label{perturbation1}
\end{align}
We describe the items appearing in (\ref{perturbation1}):
\begin{enumerate}[(i)]

\item we choose an admissible (\cite{KM}, Definition 22.1.1) perturbation $\mathfrak{q}$ of the Chern-Simons-Dirac functional on $(Y, g_{\xi_0 , \alpha_0 , j_0})$ . This induces a translation-invariant perturbation $\hat{\mathfrak{q}}(A , \Phi )$ over $\R \times Y$, as in [\cite{KM}, \S 10.1]. Then $\varphi^1$ is a smooth cutoff function on $[0 , + \infty )$, which is identically $1$ on a neighbourhood of $ 0$, and vanishes on a neighbourhood of $[1/2 , + \infty)$

\item $\mathfrak{p} \in \mathcal{P}$ induces, as before, a translation-invariant perturbation $\hat{\mathfrak{p}}$ over $\R \times Y$. We choose $\varphi^2$ to be a bump function compactly supported in $(0, 1/2) $, and identically $1$ at some interval in the interior

\item $\varphi^3$ is a cutoff function on $[ 0 , + \infty)$ which is identically $1$ over $ [1, + \infty ) $ and vanishing on a neighbourhood of $[0 , 1/2]$. We take the family of sections of $\Upupsilon_{k-1 , \gamma_0 , U}$ given by 
\begin{align*}
\hat{\mathfrak{p}}_{K,t} = ( - \frac{1}{2} \rho_t (F_{A_{t}}^{+ , g_{t} }) + ( \Phi_{t} \Phi_{t}^{\ast})_0  , 0 ).
\end{align*}
 \end{enumerate}
 The sections $\mu_{\gamma , U}$ glue to a section $\mu : \mathcal{C}_{k}(K^{\prime} ) \times \mathcal{P} \rightarrow \Upupsilon_{k-1}$, which we combine with $\mathrm{sw} : \mathcal{C}_{k}(K^{\prime} )  \rightarrow \Upupsilon_{k-1}$ to obtain the \textit{perturbed Seiberg--Witten map}:
 \begin{align}
 \mathrm{sw}_{\mu} = \mathrm{sw} + \mu : \mathcal{C}_{k}(K^{\prime} ) \times \mathcal{P} \rightarrow \Upupsilon_{k-1}. \label{SWmap}
  \end{align}

The motivation for choosing the perturbation $\hat{\mathfrak{p}}_{K}$ comes from Taubes' work \cite{taubsymp}. This perturbation term forces the canonical configurations to solve the equations $\mathrm{sw}_\mu = 0$ over the symplectic end $K \subset Z^+$. We include the  perturbations $\hat{\mathfrak{q}}$, $\hat{\mathfrak{p}}$ to achieve the necessary transversality later on.\\

\begin{Definition}The \textit{universal} moduli space of Seiberg--Witten monopoles over $K^{\prime}$ is
\begin{align*}
\mathfrak{M}_{k}(K^{\prime} ) := \mathrm{sw}_{\mu}^{-1} ( 0 ) / \mathcal{G}_{k+1}(K^{\prime} ) \cong\dirlim_{(\gamma , U )} (\mathrm{sw}_{\gamma , U } + \mu_{\gamma , U } )^{-1}(0) / \mathcal{G}_{k+1}(K^{\prime} ). 
\end{align*} 
\end{Definition}
The perturbed Seiberg--Witten map $\mathrm{sw}_{\mu}$ descends to a section on the quotient bundle $\Upupsilon_{k-1}/ \mathcal{G}_{k+1}(K^{\prime} ) \rightarrow \mathcal{B}_k (K^{\prime} ) \times \mathcal{P} .$

In \S \ref{transversality}  we will show a general transversality result (based on those of \cite{KM} and \cite{marianothesis}) which applies to the various moduli spaces that appear in this article. In particular it will give us:

\begin{Proposition}
The Seiberg--Witten map is a $C^{l-k-2}$ section of $\Upupsilon_{k-1}/ \mathcal{G}_{k+1}(K^{\prime} ) \rightarrow \mathcal{B} (K^{\prime} ) \times \mathcal{P}$ which is transverse to the zero section. Thus $\mathfrak{M}_{k}(K^{\prime} )$ is a $C^{l-k-2}$ Banach submanifold of $\mathcal{B} (K^{\prime} ) \times \mathcal{P}$.
\end{Proposition}

\subsubsection{The moduli space as a fibre product}

Using the metric $g_{\xi_0 , \alpha_0 , j_0 }$ on $Y$ and the perturbation $\mathfrak{q} \in \mathcal{P}$, one can construct the moduli space of Seiberg--Witten monopoles over the half-infinite cylinder $((- \infty , 0] \times Y , dt^2 + g_{\xi_0 , \alpha_0 , j_0 } )$ asymptotic to a critical point $[\afr]$ for the flow of the $\mathfrak{q}$-perturbed Chern-Simons-Dirac functional in the blowup. It follows that $[\afr]$ is either irreducible or unstable. This moduli is denoted $M_{k}([\afr] , (- \infty , 0] \times Y )$ and it is a Hilbert manifold. We refer the reader to \cite{KM} for details.

There are restriction maps onto the blown-up configuration space of the slice $0 \times Y$
\begin{align*}
M_{k}([\afr] , (- \infty , 0] \times Y ) \xrightarrow{R_{+}} \mathcal{B}_{k-1/2}^{\sigma}(Y , \s_{\xi_0 , \alpha_0 , j_0 } )\\
\mathfrak{M}_{k}(K^{\prime} ) \xrightarrow{\mathfrak{R}_{-}} \mathcal{B}_{k-1/2}^{\sigma}(Y , \s_{\xi_0 , \alpha_0 , j_0 } ).
\end{align*}
That the restriction maps are indeed well-defined follows by a unique continuation principle for the Seiberg--Witten equations (Proposition 10.8.1 \cite{KM}). We will see in \S \ref{transversality} that the sum of the derivatives of the restriction maps along the spinor and connection direction $$d R_{+} + d \mathfrak{R}_{-} (- , - , 0 , 0 )$$ is a Fredholm map and we will establish a transversality result:
\begin{Proposition}\label{transv1}The restriction maps $R_{+}$ and $\mathfrak{R}_{-}$ are transverse. Thus, the fibre product $\mathrm{Fib}(R_{+} , \mathfrak{R}_{+} )$ is a $C^{l-k-2}$ Banach manifold together with a Fredholm map $$\mathrm{Fib}(R_{+} , \mathfrak{R}_{+} ) \xrightarrow{\pi} \mathcal{CM}(Y, \xi_0 )\times \mathcal{P}.$$
\end{Proposition}

\begin{Definition}
The \textit{universal} moduli space of Seiberg--Witten monopoles over $Z^+ = Z \cup K^{\prime}$ associated to the triple $(\xi_0 , \alpha_0 , j_0 ) \in \mathcal{CM}(Y , \xi_0 )$ is the Banach manifold $$\mathfrak{M}([\afr] , Z^+  ) = \mathrm{Fib}(R_{+} , \mathfrak{R}_{-} ).$$
By $\mathfrak{M}(Z^+ )$ we denote the union over all critical points $[\afr]$ of the $\mathfrak{M}([\afr] , Z^+ )$. 
\end{Definition}

\begin{Remark}
    By a standard argument (see \cite{KM}, Lemma 24.2.2 and Lemma 19.1.1) one can see that any element in $\mathfrak{M}([\afr] , Z^+ ) = \mathrm{Fib}(R_+ , \mathfrak{R}_-)$ is represented by a solution $\gamma = (A, \Phi , t )$ to the Seiberg--Witten equations over the whole $Z^+$ (modulo gauge transformations $v$ with $1-v \in L^{2}_{k+1 , g_t}$ on both ends of $Z^+$) such that 
    \begin{align*}
        \gamma - (A_t , \Phi_t ) \in L^{2}_k (K )\\
        \gamma - \gamma_{w \cdot \afr} \in L^{2}_k (Z)
    \end{align*}
    where $w \in \mathcal{G}_{k+1/2}(Y)$, $\afr$ is a critical point of the $\mathfrak{q}$-perturbed Chern-Simons-Dirac functional, $\gamma_{w \cdot \afr}$ is the translation-invariant solution over the cylindrical end $Z$ determined by $w \cdot \afr$, and $(A_t , \Phi_t )$ is a canonical configuration over $K$ (in translation-invariant form).
\end{Remark}



\subsubsection{Components of the moduli space of constant index}

As with the moduli spaces that are studied in \cite{KM}, the index of $\pi$ will vary with the connected component of $\mathfrak{M}([\afr] , Z^+ )$. We give a more precise statement of this fact, following the ideas of \S 24.4 in \cite{KM}. Denote by $\mathcal{B}_{k}([\afr], K^{\prime})$ the preimage of $[\afr]$ under the partially defined restriction map to the slice $0 \times Y$
\begin{align}
R_- : \mathcal{B}_{k}(K^{\prime}) \dashrightarrow \mathcal{B}_{k-1/2}^{\sigma}(Y, \s_{\xi_0 , \alpha_0 , j_0 }) .\label{restriction1}
\end{align}
Any element of $\mathfrak{M}([\afr] , Z^+)$ is, by definition, given by a quadruple $([\gamma_{Z} ] , [\gamma_{K^{\prime}}] , t , \mathfrak{p} )$ with $[\gamma_{Z}|_{Y}] = [\gamma_{K^{\prime}}|_Y]$. The cylinder configuration $[\gamma_Z ]$ provides a path in $\mathcal{B}_{k}(K^{\prime})$ (canonical up to homotopy) from $[\gamma_{K^\prime}] \in \mathcal{B}_k ( K^\prime )$ to the subspace $\mathcal{B}_{k}([\afr] , K^{\prime} ) \subset \mathcal{B}_k ( K^\prime )$. Hence, each element of $\mathfrak{M}([\afr] , Z^+ )$ determines a connected component of $\mathcal{B}_{k}([\afr] , K^{\prime} )$, giving a map 
\begin{align}
\pi_0 \mathfrak{M}([\afr] , Z^+ ) \rightarrow \pi_0 \mathcal{B}_{k}([\afr] , K^{\prime} ). \label{connectedcomps}
\end{align}

By the homotopy invariance of the index of a Fredholm operator we have:
\begin{Proposition}\label{constantindex}
The index of $\pi : \mathfrak{M}([\afr] , Z^+ ) \rightarrow \mathcal{CM}(Y , \xi_0 ) \times \mathcal{P}$ is constant on the fibres of (\ref{connectedcomps}).
\end{Proposition}

Next we provide further information on $\pi_0 \mathcal{B}_k ([\afr] , K^{\prime} )$. Consider the natural projection $p :\mathcal{B}_k ([\afr] , K^{\prime} ) \rightarrow \mathcal{CM}(K , \xi_0 )$, and denote the fibre over a point $t$ by $\mathcal{B}_{k}([\afr] , Z^+ )_t$.

\begin{Lemma}\label{componentslem}
\begin{enumerate}[(i)]
 \item there is a bijection $\pi_0 \mathcal{B}_{k}([\afr] , K^{\prime} )_t \approx  H^{1}(Y ; \Z )$
 \item $p$ is a Serre fibration
 \item the map $\pi_0  \mathcal{B}_k ([\afr] , K^{\prime})_t \rightarrow \pi_0  \mathcal{B}_k ([\afr] , K^{\prime})$ induced by inclusion is surjective.
  \end{enumerate}
\end{Lemma}
\begin{proof}

For (i), we fix a canonical configuration $\gamma_t := ( A_{t} , \Phi_{t})$ at $t$, and fix a representative $\afr$ of $[\afr]$. We consider the space $\mathcal{C}_k (\afr ,  K^{\prime} , \gamma_t )$ which is the fibre of the partially-defined restriction map $ \mathcal{C}_{k}(K^{\prime} , \gamma_t ) \dashrightarrow \mathcal{C}_{k-1/2}^{\sigma}(Y , \s_{\xi_0 , \alpha_0 , j_0 } )$ over $\afr$. We choose a representative $v_{z}$ from every connected component $z \in \pi_0 \mathcal{G}_{k+1/2}(Y)$. Then we have a decomposition into disjoint closed subspaces $$ \mathcal{B}_{k}([\afr] , K^{\prime} )_t = \bigcup_{z \in \pi_0 \mathcal{G}_{k+1/2} (Y) } \mathcal{C}_{k}(v_z \cdot \afr , K^{\prime} , \gamma_t )/ \mathcal{G}_{k+1}(K^{\prime}, Y )$$
where $\mathcal{G}_{k+1}(K^\prime , Y )$ stands for the subgroup of $\mathcal{G}_{k+1}(K^\prime )$ consisting of gauge transformations which are the identity over $\{ 0 \} \times Y = \partial K^\prime $. Note that each of the disjoint subspaces above is connected (because $\mathcal{G}_{k+1}(K^{\prime})$ is connected). This sets up a bijection $\pi_0 \mathcal{B}_{k}([\afr] , K^{\prime})_t \cong \pi_0 \mathcal{G}_{k+1/2}(Y)$ only depending on the representative $\afr$. Finally, the latter set is identified with the group $H^{1}(Y ; \Z ) = [Y , S^1]$. 

Part (iii) follows from (ii) and the connectedness of the base of the Serre fibration $p$. For (ii) we must show: if $h : D \times [0,1] \rightarrow \mathcal{CM}(Y, \xi_0 )$ is any given homotopy, where $D$ is a compact disc, and we are given a lift of $h$ over $D \times 0$ to $\mathcal{B}_k ([\afr] , K^{\prime} )$, then there exists a lift of $h$ over $D \times [0,1]$ agreeing with the given one over $D \times \{0\}$. It suffices to consider the fibration $p^\prime$ obtained by pullback of $p$ along $h$ 
\[
p^\prime : h^\ast \mathcal{B}_k ([\afr] , K^\prime ) \rightarrow D \times [0,1]
\]
and construct a section of $p^\prime$ over $D\times [0,1]$ which extends the given one over $D \times \{0\}$. This has the advantage that $D \times [0,1]$ is contractible, hence we can find a family of canonical configurations $(A_{t,\tau }, \Phi_{t, \tau} )$ carried by $(t, \tau ) \in D \times [0,1]$. The given section over $D \times \{0\}$ can be represented by a family of configurations of the form 
\[
\gamma_t = (A_{t, 0} , \Phi_{t, 0}) + ( a_t ,  \phi_t ) 
\] 
such that $\gamma_t |_{\{0\} \times Y} = \mathfrak{a}$ for some representative $\afr$ of $[\afr]$, and $(a_t , \phi_t ) \in L^{2}_k$. Let $\bfr_{t,\tau}$ be the restriction to $\{0\} \times Y$ of $ (A_{t,\tau}, \Phi_{t, \tau}) + (a_t , \phi_t ) $. We extend the section $\gamma_t$ over to $D \times [0,1]$ by setting
\[
\gamma_{t, \tau} = (A_{t, \tau} , \Phi_{t , \tau} ) + (a_t, \phi_t ) + \beta \cdot (\gamma_{\afr} - \gamma_{\bfr_{t, \tau}} )
\]
where $\beta$ is a smooth bump function $\beta : [0,+ \infty) \rightarrow \mathbb{R}_{\geq 0}$ which is identically $1$ on a small neighbourhood of $0$ and vanishes outside of a compact set, and $\gamma_{\afr}$ stands for the translation-invariant configuration on the cylinder $\mathbb{R} \times Y$ associated to $\afr$ (and likewise for $\gamma_{\bfr_{t, \tau}}$). The proof is now complete.
\end{proof}

Thus, by the previous results we can decompose $\mathfrak{M}([\afr] , Z^+)$ into pieces where $\pi$ has constant index
\begin{align*}
\mathfrak{M}([\afr] , Z^+ ) = \bigcup_{z} \mathfrak{M}_z ([\afr] , Z^+) .
\end{align*}
which are parametrised by the connected components $z \in \pi_0  \mathcal{B}_k ([\afr] , K^{\prime})$. Note that it does \textit{not} hold necessarily that each $\mathfrak{M}_z ([\afr] , Z^+)$ is connected.

\begin{Remark}
The map from Lemma \ref{componentslem}(iii) is not injective, in general. More precisely, it is injective if and only if any loop (i.e. $S^1$-family) in $\mathcal{CM}(Y, \xi_0  )$ has a corresponding $S^1$-family of canonical configurations over $K$. 
\end{Remark}

 
\subsubsection{Orientability}\label{orientability1}

In order to define the families invariant when the coefficient ring $R$ is not of characteristic $2$, we will need to orient the Seiberg--Witten moduli spaces. In order to do so we need to orient the \textit{determinant line bundle} of the Fredholm map $$\mathfrak{M} (Z^+ ) = \bigcup_{[\afr] } \mathfrak{M}([\afr] , Z^+ ) \xrightarrow{\pi} \mathcal{CM}(Y , \xi_0 ) \times \mathcal{P}. $$For the precise construction of this real line bundle $\mathrm{det} \pi \rightarrow \mathfrak{M}(Z^+ )$ we refer to [\cite{KM} , \S20.2]. Its fibre over a given $m \in \mathfrak{M}(Z^+)$ can be identified as $$(\mathrm{det} \pi)_m  = \Lambda^{\mathrm{max}} \mathrm{ker} (d \pi )_m \otimes \Lambda^{\mathrm{max}} (\mathrm{coker} (d \pi)_m  )^{\ast}. $$ We now describe what goes into orienting the determinant line bundle.

The first ingredient is to orient the moduli spaces of trajectories in monopole Floer homology. This is formally analogous to the finite-dimensional Morse theory case. Given a critical point $[\afr] \in \mathcal{B}^{\sigma}(Y , \s )$ in the blowup, a $2$-element set $\Lambda ([\afr] ) $ is associated in [\cite{KM}, \S 20.3], playing the role of the set of orientations for the unstable manifold of $[\afr]$ in the Morse theory picture. 

The second ingredient is the following construction of a double covering $\Lambda$ of $\mathcal{C}(Y, \xi_0 )$, in the spirit of [\cite{KM}, \S 24.8]. Let $\gamma_0 = (t, A , \Phi  ) \in \mathcal{B}_k (K^\prime )$ be a configuration which restricts along $\{0\} \times Y $ onto a \textit{reducible} configuration $[\afr_0 ] \in \mathcal{B}^{\sigma}_{k-1/2}(Y, \s_{\xi_0 , \alpha_0 , j_0 } )$. In other words, the configuration $\gamma_0$ has image $[\afr_0]$ under the partially-defined map from (\ref{restriction1})
\[
R_- : \mathcal{B}_k (K^\prime ) \dashrightarrow \mathcal{B}_{k-1/2}^\sigma (Y, \s_{\xi_0 , \alpha_0 , j_0 } ).
\]
The vector bundle bundle $\Upsilon_{k-1}$ on the right-hand side of (\ref{SWmap}) descends onto a vector bundle denoted $[\Upsilon_{k-1}]$ over $\mathcal{B}_k (K^\prime )$, and the (unperturbed) Seiberg-Witten map gives a section of this bundle
\[
\mathrm{sw} : \mathcal{B}_k (K^\prime  ) \rightarrow [\Upsilon_{k-1}].
\]
Consider also the map 
\[\pi_{[\afr_0]}: T_{[\afr_0]}\mathcal{B}_{k-1/2}^{\sigma}(Y, \s_{\xi_0 , \alpha_0 , j_0 }) \rightarrow \mathcal{K}^{+}_{[\afr_0]}
\]
given by orthogonal projection onto the subspace $\mathcal{K}^{+}_{[\afr]}$ which is defined as the closure of the span of the \textit{non-negative} eigenvectors of the (unperturbed) Hessian $\mathrm{Hess}^{\sigma} :T_{[\afr_0]}\mathcal{B}_{k}^{\sigma}(Y, \s_{\xi_0 , \alpha_0 , j_0 }) \rightarrow T_{[\afr_0]}\mathcal{B}_{k-1}^{\sigma}(Y, \s_{\xi_0 , \alpha_0 , j_0 })$ of the Chern--Simons--Dirac functional. From these we assemble the Seiberg-Witten map with \textit{Atiyah--Patodi--Singer boundary condition}
\[
P_{\gamma_0} = (\mathcal{D} \mathrm{sw})_{\gamma_0} \oplus \pi_{[\afr_0]}\circ (d R_- )_{\gamma_0} : T_{\gamma_0} \mathcal{B}_k (K^\prime ) \rightarrow [\Upsilon_{k-1}]_{\gamma_0} \times \mathcal{K}^{+}_{[\afr_0]}
\]
where $(\mathcal{D} \mathrm{sw})_{\gamma_0}$ stands for the vertical component of the derivative of the section $\mathrm{sw}$ at $\gamma_0$ (taken with respect to the natural connection that the vector bundle $[\Upsilon_{k-1}]$ carries). As in \cite{KM}, the Atiyah--Patodi--Singer theory establishes the Fredholm property of the linear operator $P_{\gamma_0}$.

\begin{Definition}\label{orientationlocalsystem}
We define a double covering $\Lambda$ of $\mathcal{C}(Y , \xi_0 )$ with fibers $\Lambda(\xi) $ as follows. For a given $\xi \in \mathcal{C}(Y , \xi_0 )$ choose any configuration $\gamma_0 = (t , A , \Phi ) \in \mathcal{B}_k (K^\prime )$ lying over $\xi $ (i.e. $t = (\xi, \alpha , j )$ for some $\alpha$ and $j$) and which restricts onto a \textit{reducible} configuration $[\afr_0] \in \mathcal{B}_{k-1/2}^{\sigma}(Y , \s_{\xi_0 , \alpha_0 , j_0 })$ along $\{0\} \times Y$. We define $\Lambda(\xi) $ to be the two-element set of orientations of $\mathrm{det}(P_{\gamma_0} )$.
\end{Definition}

\begin{Remark}
 The two-element set $\Lambda (\xi )$ is independent of the choice of $\gamma_0$ or $[\afr_0 ]$, up to canonical bijection. Furthermore, it is also independent of our chosen base configuration $(\xi_0 , \alpha_0 , j_0 )$, up to canonical bijection. These assertions all follow from [\cite{KM}, Lemma 20.3.3] and standard homotopy arguments as those found in [\cite{KM}, \S 20.3 and \S 24.8].
\end{Remark}

Associated to our double cover $\Lambda$ there is a local system $\Lambda_\Z $ whose fibers are $\Z$-modules of rank 1. Explicitly, we can take the fiber $\Lambda_\Z (\xi )$ to be the quotient of the free $\Z$-module on the two-element set $\Lambda(\xi) = \{ \mathfrak{o} , \mathfrak{o}^{\prime}\}$ by the submodule generated by the element $\mathfrak{o}+\mathfrak{o}^{\prime} $, and the monodromy action of paths is inherited from that on $\Lambda$. We write $\Lambda _R $ for the local system of free $R$-modules of rank 1 obtained by taking the tensor product $\Lambda_\Z\otimes_{\Z} R$.

The proof of the next result follows the same argument as in \S 24.8 of \cite{KM}
\begin{Proposition}\label{orientability}
Given a choice of an element in each orientation set $\Lambda ([\afr])$ for each critical point $[\afr]$, there is a canonical homotopy class of isomorphism of real line bundles $\mathrm{det} \pi \cong \pi^{\ast}\Lambda_\R$ over $\mathfrak{M}(Z^+)$. Here $\mathrm{det} \pi$ is the determinant line bundle of $\pi : \mathfrak{M}(Z^+ ) \rightarrow \mathcal{CM}( Y , \xi_0 ) \times \mathcal{P}$, and $\pi^{\ast} \Lambda_{\R}$ is the pullback of $\Lambda_{\R}$ by $ \pi : \mathfrak{M}(Z^+) \rightarrow \mathcal{CM}( Y , \xi_0 ) \times \mathcal{P} \simeq \mathcal{C}(Y , \xi_0 )$.
\end{Proposition}

We explain how this orients the moduli spaces that will be relevant. Consider a $C^2$ map $\sigma : \Delta^n \rightarrow \mathcal{CM}( Y , \xi_0 ) \times \mathcal{P}$ from the standard $n$-simplex $\Delta^n = \{ x \in (\R_{\geq 0})^{n+1} \, | \, \sum_{i = 1}^{n+1} x_i = 1\}$. We equip $\Delta^n$ with its canonical orientation. Suppose $\sigma$ is transverse to $\pi : \mathfrak{M}_{z}(Z^+) \rightarrow \mathcal{CM}( Y , \xi_0 ) \times \mathcal{P}$ along each stratum of $\Delta^n$. Then we obtain a $C^2$ manifold $M_z ([\afr] , \sigma ) = \mathrm{Fib}(\pi , \sigma ) $ as the fibre product of $\pi :\mathfrak{M}_{z}([\afr] , Z^+ ) \rightarrow \mathcal{CM}( Y , \xi_0 ) \times \mathcal{P} $ with $\sigma $, which is of dimension $\mathrm{ind} \pi + n$, where $\mathrm{ind} \pi$ is computed over the component $\mathfrak{M}_{z}([\afr] , Z^+ )$. 

If choices in each $\Lambda ([\afr])$ are made and we are given an orientation in $\Lambda (\sigma(b))$, where $b$ stands for the barycenter of $\Delta^n$, then Proposition \ref{orientability} picks out preferred orientations of all the moduli spaces $M_{z}([\afr] , T )$. This is is a matter of linear algebra:

\begin{Lemma}\label{signrule}
Consider transverse linear maps $ M \xrightarrow{\pi} C \xleftarrow{\sigma} \Delta $ of Banach spaces, with $\pi$ Fredholm and $\Delta$ of finite dimension. Let $F = \pi- \sigma : M \oplus \Delta \rightarrow C$. Then $F$ is Fredholm and there is a canonical isomorphism $$ \mathrm{det} F \cong  \mathrm{det}\pi \otimes \Lambda^{\mathrm{max}} \Delta  .$$
\end{Lemma}
 \begin{proof}
Because of the transversality assumption, one has the canonical isomorphism  (see the construction of [\cite{KM},\S 20.2], and put $J = \mathrm{Im}\sigma$)
$$\mathrm{det} \pi \cong \Lambda^{\mathrm{max}} \pi^{-1} (\mathrm{Im}\sigma ) \otimes \Big( \Lambda^{\mathrm{max}} \mathrm{Im} \sigma \Big)^{\ast}.$$
Then the short exact sequences 
\begin{align*}
& 0 \rightarrow \mathrm{Ker}\sigma \rightarrow \mathrm{Ker}F \rightarrow \pi^{-1} ( \mathrm{Im}\sigma ) \rightarrow 0 \\
& 0 \rightarrow \mathrm{Ker}\sigma \rightarrow \Delta \rightarrow \mathrm{Im} \sigma \rightarrow 0
\end{align*}
provide us with canonical isomorphisms
\begin{align*}
\Lambda^{\mathrm{max}} \pi^{-1} (\mathrm{Im}\sigma ) \cong   \Lambda^{\mathrm{\max}}\mathrm{Ker} F \otimes \Big( \Lambda^{\mathrm{max}} \mathrm{Ker}\sigma \Big)^{\ast} \\ \cong \Lambda^{\mathrm{\max}}\mathrm{Ker} F \otimes \Big( \Lambda^{\mathrm{max}} \Delta \Big)^{\ast} \otimes \Lambda^{\mathrm{max}} \mathrm{Im} \sigma .
\end{align*}
This says $\mathrm{det} \pi \cong \mathrm{det} F \otimes \Big( \Lambda^{\mathrm{max}} N \Big)^{\ast}$.
\end{proof}

More precisely, one orients $M_z ([\afr] , \sigma)$ by following the proof of Lemma \ref{signrule} above using the \textit{fibre-first convention} for orienting vector spaces in a short exact sequence. This agrees with orientation conventions in \cite{KM} (see p.525) for parametrised moduli spaces over an oriented manifold.

We refer to this as the \textit{canonical orientation} of $M_{z}([\afr] , \sigma )$ (depending on the choices of elements in $\Lambda([\afr])$ and $\Lambda(\sigma(b))$). Whenever these moduli are $0$-dimensional and we use them to make counts of points, each point is counted with a sign corresponding to its canonical orientation (relative to the natural orientation of a point).
 
\subsection{The families contact invariant}\label{theinvariant}

We describe now the construction of the homomorphism (\ref{invtmap}).  We will write $\mathcal{C}$ for the Banach manifold $\mathcal{CM}( Y , \xi_0 ) \times \mathcal{P}$ for ease in notation. This space has the weak homotopy type of the space of contact structures $\mathcal{C}(Y , \xi_0 )$.

We fix orientations in $\Lambda ([\afr])$ for all critical points $[\afr]$. We fix a ring $R$ (commutative, unital).

\subsubsection{Transverse singular chains}\label{transversechains}

Let $M \xrightarrow{\pi} C$ be a $C^r$ Fredholm map of $C^r$ Banach manifolds. We assume that $C$ is connected but allow $M$ disconnected, with at most countably many components. The index of $\pi$, $\mathrm{ind} \pi \in \Z$, depends on the chosen connected component of $M$. 

Below we view the standard $n$-simplex $\Delta^n = \{ x \in (\R_{\geq 0})^{n+1} \, | \, \sum_{i = 1}^{n+1} x_i = 1\}$ as a manifold with corners, and by a $C^r$ map with domain $\Delta^n$ we mean a map which extends to a $C^r$ map on an open neighbourhood of $\Delta^n \subset \R^{n+1}$.

\begin{Definition}
A $C^r$ singular $n$-simplex $\sigma : \Delta^n \rightarrow C$ is \textit{transverse} to $\pi$ if the restriction of $\sigma$ to each stratum (i.e. face) of the $n$-simplex $\Delta^n$ is transverse to $\pi$. In particular, the image of each vertex of $\Delta^n$ under $\sigma$ is a regular value of $\pi$.
\end{Definition}

For our purposes it suffices to take $r = 2$. Next we set up a version of the complex of singular chains on $\mathcal{C}$ with coefficients in the local system $\Lambda_R$, made up of transverse chains

\begin{Definition}\label{transvsimplex}
Let $( S_{\ast}^{\pi}(C ; \Lambda_R ) , \partial )$ be the chain complex over $R$ given by finite formal sums$$
\sum a \cdot \sigma
$$where $\sigma$ is a $C^2$ singular simplex $\sigma : \Delta^{n} \rightarrow C$ (with $n \geq 0$) which is \textit{transverse} to $\pi$ along components of $M$ with $\mathrm{ind} \pi  \leq 1-n$; and $a$ is an element of the ring $\Lambda_{R}(\sigma (b) )$, where $b \in \Delta^n$ is the barycenter of $\Delta^n$. The differential $\partial$ is the singular differential coupled to the isomorphism $\Lambda(\sigma(b)) \rightarrow \Lambda (\sigma_i (b_i) )$ associated with the straight line segment from $b$ to $b_i$, where $\sigma_i $ denotes the restriction of $\sigma$ to the $i$th codimension 1 face $\Delta^n_{i}$ of $\Delta^n$ and $b_i$ the barycenter of $\Delta^{n}_{i}$.
\end{Definition}

\begin{Remark}
For ease in notation, whenever we refer to a singular $n$-simplex $\sigma : \Delta^n \rightarrow \mathcal{C} $ we will assume it is equipped with an element in $\Lambda(\sigma(b))$, and regard instead the coefficient $a$ as an element in the ring $R$.
\end{Remark}

The restriction to components of $M$ with $\mathrm{ind} \pi \leq 1-n$ is imposed on us by the Thom-Smale transversality theorem \cite{sardsmale}. This result states that for $C^r$ maps ($r \geq 1$) of $C^r$ Banach manifolds $X \xrightarrow{f} Y \xleftarrow{g} Z $ with $\mathrm{dim} X = n < + \infty$ and $g$ Fredholm, one can always $C^r$-approximate $f$ by a map $f^{\prime}$ which is transverse to $g$, provided that $r > \mathrm{max}  (\mathrm{ind}g+ n , 0 )$. Furthermore, if $f$ was already transverse to $g$ along a closed subset $X^\prime \subset X$ then one can choose $f^{\prime}$ to agree with $f$ along $X^{\prime}$. Then, by the Thom-Smale transversality theorem we learn that the inclusion of $S_{\ast}^{\pi}(C; \Lambda_R)$ into the chain complex of (continuous) singular chains on $C$ with coefficients in the local system $\Lambda_R$ induces a quasi-isomorphism, so that $( S_{\ast}^{\pi}(C; \Lambda_R) , \partial )$ computes the singular homology $H_{\ast}(C ; \Lambda_R ) \cong H_{\ast} (\mathcal{C}(Y , \xi_0 ) ; \Lambda_R )$.

\subsubsection{Counting solutions to the Seiberg--Witten equations}


Consider a $C^2$ singular $n$-simplex $\sigma : \Delta^n \rightarrow \mathcal{C}$ satisfying the transversality condition of Definition \ref{transvsimplex} with respect to the Fredholm map $\pi : \mathfrak{M}(Z^+ ) \rightarrow \mathcal{C}$ (of regularity $C^{l-k-2} \subset C^2$). For such $\sigma$ and each pair $([\afr] , z)$ we have the space $M_{z}([\afr] , \sigma )$  consisting of solutions of the Seiberg--Witten equations over the singular simplex $\sigma$. Namely, $M_z ([\afr] , \sigma ) = \mathrm{Fib}(\pi , \sigma )$ is the fibre product of $\pi : \mathfrak{M}_z ([\afr] , Z^+ ) \rightarrow \mathcal{C}$ and $\sigma$. Whenever the expected dimension of $M_{z}([\afr], \sigma )$ is $\leq 1$, i.e. $\mathrm{ind} \pi \leq 1-n$, we can guarantee that this fibre product is transverse, and hence that $M_{z}([\afr] , \sigma )$ will be a $C^2$-manifold with corners. We denote by $\# M_z ([\afr] , \sigma )$ the count of points in the discrete ($0$-dimensional) moduli space $M_{z}([\afr] , \sigma )$ when $\mathrm{ind} \pi = -n$, counted with the signs corresponding to their canonical orientation (see \S \ref{orientability1}); and we set $\# M_z ([\afr] , \sigma ) = 0$ if $\mathrm{ind} \pi \neq - n $. The possibility to make such count relies on the fact that the $0$-dimensional moduli spaces $M_{z}([\afr] , \sigma )$ are indeed finite, which we will address momentarily.

We can now assemble the counts of solutions to the Seiberg--Witten equations into a homomorphism of $R$-modules
\begin{align}
\psi : S_{\ast}(\mathcal{C} ; \Lambda_R )&  \rightarrow \widehat{C}^{\ast}(Y , \s_{\xi_0 , \alpha_0 , j_0 } ; R ) \label{chainmap}\\
\sigma  \mapsto & \, \mathfrak{M}(Z^+ ) \cdot \sigma :=  \sum_{[\afr] , z } \big(  \# M_{z}([\afr] , \sigma )\big) \cdot [\afr] . \nonumber
\end{align}

 The right side of (\ref{chainmap}) is the monopole Floer \textit{cochain} complex of $Y$ (in the \textit{from} version), obtained by taking the dual of the monopole Floer chain complex $\widehat{C}_{\ast} (Y , \s_{\xi_0 , \alpha_0 , j_0} ; R)$ with differential $\widehat{\partial}$. The latter complex  is constructed from the spin-c structure $\s_{\xi_0 , \alpha_0 , j_0}$ and admissible perturbation $\mathfrak{q}$. It is freely generated over $R$ by the union of the sets $\mathfrak{C}^{o}$, $\mathfrak{C}^{u}$ of irreducible and unstable critical points, which gives a decomposition $\widehat{C}_{\ast}(Y, \s_{\xi_0 , \alpha_0 , j_0} ) = C^{o}_{\ast} \oplus C^{u}_{\ast} $. The Floer differential is given by the following matrix (see \cite{KM}, Definition 22.1.3)
\begin{align}
&\widehat{\partial} = \begin{pmatrix} \partial^{o}_{o} & \partial^{u}_{o} \\ - \overline{\partial}^{s}_{u} \partial^{o}_{s} & - \overline{\partial}^{u}_{u} - \overline{\partial}^{s}_{u} \partial^{u}_{s} \end{pmatrix}. \label{floerdiff}
\end{align}

\begin{Remark}
For the expression (\ref{chainmap}) to be well-defined, we require the fact that there are only finitely many pairs $([\afr] , z)$ for which $M_{z}([\afr] , \sigma )$ is of dimension $0$ and non-empty. This can be shown following the standard arguments in \cite{KM}, and we defer a discussion of this fact to \S \ref{compactness}. 
\end{Remark}

\begin{Proposition}\label{chainmapprop}
Up to signs, $\psi$ is a chain map. Precisely, $\psi( \partial \sigma ) = (-1)^{n}\widehat{\partial}^{\ast} \psi ( \sigma )$, where $\sigma$ is a singular $n$-simplex.
\end{Proposition}

To see this, we first make some remarks on the compactness properties of the moduli spaces $M_{z}([\afr] , \sigma)$. We restrict to the case of the moduli spaces of expected dimension $\leq 1$, since for the higher dimensional ones we cannot guarantee that they are transversely cut out. The $M_{z}([\afr] , \sigma )$ are, in general, non-compact manifolds with corners. However, we have
\begin{Proposition}\label{strata1}
The $0$-dimensional moduli spaces $M_{z}([\afr] , \sigma )$ consist of finitely-many points. The $1$-dimensional moduli spaces $M_{z}([\afr] , \sigma )$ admit a compactification into a space $M^{+}_{z}([\afr] , \sigma )$ stratified by manifolds. The top stratum consists of $M_{z}([\afr] , \sigma )$ itself, and the boundary of the top stratum consists of ``broken'' configurations of the form
\begin{enumerate}[(a)]
\item $\breve{M}_{z_1}([\afr] , [\bfr] ) \times M_{z_0}([\bfr] ,\sigma )$
\item $\breve{M}_{z_2}([\afr] , [\bfr] ) \times \breve{M}_{z_1}([\bfr] , [\cfr] ) \times M_{z_0}([\cfr] , \sigma )$
\end{enumerate}
where the middle factor in (b) is boundary-obstructed; together with configurations arising from the boundary stratum of $\Delta^n$, which is the union of the $(n-1)$-simplices $\Delta_{0}^{n-1} , \Delta_{1}^{n-1} ,  \ldots ,  \Delta_{n}^{n-1}$ that are codimension-$1$ faces of $\Delta^n$:
\begin{enumerate}[(c)]
\item $\bigcup_{i = 0 , n } M_{z}([\afr] , \sigma_{|  \Delta_{i}^{n-1}} )$ .
\end{enumerate}
For each boundary stratum above, the homotopy classes must concatenate to $z$ (e.g. for (a) we need $z_1 \circ z_0 = z$). Furthermore, the structure near each boundary stratum is: $C^0$ manifold-with-boundary structure at (a); a \textit{codimension}-$1$ $\delta$-structure (a more general structure than $C^0$ manifold-with-boundary, see [\cite{KM}, Definition 19.5.3]) at (b); and a $C^2$ manifold-with-boundary structure at (c).
\end{Proposition}

All the analysis required to deduce these results is provided by the techniques in \cite{KM}, \cite{monocont}, \cite{boyu}. We discuss in \S \ref{compactness} some technical results that are involved. 

\begin{proof}[Proof of Proposition \ref{chainmapprop}]
In general, for any singular simplex $\sigma$ transverse to $\pi : \mathfrak{M}_{z}([\afr] , Z^+) \rightarrow \mathcal{C}$, one can construct a compactification $M_{z}^{+}([\afr] , \sigma )$ of $M_{z}([\afr] , \sigma )$ by adding broken configurations as in \cite{KM}. In the case when $M_{z}([\afr] , \sigma)$ is transversely cut out and of dimension $0$, it follows from index reasons that no broken configurations are added, and so the moduli consists of finitely-many points. In the case where the dimension of $M_{z}([\afr] , \sigma )$ is $1$, the corresponding compactification $M^{+}_{z}([\afr] , \sigma)$ is a $1$-dimensional stratified space with a codimension-$1$ $\delta$-structure along its boundary. Such a space enjoys the nice property that the enumeration of its boundary points  gives total count zero [\cite{KM}, Corollary 21.3.2]. Thus, enumerating the boundary points, of types (a), (b) and (c) as above, yields corresponding identities 
\begin{align*}
\inner{ \psi(\sigma)^{o}}{\partial^{o}_{o} [\afr] } - \inner{ \psi (\sigma)^{u}}{ \overline{\partial}^{s}_{u} \partial^{o}_{s} [\afr] }  + (-1)^{n-1}\inner{ \psi (\partial \sigma )^{o}}{ [\afr] } = 0\, , \quad &\forall [\afr] \in \mathfrak{C}^{o} \\
\inner{\psi (\sigma)^{o}}{\partial^{u}_{o}[\afr] } + \inner{\psi (\sigma)^{u}}{\overline{\partial}^{u}_{u} [\afr] } - \inner{ \psi(\sigma)^{u}}{\overline{\partial}^{s}_{u} \partial^{u} [\afr]} + (-1)^{n-1}\inner{ \psi (\partial\sigma )^{u} }{[\afr] } = 0 \, , \quad &\forall [\afr] \in \mathfrak{C}^{u}
\end{align*}
which give the required equality $\psi ( \partial \sigma ) =(-1)^{n} \widehat{\partial}^{\ast} \psi ( \sigma )$. For the origin of the signs see Lemma \ref{signs1} \footnote{A rather technical point is that the sign of $+ \inner{\psi (\sigma)^{u}}{\overline{\partial}^{u}_{u} [\afr] }$ written above should be flipped if one follows the \textit{reducible convention} for orienting the moduli $M_{z_1}([\afr] , [\bfr] )$ when both $[\afr], [\bfr]$ are boundary-unstable (see \S 20.6 \cite{KM}). This reducible convention is meant when writing the term $- \overline{\partial}^{u}_{u}$ in the Floer differential (\ref{floerdiff}). The signs listed in Lemma \ref{signs1} follow the usual convention. These two conventions differ by the sign $(-1)^{\mathrm{dim}M_{z_1}([\afr] , [\bfr] )} = -1 $.} in \S \ref{compactness}.
\end{proof}

\begin{Definition}
The \textit{families contact invariant} of $(Y , \xi_0 )$ is the homomorphism induced by the chain map $\psi$
\begin{align}
 \fc := \psi_\ast : H_{\ast}(\mathcal{C}(Y , \xi_0 ) ; \Lambda_R ) \cong H_{\ast}(\mathcal{C} ; \Lambda_R ) \rightarrow \widehat{HM}^{\ast}(Y , \s_{\xi_0} ; R ).  \label{invtmap2}
\end{align}
\end{Definition}

Some observations are relevant now:

\begin{Remark}
\begin{enumerate}[(i)]
\item \textit{Invariant for a single contact structure}. Fixing an element of the $2$-element set $\Lambda(\xi_0 )$ fixes the sign of the contact invariant $\mathbf{c}(\xi_0)$ of Kronheimer--Mrowka-Ozsv\'{a}th-Szab\'{o} \cite{monolens}. In turn, this also picks out a canonical generator $1 \in  H_{0}( \mathcal{C}(Y , \xi_0 ) ; \Lambda_R )$($= R$ or $R/2R$ accoding as to whether the local system $\Lambda$ is trivial or not). It is clear from our construction that  $\mathbf{c}(\xi_0 )$ agrees with $\fc(1)$. Part (A) of Theorem \ref{mainthm} is then proved.

\item \textit{Gradings}. With respect to the natural grading of the Floer cohomology groups by the set of homotopy classes of oriented $2$-plane fields, the map $\psi$ defining $(\ref{invtmap})$ has the form 
\begin{align*}
\psi : S^{\pi}_n (\mathcal{C}(Y , \xi_0 ) ; \Lambda_R ) \rightarrow \widehat{C}^{[\xi_0 ] - n }(Y , \s_{\xi_0} ; R ), \quad n \geq 0. 
\end{align*}
For $n = 0$, i.e. for the contact invariant $\mathbf{c}(\xi_0 )$, a proof of this fact can be found in [\cite{marianothesis}, \S 7.1]. For higher $n \geq 0$ the statement follows in a straightforward way from the $n = 0$ case and the identity of expected dimensions $\mathrm{dim} M_{z}([\afr] , \sigma ) = n + \mathrm{dim} M_{z}([\afr] , \ast )$, with $\sigma : \Delta^n \rightarrow \mathcal{C}$ an $n$-simplex and $\ast : \{ \ast \} \rightarrow \mathcal{C}$ the inclusion of a point.

\item \textit{Criterion for triviality of $\Lambda$}. It is unclear to the author whether the double cover $\Lambda$ can be non-trivial in general. However, under the assumption that the contact invariant $\mathbf{c}(\xi_0 )$ with $R = \Z$ coefficients is \textit{not} a 2-torsion element, then we can conclude that $\Lambda$ is trivial (Corollary \ref{triviallocalsystem}). This criterion applies in many cases of interest, e.g. whenever the contact structure admits a strong symplectic filling.

\item \textit{Sign-ambiguity}.
Even when the double cover $\Lambda$ of $\mathcal{C}(Y , \xi_0 )$ is trivial, there is no canonical choice in the $2$-element set $\Lambda(\xi_0 )$. In fact, Lin--Ruberman--Saveliev \cite{lin-ruberman-saveliev} have shown that one \textit{cannot} associate canonically an element in $\Lambda (\xi_0 )$ to each isotopy class of a contact structure $\xi_0$ in such a way that the contact invariant $\mathbf{c}(\xi_0)$ is natural with respect to orientation-preserving diffeomorphisms of $Y$. This is done by showing that the unique tight contact structure on $-\Sigma(2,3,7)$ admits a contactomorphism which reverses the sign of $\mathbf{c}(\xi_0 )$. We also note that the local system $\Lambda$ is trivial, because this contact structure is strongly (and in fact Stein) fillable.

\item \textit{Invariance}. The construction of (\ref{invtmap2}) involved choices. The main ones were the base contact structure $\xi_0$ together with a base triple $(\xi_0 , \alpha_0 , j_0 ) \in \mathcal{CM}(Y , \xi_0 )$ and an admissible perturbation $\mathfrak{q} \in \mathcal{P}$. The remaining ones were rather inessential choices of cutoff functions (\S \ref{AKend} , \S \ref{moduliandpert}). Given two choices $(\xi_i , \alpha_i , j_i ) \in \mathcal{CM}(Y , \xi_0 )$, $i = 0 , 1$, together with perturbations and cutoff functions that we omit from the notation, we obtain two corresponding chain maps $$\psi_i : S_{\ast}^{\pi_i} (\mathcal{C};\Lambda_R) \rightarrow \widehat{C}^{\ast}(Y , \s_{\xi_i , \alpha_i , j_i };R), \quad i = 0 , 1 .$$ If we choose a generic path from $(\xi_0 , \alpha_0 , j_0 )$ to $(\xi_1 , \alpha_1 , j_1 )$ this yields a spin-c structure $\mathfrak{s}_W$ on $W = [0,1] \times Y$ and after a further choice of perturbations there is an associated cobordism map $\widehat{m}(W, \mathfrak{s}_W) : \widehat{C}^{\ast}(Y , \s_{\xi_1 , \alpha_1 , j_1 };R) \rightarrow \widehat{C}^{\ast}(Y , \s_{\xi_0 , \alpha_0 , j_0 };R)$ (see \cite{KM}, p.518). We also define a subcomplex $S_{\ast} \subset S^{\pi_i}_{\ast}(\mathcal{C}; \Lambda_R)$ of chains transverse to both $\pi_0$ and $\pi_1$ (in the same index range as before). The inclusion of this subcomplex is a quasi-isomorphism. Then one concludes by showing that the following diagram is homotopy-commutative, which is a standard argument

\begin{tikzcd}
& S^{\pi_0}_{\ast}(\mathcal{C}; \Lambda_R ) \arrow{r}{\psi_0} & \widehat{C}^{\ast}(Y , \s_{\xi_0 , \alpha_0 , j_0 };R)\\
S_{\ast} \arrow{ur} \arrow{dr} & & \\
& S^{\pi_1}_{\ast}(\mathcal{C};\Lambda_R ) \arrow{r}{\psi_1} & \widehat{C}^{\ast}(Y , \s_{\xi_1 , \alpha_1 , j_1 } ; R) \arrow{uu}{\widehat{m}(W,\mathfrak{s}_W )}\\
\end{tikzcd}

\item \textit{Naturality}. The assertion on naturality from Theorem \ref{mainthm} (see the Remark after the aforementioned Theorem) readily follows from the construction from this section. 
\end{enumerate}
\end{Remark}

\section{Module structures}\label{modulesection}

In this section we define the module structures that Theorem \ref{mainthm} (B) refers to. We consider the graded ring 
\begin{align*}
\mathbb{A}(Y ; \Z) = \Z[U] \otimes_{\Z} \Lambda^{\ast} \big( H_{1}(Y ; \Z )/ \mathrm{torsion} \big) \\
|U| = 2 , \quad |\gamma| = 1  \quad \gamma \in H_1 (Y ; \Z ) /\mathrm{torsion}. 
\end{align*}
We write $\mathbb{A}^{\dagger}(Y ; \Z )$ for the opposite ring, with the opposite grading: $|U| = -2$, $|\gamma| = -1$ for $\gamma \in H_1 (Y ; \Z ) / \mathrm{torsion}$. For a given (commutative, unital) ring $R$, we obtain graded $R$-algebras $\mathbb{A}(Y ; R ) := \mathbb{A} (Y ; \Z ) \otimes R$ and $\mathbb{A}^{\dagger}(Y ; R ) := \mathbb{A}^{\dagger}(Y ; \Z ) \otimes R $.

\begin{Remark}
A different notation was used earlier, namely $\mathbb{A}(R) = \mathbb{A}^{\dagger}(Y ; R )$ (see (\ref{algebra})).
\end{Remark}

The Floer cohomology groups $\widehat{HM}^{\ast}(Y , \s ; R)$ carry a natural module structure over the graded $R$-algebra $\mathbb{A}(Y ; R )$ \cite{KM}. In this section, we first give a chain level description of this module structure which is well-suited to our purposes. We make no claim of originality here, as the material presented here is surely known to the experts. Our approach is "dual" to that of [\cite{KM}, \S VII], and in a similar spirit to the construction of the $U$ map given in [\cite{monolens}, \S 4.11].  Finally, we introduce the analogous $\mathbb{A}^{\dagger}(Y ; R )$-module structure on $H_{\ast}(\mathcal{C}(Y , \xi_0 ) ; \Lambda_R )$. The geometric interpretation of these algebraic structures that we provide in this section will be a key ingredient in the proof of Theorem \ref{mainthm} (B).

\subsection{The module structure on $\widehat{HM}^{\ast}(Y , \s  )$ }

Throughout this subsection, we fix a closed oriented $3$-manifold $Y$, and a spin-c structure $\s = (g , S , \rho)$ on $Y$. A construction reminiscent of the cup product pairing on the cohomology of $\mathcal{B}^{\sigma}(Y, \s )$ yields a pairing 
\begin{align}
H^{k}(\mathcal{B}^{\sigma}(Y , \s ) ; R) \otimes \widehat{HM}^{\ast}(Y , \s ; R) \xrightarrow{\cup} \widehat{HM}^{\ast+k}(Y , \s ; R ).\label{pairing}
\end{align}
The $\mathbb{A}(Y ; R )$-module structure on $\widehat{HM}^{\ast}(Y , \s;R )$ is then obtained from a canonical isomorphism $\mathbb{A}(Y ; R ) \cong H^{\ast}(\mathcal{B}^{\sigma}(Y , \s ) ; R )$. In what follows, our goal is to first describe this isomorphism (Proposition \ref{cohring2}) and later define the pairing (\ref{pairing}).

\subsubsection{The cohomology ring of configuration space}

We consider the blown-up configuration space $\mathcal{B}^{\sigma}(Y , \s )$ as in \S \ref{confspaceY}, where we have dropped the $k-1/2$ subscript for ease in notation. Its homotopy type is that of $\C P^\infty \times T$, where $T$ is a torus of dimension $b_1 (Y) = \mathrm{rank} \, H_1 (Y ; \Z )$. This fact is proved in [\cite{KM}, \S 9.7]. Because we will use it later, we present here a short argument (in the same spirit) that proves a weaker statement.

\begin{Proposition}[\cite{KM}]\label{cohring1}
There is an isomorphism of graded algebras $$H^{\ast}(\mathcal{B}^{\sigma} (Y, \s ); \Z) \cong \mathbb{A}(Y ; \Z ).$$
\end{Proposition}

\begin{Remark}The isomorphism given in the proof below is not canonical. We will obtain a canonical isomorphism in Proposition \ref{cohring2} using a different approach.
\end{Remark}

\begin{proof}
The inclusion of $\mathcal{B}^\ast (Y , \s )$ into the blown-up configuration space $\mathcal{B}^{\sigma}(Y , \s )$ induces a homotopy-equivalence, so we work with the former. We fix a spin-c connection $B_0$ on $S$. For another spin-c connection $B$ we have the Hodge decomposition $B- B_0 = h + d \alpha + d^{\ast}\beta$ where $h$ is harmonic. The projection $(B , \Psi )\mapsto h + d^\ast \beta$ induces a well-defined fibre bundle projection 
\begin{align}
\mathcal{B}^{\ast}(Y , \s ) \rightarrow \big\{ a \in \Omega^1 (Y ; i \R ) \, | \, d^{\ast}a = 0 \big\} / \mathcal{G}^{h}(Y) \label{fibrebundle}
\end{align} 

Here $\mathcal{G}^{h} (Y)$ stands for the group of harmonic maps $Y \rightarrow \mathrm{U}(1)$. The fiber of (\ref{fibrebundle}) is given by the projectivisation of the complex vector space of $L^{2}_{k-1/2}$ sections of $S$, which has the weak homotopy-type of $\mathbb{C}P^{\infty}.$ By further projecting to the harmonic part, we obtain a homotopy equivalence of the base of (\ref{fibrebundle}) with the torus of harmonic $1$-forms $\mathcal{H}^{1}(Y ; i \R ) / \mathcal{H}^{1}(Y ; 2 \pi i \Z ) $, which is diffeomorphic to a torus $T$ of rank $b_1 (Y)$.

We next argue that this fibre bundle is cohomologically trivial, which completes the proof. The space $\mathcal{B}^{\ast}(Y , \s )$ is the base of a principal $\mathcal{G}(Y)$-bundle with weakly contractible total space $\mathcal{C}^\ast (Y , \s ) \simeq  \ast$ and so $\mathcal{B}^{\ast}(Y , \s )$ is a model for the classifying space $B \mathcal{G}(Y)$. The inclusion of the fibre agrees with the map on classifying spaces induced by the inclusion map $\mathrm{U}(1) \rightarrow \mathcal{G}(Y)$ by the constant gauge transformations. Now, the inclusion followed by evaluation at a fixed point $\mathrm{U}(1) \rightarrow \mathcal{G}(Y) \rightarrow \mathrm{U}(1)$ has degree $1$, which shows that the inclusion of the fibre induces a surjective map on cohomology. This shows that the bundle is cohomologically trivial by the Theorem of Leray-Hirsch.
\end{proof}

\subsubsection{The slant product construction}

A standard construction [\cite{DK}, \S 5] involving the \textit{slant product}
\begin{align*}
\backslash : H_{k}(C_{\ast}  ) \otimes H^{n}((C_{\ast})^{\check{}} \otimes B^{\ast})  & \rightarrow H^{n-k}(B^{\ast}  )\\
\alpha \otimes c & \mapsto \alpha \backslash c  
\end{align*}
can be used to produce cohomology classes on $\mathcal{B}^\sigma (Y , \s )$ from homology classes in $Y$, by taking the slant product with characteristic classes of bundles over $Y \times \mathcal{B}^{\sigma}(Y , \s )$. We now describe this construction adapted to our setting.

\begin{Definition}
The \textit{canonical line bundle} $\mathcal{U}$ over $ Y \times \mathcal{B}^\sigma (Y , \s )$ is constructed from the trivial complex line bundle $\C \times Y \times \mathcal{C}^\sigma (Y , \s )$ over $Y \times \mathcal{C}^{\sigma}(Y , \s )$ as follows: make this vector bundle into a $\mathcal{G}(Y)$-equivariant vector bundle by acting on the base in the standard way and on the total space by $v \cdot (\lambda , p , B , s , \Psi ) := (v(p) \lambda , p , B - v^{-1} dv , s , v  \Psi )$, where $v \in \mathcal{G}(Y)$. The bundle $\mathcal{U}$ is obtained by taking the quotient by the $\mathcal{G}(Y)$-action.
\end{Definition}

\begin{Definition}
The \textit{slant product map} is defined for $k = 0 , 1$ by
\begin{align}
\mu : H_{k}(Y ; \Z )/ \mathrm{torsion} & \rightarrow H^{2-k}(\mathcal{B}^{\sigma}(Y , \s ); \Z )\label{slant1} \\
\alpha & \mapsto \alpha \backslash c_1 (\mathcal{U} ) \nonumber
\end{align}
\end{Definition}
\begin{Remark}
Observe that the torsion in $H_{k}(Y ; \Z )$ is not in play in (\ref{slant1}) because the cohomology of $\mathcal{B}^{\sigma}(Y , \s )$ has no torsion (Proposition \ref{cohring1}).
\end{Remark}

Recall from \S \ref{familiesirred} that there is a \textit{universal} family of spin-c structures and irreducible configurations on $Y$ parametrised by $\mathcal{B}^{\sigma}(Y , \s )$. We denote by $\mathbb{S} := \underline{S_{\infty}} \rightarrow   Y \times \mathcal{B}^{\sigma}(Y , \s )$ the universal family of spinor bundles, which arises from the quotient by the natural action of $\mathcal{G}(Y)$ on the fibres and base of the bundle $\mathrm{pr}_{1}^{\ast} S \rightarrow Y \times \mathcal{C}^{\sigma}(Y , \s )$. We denote by $L := \Lambda^2 S \rightarrow Y$ the line bundle associated to the spin-c structure $\s$ on $Y$. Of course, there is the splitting $\mathbb{S} = \mathbb{S}^+ \oplus \mathbb{S}^-$, and we observe that

\begin{Lemma}\label{tautbundles}
There is an isomorphism $ \Lambda^2 \mathbb{S}^+ \cong (\mathrm{pr}_{1}^{\ast}L) \otimes \mathcal{U}^{\otimes 2} $ of $\mathrm{U}(1)$-bundles over $  Y \times \mathcal{B}^{\sigma}(Y)$.
\end{Lemma}
\begin{proof}
We have 
\begin{align*}
\mathbb{S}^+ = ( \mathrm{pr}_{1}^\ast S^+ )/\mathcal{G}(Y) \cong \big(\mathrm{pr}_{1}^{\ast} (S^+ \otimes \mathbb{C} )\big)/\mathcal{G}(Y) \cong (\mathrm{pr}_{1}^\ast S^+) \otimes \big( ( \mathrm{pr}_{1}^\ast \mathbb{C}) /\mathcal{G}(Y) \big) = (\mathrm{pr}_{1}^\ast S^+ ) \otimes \mathcal{U}
\end{align*}
and thus $\Lambda^2 \mathbb{S}^+ \cong (\mathrm{pr}_{1}^\ast \Lambda^2 S^+ ) \otimes \mathcal{U}^{\otimes 2} = (\mathrm{pr}_{1}^\ast L )\otimes \mathcal{U}^{\otimes 2}$.
\end{proof}
Thus, since $\mathrm{pr}_{1}^{\ast} L$ is pulled back from $Y$, one could have defined $\mu$ in terms of the bundle $\mathbb{S}^+$ instead, as
\begin{align}
\mu (\alpha ) = \frac{1}{2}\big( \alpha \backslash  c_1 (\mathbb{S}^+ ) \big)  = \frac{1}{2} \big( \alpha \backslash c_1 (\Lambda^2 \mathbb{S}^+ )  \big). \label{slant12}
\end{align}

Below we provide geometric interpretations for the maps $\mu : H_k ( Y ; \Z ) / \mathrm{torsion} \rightarrow H^{2-k} ( \mathcal{B}^{\sigma}(Y ,\s ) ; \Z )$ for $k = 0 , 1$. The ultimate goal in doing so is to describe the image of $\mu$ from a dual point of view.

\subsubsection{The case $k = 0$}\label{geomdescription1}

From (\ref{slant12}) it is clear that $\mu(1)$ agrees with the first Chern class of the restriction of the bundle $\mathcal{U}$ to a slice $p \times \mathcal{B}^{\sigma}(Y , \s )$ : $\mu (1 ) = c_1 (  \mathcal{U}|_{p \times\mathcal{B}^{\sigma}(Y , \s )} ) .$
This class can be understood from a dual point of view, which we do now. 

\begin{Remark}
Below, any inclusion $M \subset B$ of manifolds $M$ and $B$ with boundary is assumed to provide an inclusion of the boundaries $\partial M \subset \partial B$ as well. 
\end{Remark}

\begin{Definition}
Let $B$ be a Banach manifold with boundary, and let $Z \subset B$ be a Banach submanifold with boundary $\partial Z \subset \partial B$ and which is of finite codimension and cooriented. We say that $Z$ is \textit{Poincar\'{e} dual} to a given cohomology class $c \in H^{q}(B ; \Z )$ if for any finite-dimensional compact oriented submanifold $M \subset B$ with boundary $\partial M \subset \partial B$ such that $M$ is embedded transversely to $Z$ (in the sense that $M \cap Z$, $\partial M \cap Z$ and $M \cap \partial Z$ are transverse intersections in the ambient $B$) then the oriented submanifold $M \cap Z \subset M$ is Poincar\'{e} dual to the cohomology class $c$ restricted onto $M$. Namely, $$\mathrm{PD}(c_{| M} ) = [M \cap Z , \partial (M \cap  Z ) ] \in H_{\mathrm{dim}M - q} ( M , \partial M ; \Z  ).$$
\end{Definition}

\begin{Remark}
Above, $M \cap Z$ is oriented in the standard way: by the exact sequence $0 \rightarrow TM \cap TZ \rightarrow TM \rightarrow TB / TZ \rightarrow 0. $
\end{Remark}

Going back to our case of interest, a section of $\mathcal{U}|_{p \times \mathcal{B}^{\sigma}(Y , \s ) }$ is provided by a $\mathcal{G}(Y)$-equivariant map $f : \mathcal{C}^{\sigma}(Y , \s ) \rightarrow \C$, with $v \in \mathcal{G}(Y)$ acting on $\C$ by the element $v(p ) \in \mathrm{U}(1)$. A concrete example of such map can be obtained as follows: fix a unitary trivialisation of the fibre of $S$ at the point $p \in Y$, denoted by $\tau = (\tau_1 , \tau_2 ): S_{p} \xrightarrow{\cong} \C^2$, and set $f_{\tau} (B , s, \Psi ) = \tau_1 \Psi (p ) .$ The section $f_{\tau}$ just constructed is transverse to the zero section of $\mathcal{U}|_{p \times \mathcal{B}^{\sigma}(Y , \s )  }$. We obtain:
\begin{Lemma}\label{dualU}
The oriented submanifold with boundary $Z_{\tau} := f_{\tau}^{-1}(0 ) \subset \mathcal{B}^{\sigma}(Y , \s )$ is Poincar\'{e} dual to $\mu (1) \in H^{2}(\mathcal{B}^{\sigma}(Y , \s ) ; \Z)$.
\end{Lemma}




\subsubsection{The case $k = 1$} The dual interpretation of the classes $\mu ([\gamma] ) \in H^1 ( \mathcal{B}^{\sigma}(Y , \s ) ; \Z )$ for a homology class $[\gamma] \in H_1 ( Y ; \Z )$ brings in the holonomy of $\mathrm{U}(1)$-connections, as follows. The universal family of spinor bundles $\mathbb{S} \rightarrow  Y \times \mathcal{B}^{\sigma}(Y ; \Z )$ carries a tautological family of connections along the $Y$-slices. For a given closed curve $\gamma : S^1 \rightarrow Y $, we obtain a \textit{half-holonomy evaluation map}  
\begin{align*}
\mathcal{B}^{\sigma}(Y , \s ) & \xrightarrow{h_{\gamma} } \mathrm{U}(1) \\
[ B , s, \Psi ] & \mapsto \mathrm{exp} \Big( \frac{1}{2} \int_{\gamma} \hat{B} \Big)  .
\end{align*}

We now explain this formula. As before, $\hat{B}$ stands for the $\mathrm{U}(1)$ connection induced by $B$ on $L= \Lambda^2 S$. By the integral above we mean the following: choose a trivialisation of $S$ along the closed curve $\gamma$, so as to identify $\hat{B}$ with a $1$-form $b$ on $\gamma$ with values in $i \R$, and evaluate $\mathrm{exp} \frac{1}{2} \int_{\gamma} b$. This element of $\mathrm{U}(1)$ does not depend on the chosen trivialisation of $S$ along $\gamma$ and the chosen representative of the gauge-equivalence class, since for different choices the $1$-form $b$ changes by adding $- 2 v^{-1}dv$ for some smooth function $v : S^1 \rightarrow \mathrm{U}(1)$ (and note that $\int_{S^1} v^{-1} dv \in 2 \pi i \mathbb{Z}$).

Thus, the half-holonomy evaluation map $h_\gamma$ provides a square root of the "usual" holonomy map $ \mathcal{B}^\sigma (Y, \mathfrak{s} ) \rightarrow \mathrm{U}(1) $ given by $[B] \mapsto \mathrm{exp} \int_\gamma \hat{B}$. The geometric content of the slant map for $k = 1$ is described by the next result:

\begin{Proposition}\label{geom2}
Let $\gamma$ be an oriented closed curve in $Y$. The class $ \mu ( [\gamma] ) \in H^1 ( \mathcal{B}^{\sigma}(Y , \s ) ; \Z )$ is represented by the half-holonomy map $h_\gamma : \mathcal{B}^{\sigma}(Y , \s ) \rightarrow \mathrm{U}(1)$. Thus, $\mu (\gamma )$ is Poincar\'{e} dual to the fibres of the submersion $h_\gamma  .$
\end{Proposition}

To show this, we consider a hermitian line bundle $L$ over a finite-dimensional manifold $X$. Denote by $\mathcal{A}$ the affine space of unitary connections on $L$, and by $\mathcal{G} = C^\infty (X,\mathrm{U}(1))$ the gauge group of $L$. As before, there is a tautological unitary line bundle $\mathcal{L}$ over $(\mathcal{A}/\mathcal{G} ) \times X$ carrying a tautological family of unitary connections on the $X$-slices.

\begin{Lemma}\label{lemmaholonomy}
For each $\gamma \in H_1 (X ; \Z )$, the class $ \gamma \backslash c_1 ( \mathcal{L} )  \in H^1 ( \mathcal{A}/\mathcal{G} ; \Z )$ is the cohomology class represented by the holonomy map $\mathrm{hol}_\gamma :  \mathcal{A}/\mathcal{G} \rightarrow \mathrm{U}(1)$ given by $\mathrm{hol}_\gamma ([B] )= \mathrm{exp} \int_\gamma B $.
\end{Lemma}

\begin{proof}
We view $\mathrm{U}(1)$ as $i \R / 2 \pi i \Z$, and denote  by $\omega = [ \frac{1}{2\pi} dx ] \in H^1 (\mathrm{U}(1) ; \Z )$ the fundamental cohomology class. We must establish the identity $\mathrm{hol}^{\ast}_{\gamma} \omega = \gamma \backslash c_1 ( \mathcal{L} )$, which is equivalent to the following: for any integral $1$-cycle $\delta$ in $\mathcal{A}/\mathcal{G}$ we have 
\begin{align}
\langle \omega , ( \mathrm{hol}_{\gamma} )_{\ast} \delta \rangle = \langle c_1 (\mathcal{L}) ,   \gamma \times \delta \rangle. \label{holonomyidentity}
\end{align}
That it suffices to show (\ref{holonomyidentity}) follows from the fact $\mathcal{A}/\mathcal{G}$ has no torsion in its cohomology.

We may suppose that $\delta $ is a smooth map $\delta : S^1 \rightarrow \mathcal{A}/\mathcal{G}$. This can be viewed as a path $t \mapsto B(t)$ of unitary connections on $L$ with $B(0)$ gauge-equivalent to $B(1)$. We see that 
$$\langle \omega , ( \mathrm{hol}_{\gamma} )_{\ast} \delta \rangle = \frac{1}{2 \pi i} \int_{t = 0}^{1} \Big( \int_{\gamma} \frac{\partial B(t)}{\partial t} \Big) dt  = -i \int_{\gamma \times \delta} \frac{\partial B(t)}{\partial t}  \wedge \omega . $$

We now provide a representative for the class $c_1( \mathcal{L} )|_{\gamma \times \delta}.$ The bundle $\mathcal{L}|_{\gamma \times \delta}$ carries the family of connections $B(t)$ on the $\gamma$-slices, and these induce a well-defined connection $\mathbf{B}$ on $\mathcal{L}|_{\gamma \times \delta}$ by setting $\nabla_{\mathbf{B}} = \frac{d}{dt} + \nabla_{B(t)}$. The class $c_1( \mathcal{L} )|_{\gamma \times \delta}$ is represented by the Chern-Weil form $$\frac{i}{2 \pi }F_{\mathbf{B}} = \frac{i}{2 \pi} \Big( F_{B(t)} + dt \wedge \frac{\partial B(t)}{\partial t} \Big) $$
and hence 
$$\langle c_1 ( \mathcal{L}  ) , \gamma \times \delta \rangle = i \int_{\gamma \times \delta }\omega \wedge \frac{ \partial B(t)}{\partial t}  = -i \int_{\gamma \times \delta } \frac{ \partial B(t)}{\partial t} \wedge \omega .$$\end{proof}

\begin{proof}[Proof of Proposition \ref{geom2}]
We specialise our discussion from above to $X := Y$ and $L := \Lambda^2 S^+$. We have a natural map $p: \mathcal{B}^\sigma (Y, \s ) \rightarrow \mathcal{A}/\mathcal{G}$ given by $p ( [B,s,\Psi]  ) = [\hat{B}]$, and the identity $\mathrm{hol}_\gamma \circ p = (h_\gamma )^2$. It follows from this and Lemma \ref{lemmaholonomy} that
\[
2 \cdot h_{\gamma}^\ast [\mathrm{U}(1)]^\vee = p^\ast \mathrm{hol}_{\gamma}^\ast [\mathrm{U}(1)]^\vee = p^\ast ( \gamma \backslash c_1 (\mathcal{L}  ) ) = \gamma \backslash (\mathrm{id}_Y \times p )^\ast c_1 (\mathcal{L} ).
\]
From Lemma \ref{tautbundles} we have $\Lambda^2 \mathbb{S}^+ \cong (\mathrm{pr}_{1}^\ast \Lambda^2 S^+) \otimes \mathcal{U}^{\otimes 2}$, and likewise one has $\mathcal{L} \cong (\mathrm{pr}_{1}^\ast L) \otimes \mathcal{U}$ where $\mathcal{U} = (\mathbb{C} \times Y \times \mathcal{A} )/\mathcal{G}$ is the canonical line bundle over $Y \times \mathcal{A}/\mathcal{G}$. Because the fiberwise $\mathrm{U}(1)$-action on $S^+$ induces the \textit{weight two} $\mathrm{U}(1)$-action on $L = \Lambda^2 S^+$ it follows that $(\mathrm{id}_Y \times p )^\ast \mathcal{U} \cong \mathcal{U}^2$, and hence 
\[
(\mathrm{id}_Y \times p )^\ast \mathcal{L} \cong (\mathrm{id}_Y \times p )^\ast \big( (\mathrm{pr}_{1}^\ast L ) \otimes \mathcal{U}  \big) = (\mathrm{pr}_{1}^\ast L ) \otimes \mathcal{U}^{\otimes 2} \cong \Lambda^2 \mathbb{S}^+ .
\]
Putting everything together we have $2 \cdot h_{\gamma}^\ast [\mathrm{U}(1)]^\vee = \gamma \backslash c_1 (\Lambda^2 \mathbb{S}^+ ) = 2 \mu (\gamma ) $. This gives $h_{\gamma}^\ast [\mathrm{U}(1)]^\vee = \mu (\gamma )$ and from this identity the result follows.
\end{proof}

\subsubsection{The cohomology ring of the configuration space, again}

We can now upgrade the isomorphism in Proposition \ref{cohring1} to a canonical one:

\begin{Proposition}\label{cohring2}
The slant map $\mu$ induces an isomorphism of graded rings 
\begin{align*}
 \mathbb{A}(Y ; \Z ) & \xrightarrow{\cong} H^{\ast} ( \mathcal{B}^{\sigma} (Y , \s ) ; \Z ) 
 \end{align*}
determined by sending $U  \mapsto \mu (1)$, and $[\gamma] \mapsto \mu ([\gamma] )$ for $[\gamma] \in H_1 ( Y ; \Z )/\mathrm{torsion}$.

\end{Proposition}

\begin{proof}
We consider the fibre bundle (\ref{fibrebundle}) from the proof of Proposition \ref{cohring1}. Its fibre has the weak homotopy-type of $\C P^{\infty}$, and the line bundle $\mathcal{U}|_{p \times \mathcal{B}^{\ast}(Y , \s )  }$ restricts to the canonical line bundle $\mathcal{O}(1)$ over $\C P^{\infty}.$ Hence the class $\mu (1 ) = c_1 ( \mathcal{U}|_{p \times \mathcal{B}^{\ast}(Y , \s ) } )\in H^2 (\mathcal{B}^{\ast} (Y , \s ) ; \Z )$ restricts to a generator of the cohomology ring of the fibres.

On the other hand, the base of the fibre bundle has the homotopy type of the torus $\mathcal{H}^1 ( Y ; i \R ) / \mathcal{H}^1 (Y ; 2 \pi i \Z )$. Choosing a $\Z$-basis of oriented closed curves $(\gamma_i )_{i = 1 , \ldots , b_1 (Y)}$ for $H_1 ( Y ; \Z ) / \mathrm{torsion}$ we obtain an explicit identification with the torus $T = \mathrm{U}(1)^{\times b_1(Y)}$  
\begin{align*}
\mathcal{H}^1 ( Y ; i \R ) / \mathcal{H}^1 (Y ; 2 \pi i \Z )  \xrightarrow{\cong} T  \, ,  \, \quad [b] \mapsto  \Big(  \mathrm{exp}  \int_{\gamma_i} b  \Big)_{i = 1, \ldots , b_1 (Y) }
\end{align*}
and the bundle projection $\mathcal{B}^{\ast}(Y , \s ) \rightarrow T$ is then identified with $$[B , \Psi ] \mapsto \Big( \mathrm{exp}  \int_{\gamma_i} (B - B_0 )^{\mathcal{H}}  \Big)_{i = 1, \ldots , b_1 (Y)}.$$ The latter map is easily seen to be homotopic to the product of the half-holonomy maps $\mathrm{hol}_{\gamma_i}$, and hence a basis for the cohomology of the base of the fibre bundle pulls back to the classes $\mu ( \gamma_i )$ (using Proposition \ref{geom2}). The fact that the fibre bundle is cohomologically trivial was shown in the proof of \ref{cohring1}, so the result follows.\end{proof}

\subsubsection{The module structure in Floer cohomology}\label{modulefloer}


The cup product pairing $(\ref{pairing})$ in Floer cohomology is obtained, roughly speaking, by integrating cohomology classes in $\mathcal{B}^{\sigma}(Y , \s )$ over the moduli spaces $M_{z}([\afr],[\bfr])$. A general definition using \v{C}ech cohomology is given in [\cite{KM}, \S 25]. Using our dual description of the generators of the cohomology of $\mathcal{B}^{\sigma}(Y , \s )$ we now give an equivalent description of this pairing which will serve better our purposes. We note that the construction given here for the $U$ map is essentially identical to the construction found in [\cite{monolens},\S 4.11], except that here we make use of the explicit sections $f_\tau$ of the canonical line bundle $\mathcal{U}$, rather than arbitrary ones.

After choosing a metric and admissible perturbation $(g , \mathfrak{q} )$, there is a (universal) Seiberg--Witten moduli space $\mathfrak{M}^{\prime}([\afr] , [\bfr] ) \rightarrow \mathcal{P}$ over the cylinder $(\R \times Y , dt^2 + g )$ . This is constructed in [\cite{KM} , \S 25] in the more general setting of cobordism maps, as a fibre product of moduli spaces over the cylinders $(- \infty , -1/2 ] \times Y$ , $[-1/2 , 1/2] \times Y$ and $[1/2 , + \infty ) \times Y $. Here, the moduli space over $[-1/2 , 1/2] \times Y$ consists of configurations $( A, \Phi  ,  \mathfrak{t} )$ where $\mathfrak{t} \in \mathcal{P}$ is used to construct a perturbation term supported in a collar neighbourhood of the boundary, by taking $\eta \hat{\mathfrak{t}}$, with $\eta(t)$ a bump function compactly supported in $(-1/2 , 0) \cup (0 , 1/2)$. 

By the unique continuation principle, the moduli $\mathfrak{M}^{\prime}([\afr] , [\bfr] )$ can be regarded as a subset of the configuration space $\mathcal{B}^{\sigma}([-1/2 , 1/2] \times Y , \s ) \times \mathcal{P}$. On the latter space we have two maps defined. First, there is the section $f_{\tau}$ of the canonical line bundle $\mathcal{U}$, $f_{\tau} (A, s , \Phi ) = \tau_1 \Phi (0 , p ) \in \C$, defined by a choice of unitary splitting $\tau =(\tau_1 , \tau_2 ) : S^{+}_{(0 , p )} \xrightarrow{\cong} \C^2 $ at $(0,p) \in \R \times Y$. On the other hand, we have the half-holonomy map $h_{\gamma} (A , s , \Phi ) = \mathrm{exp} \frac{1}{2}\int_{0 \times \gamma } \hat{A} \in \mathrm{U}(1)$ obtained from an oriented closed curve $0 \times \gamma$ in the slice $0 \times Y \subset \R \times Y$.

\begin{Proposition} \label{transversalityev}
Fix oriented closed curves $\gamma_i \subset Y , i = 1, \ldots , b_{1}(Y)  $ providing a basis of $H_1 (Y ; \Z )/ \mathrm{torsion}$. Then
\begin{enumerate}[(i)]
\item the fibre product defining the moduli spaces $\mathfrak{M}^{\prime}([\afr] , [\bfr] )$ is transverse
\item $ Z_{\tau} = f_{\tau}^{-1}(0)$ is transverse to the submanifold $\mathfrak{M}^{\prime}([\afr] , [\bfr] ) \subset \mathcal{B}^{\sigma}([-1/2,1/2] \times Y , \s ) \times \mathcal{P}$
\item for each $i$, $Z_{\gamma_i , \kappa} = h_{\gamma_i}^{-1}(\kappa)$ is transverse to the submanifold $\mathfrak{M}^{\prime}([\afr] , [\bfr] ) \subset \mathcal{B}^{\sigma}([-1/2,1/2] \times Y , \s ) \times \mathcal{P}$, where $\kappa \in \mathrm{U}(1)$ is any given value.
\end{enumerate}
\end{Proposition}

Part (i) is proved in [\cite{KM}, \S 25] in a more general setting, and (ii)-(iii) follow in a similar way as the transversality results presented in \S \ref{transversality}. To define the module structure on the monopole Floer cohomology group $\widehat{HM}^{\ast}(Y , \s; \Z )$ one chooses a perturbation $\mathfrak{t} \in \mathcal{P}$ that is a regular value of the Fredholm maps 
\begin{align*}
Z_{\tau} \cap \mathfrak{M}^{\prime}([\afr] , [\bfr] )  \rightarrow \mathcal{P}\\
Z_{\gamma_i , \kappa} \cap \mathfrak{M}^{\prime}([\afr] , [\bfr] )  \rightarrow \mathcal{P}
\end{align*}
for all $i$ and all pairs of critical points $[\afr] , [\bfr]$, and we denote by $M([\afr] , U , [\bfr]  ; \tau )$ and $M([\afr] , \gamma_i , [\bfr] ; \kappa )$ the corresponding fibres over $\mathfrak{t}$, which are smooth manifolds of finite dimension.


The $\mathbb{A}(Y ; \Z )$-module structure on $\widehat{HM}^{\ast}(Y , \s ; \Z)$ is now constructed by writing down maps $\widehat{m}_\tau (U )^{\ast} , \widehat{m}_\kappa (\gamma_i  )^{\ast} : \widehat{C}^{\ast}(Y , \s ) \rightarrow \widehat{C}^{\ast}(Y , \s )$ as follows. Each enumerates trajectories between joining critical points of the three kinds, e.g.
\begin{align*}
m(U)^{u}_{o} : C^{u} \rightarrow C^{o} \, , \quad [\afr] \mapsto \sum_{[\bfr] \in \mathfrak{C}^{o}} \# M([\afr] , U , [\bfr] ; \tau ) \cdot [\bfr]
\end{align*}
and similarly for maps $m(U)^{o}_{s} , m(U)^{o}_{o}  , m(U)^{s}_{u}, m(U)^{u}_{s}$ together with similar maps $\overline{m}(U)^{s}_{u} , \overline{m}(U)^{u}_{u} , \overline{m}(U)^{s}_{u}$ for the reducible loci in the moduli spaces. These assemble into a chain map $\widehat{m}_\tau (U) : \widehat{C}_{\ast}(Y , \s ; \Z) \rightarrow \widehat{C}_{\ast}(Y , \s ; \Z)$ given by 
\begin{align*}
 \begin{pmatrix} m(U)^{o}_{o} & m(U)^{u}_{o} \\ \overline{m}(U)^{s}_{u} \partial^{o}_{s} - \overline{\partial}^{s}_{u} m(U)^{o}_{s} & \overline{m}(U)^{u}_{u} + \overline{m}(U)^{s}_{u} \partial^{u}_{s} - \overline{\partial}^{s}_{u} m(U)^{u}_{s} \end{pmatrix} 
 \end{align*}
and dualising yields the desired cochain map $\widehat{m}(U; \tau)^{\ast}$. Similarly one obtains the cochain map $\widehat{m}_\kappa (\gamma_i  , \kappa)^{\ast}$. Up to chain homotopy, these maps do not depend on $\tau$ or $\kappa$, and passing to cohomology defines the action of $U , \gamma_i \in \mathbb{A}(Y , \Z )$ on $\widehat{HM}^{\ast}(Y , \s ; \Z )$, which gives the pairing (\ref{pairing}) when $R = \Z$. For a general ring $R$, we tensor the cochain maps $\widehat{m}_\tau (U)^\ast , \widehat{m}_\kappa (\gamma_i)^\ast $ with $R$, and this induces the action of $\mathbb{A}(Y ; R ) = \mathbb{A}(Y ; \Z ) \otimes R $ on the monopole Floer cohomology $\widehat{HM}^{\ast}(Y , \s ; R )$. This completes our description of the module structure (\ref{pairing}) in monopole Floer cohomology.

\begin{Remark}
We briefly mention why the above description of the module structure on Floer homology agrees with the \v{C}ech cohomology construction given in [\cite{KM},\S 25]. As in [\cite{KM},\S 21] we consider a compact $d$-dimensional space $N^d$ stratified by oriented smooth manifolds, where we denote by $M^d \subset N^d$ the top dimensional stratum, and let $\mathfrak{U}$ be an open cover of $N^d$ transverse to the strata (in the sense of [\cite{KM}, \S 21.2]). Given a connected component $M_{\alpha}^{d} \subset M^d$ of the top dimensional stratum there is an associated integration map (see \cite{KM}, formula (21.3) ) $\langle - , [M_{\alpha}^{d}] \rangle: \widecheck{C}^d (\mathfrak{U} ; \mathbb{Z} ) \rightarrow \mathbb{Z}$, where $\widecheck{C}^d (\mathfrak{U} ; \mathbb{Z} )$ is the abelian group of degree $d$ \v{C}ech cochains (which are automatically closed and vanish when restricted to $C^{d-1}(\mathfrak{U}|_{N^{d-1}}; \mathbb{Z} )$ because of the transversality assumption). The essential point is the following: if $u \in \widecheck{C}^d (\mathfrak{U} ; \mathbb{Z} )$ is a degree $d$ cochain representing a cohomology class in $H^{d}_{c}(M^d ; \mathbb{Z} ) \cong H^d (N^d , N^{d-1} ; \mathbb{Z} )$ which is Poincaré dual to an oriented submanifold $Z \subset M$ (necessarily of dimension zero) which is transverse to the strata of $N^d$, then one has that $\langle u, [M_{\alpha}^{d}] \rangle$ agrees with the signed count of points in $Z \cap M_{\alpha}^d$. Applying this observation when $N^d$ is the compactification by broken trajectories of a $d$-dimensional Seiberg-Witten moduli space $M ([\afr] , [\bfr] )$ on the cylinder $\mathbb{R} \times Y$ (here we mean the "smaller" version of the compactification, as in Definition 24.6.9 in \cite{KM}) it follows that our definition of the module structure agrees with the one given in \cite{KM}.
\end{Remark}

\subsection{The module structure on $H_{\ast}(\mathcal{C}(Y , \xi_0 ) ; \Lambda_R )$}

We now fix a closed oriented contact $3$-manifold $(Y , \xi_0 )$. We will define a graded $\mathbb{A}^{\dagger}(Y ; R )$-module structure  $$\mathbb{A}^{\dagger}(Y ; R ) \otimes H_{\ast}(\mathcal{C}(Y , \xi_0 ); \Lambda_R ) \rightarrow H_{\ast-k}(\mathcal{C}(Y , \xi_0 ) ; \Lambda_R)$$ and describe its geometric meaning.


\subsubsection{The slant construction}

We do a similar construction as before, using the slant product $$\:  H_{k}(Y; \Z) \otimes H^{n}(Y \times \mathcal{C}(Y, \xi_0 ); \Z ) \rightarrow H^{n-k} (\mathcal{C}(Y , \xi_0 );\Z ).$$ There is a tautological family of contact structures on $Y$ parametrised by $\mathcal{C}(Y , \xi_0 )$, which provide us with a real oriented rank $2$ vector bundle $\boldsymbol{\xi} \rightarrow Y \times \mathcal{C}(Y , \xi_0 )$. The bundle $\boldsymbol{\xi}$ is a subbundle of a trivial rank $3$ bundle (since $TY$ is trivial for any closed oriented $3$-manifold), so its second Stiefel-Whitney class $w_2 ( \boldsymbol{\xi} )$ vanishes. Consequently, the Euler class $e(\boldsymbol{\xi} ) \in H^2 (Y \times \mathcal{C}(Y , \xi_0 ); \Z )$ is divisible by $2$.

\begin{Definition} \label{slantcont}For $k = 0 , 1$ we define
\begin{align}
\overline{\mu} :  H_{k} ( Y ; \Z )/\mathrm{torsion} & \rightarrow H^{2-k} ( \mathcal{C}(Y , \xi_0 ) ; \Z ) \label{slant2prime}\\
\alpha & \mapsto \frac{1}{2}\alpha \backslash  e( \boldsymbol{\xi} )\nonumber .
\end{align}
\end{Definition}
\begin{Remark}
Observe that for $k = 0$, the slant product map (\ref{slant2prime}) is, a priori, only well-defined as a map into $H^{2}(\mathcal{C}(Y , \xi_0 ) ; \Z )$ \textit{modulo} the $2$-torsion subgroup. This ambiguity arises from dividing by $2$ in (\ref{slant2prime}). However, we now explain that there is a canonical lift, which we take as the definition of (\ref{slant2prime}). Observe that taking $\alpha = 1 \in H_0 ( Y ; \Z ) = \Z$ we have $\alpha \backslash e ( \boldsymbol{\xi} ) = e ( \boldsymbol{\xi}|_{p \times \mathcal{C}(Y , \xi_0 )} )$, so the matter reduces to having a preferred square root of $\boldsymbol{\xi}|_{p \times \mathcal{C}(Y , \xi_0 )}$, up to isomorphism. This rank $2$ bundle comes with a preferred homotopy class of embeddings into the trivial rank $3$ bundle, simply obtained by fixing a positive framing $T_p Y \cong \R^3$. In other words, there is a canonical homotopy class of maps $\mathcal{C}(Y , \xi_0 ) \rightarrow \ogr{\R^3}$ into the Grassmannian of oriented $2$-planes in $\R^3$, which by pullling back the tautological $2$-plane bundle over $\ogr{\R^3}$ yield $\boldsymbol{\xi}|_{p \times \mathcal{C}(Y , \xi_0 )}$. It is now elementary to observe that there is a unique square root (i.e. spin structure) for the tautological $2$-plane bundle over $\ogr{\R^3}$.
\end{Remark}

\begin{Definition}\label{modulecontact}
We endow $H_{\ast}( \mathcal{C}(Y , \xi_0 ) ; \Lambda_R )$ with a graded $\mathbb{A}^{\dagger}(Y ; R)$-module structure
\begin{align}
\mathbb{A}^{\dagger}(Y ; R ) \otimes H_{\ast}( \mathcal{C}(Y , \xi_0 ) ; \Lambda_R ) \rightarrow H_{\ast}( \mathcal{C}(Y , \xi_0 ) ; \Lambda_R )\label{pairing2}
\end{align}
by setting: for $T \in H_{n}(\mathcal{C}(Y , \xi_0 ); \Lambda_R )$ 
\begin{align*}
U \cdot T & := \overline{\mu}(1) \cap T \in H_{n - 2 }(\mathcal{C}(Y , \xi_0 ) ; \Lambda_R ) \\
\gamma \cdot T &:= \overline{\mu} ( \gamma ) \cap T  \in H_{n - 1 }(\mathcal{C}(Y , \xi_0 ) ; \Lambda_R )\, , \quad \gamma \in H_1 (Y ; \Z ) / \mathrm{torsion} .
\end{align*}
\end{Definition}
Here $\cap$ denotes the cap product with \textit{coefficients in the local system} $\Lambda_R$: $$H^{k}(\mathcal{C}(Y , \xi_0 ) ; \Z ) \otimes H_{n}( \mathcal{C}(Y , \xi_0 ) ; \Lambda_R) \rightarrow H_{n-k}(\mathcal{C}(Y , \xi_0 ) ; \Lambda_R) .$$

We now relate the slant product maps $\mu$ and $\overline{\mu}$. The space $\mathcal{CM}(Y , \xi_0 )$ of triples $(\xi , \alpha , j )$ parametrises a family of spin-c structures and irreducible configurations on $Y$ (see \S \ref{spinccont}), and to this it corresponds a classifying map $f : \mathcal{CM}(Y , \xi_0 ) \rightarrow \mathcal{B}^{\ast} ( Y , \s_{\xi_0 , \alpha_0 , j_0} )$ (see Lemma \ref{correspondence}). Here, $(\xi_0 , \alpha_0 , j_0 )$ is a fixed triple.
\begin{Lemma}\label{pull}
We have the identity $\overline{\mu} = f^{\ast} \mu$. 
\end{Lemma}

\begin{proof}
Indeed, under the map $\mathrm{id}_Y \times f : Y \times \mathcal{C}(Y , \xi_0 ) \rightarrow Y \times \mathcal{B}^{\ast}(Y , \s _{\xi_0 , \alpha_0 , j_0 })$ the bundle $\Lambda^2 \mathbb{S}^+$ pulls back to the bundle $\boldsymbol{\xi}$, and hence $e(\boldsymbol{\xi} ) = (\mathrm{id}_Y \times f)^\ast c_1 ( \Lambda^2 \mathbb{S}^+ )$.
\end{proof}

\subsubsection{Geometric interpretations}

We conclude this section by interpreting the module action (\ref{pairing}) in geometric terms. We start with the $U$ map $$H_{\ast}(\mathcal{C}(Y , \xi_0 ); \Lambda_R ) \xrightarrow{U} H_{\ast-2}(\mathcal{C}(Y , \xi_0 ; \Lambda_R )\, , \quad T \mapsto \overline{\mu}(1) \cap T .$$ Fixing a point $p \in Y$, there is a natural evaluation map to the Grassmanian of oriented planes in $T_{p} Y$ 
\begin{align*}
\mathcal{C}(Y , \xi_0 ) \xrightarrow{ev} \ogr{ T_{p} Y} \cong S^2 \, , \quad \xi \mapsto \xi (p ) .
\end{align*}
The main geometric content is:
\begin{Proposition}\label{geominterpcont}
The class $\overline{\mu}(1) \in H^{2}( \mathcal{C}(Y , \xi_0 ) ; \Z )$ is represented by the map $ev$, i.e. $\overline{\mu}(1) = ev ^{\ast} [S^2]^{\vee}$.
\end{Proposition}

From this it follows that the more geometric description of the $U$ action on $H_{\ast}(\mathcal{C}(Y , \xi_0 ) )$ given in (\ref{Uintro}) agrees with the one just given in Definition \ref{modulecontact}.


\begin{proof}

For each unitary framing $(\tau_1 , \tau_2 ) : S_p \xrightarrow{\cong} \C^2$ of the fibre over $p \in Y$ of the spinor bundle $S := S_{\xi_0 , \alpha_0 , j_0}$ we have the section $f_{\tau}(B , \Psi ) = \tau_1 \Psi ( p ) $ of the canonical line bundle $\mathcal{U}$ restricted to $p \times \mathcal{B}^{\ast}(Y , \s_{\xi_0 , \alpha_0 , j_0} )$. This ``pencil" of sections induces a map $e$ to the projectivisation of $S_{p}$, away from the base locus $B = \{ [B ,\Psi ] \, : \, \Psi (p ) = 0 \}$
\begin{align*}
\mathcal{B}^{\ast}(Y , \s_{\xi_0 , \alpha_0 , j_0 } ) \setminus B \xrightarrow{e} \mathbb{P}(S_p ) \cong \mathbb{P}^1 \, , \quad [ B ,  \Psi ] \mapsto \C \cdot \Psi ( p ) . 
\end{align*}
Note that if $Z_{\tau} := f_{\tau}^{-1}(0)$ then $Z_\tau \setminus B   \subset \mathcal{B}^{\ast}(Y , \s_{\xi_0 , \alpha_0 , j_0 } ) \setminus B$ is a fibre of $e$. 

Observe now that the image of the classifying map $f : \mathcal{CM}(Y , \xi_0 ) \rightarrow \mathcal{B}^{\ast}(Y , \s_{\xi_0 , \alpha_0 , j_0 } )$ does not meet $B$. From this, together with the fact that $Z_\tau$ is Poincaré dual to $\mu (1)$ (Lemma \ref{dualU}) and $f^\ast \mu = \overline{\mu}$ (Lemma \ref{pull}), it follows that the cohomology class $\overline{\mu}(1) \in H^{2}(\mathcal{CM}(Y , \xi_0 ))$ is represented by the map $e \circ f: \mathcal{CM}(Y , \xi_0 ) \rightarrow \mathbb{P}(S_p ).$

The key observation is now the following

\begin{Lemma}\label{cliffordlemma}
Let $(V, g)$ be a $3$-dimensional real oriented inner product vector space, and let $S \cong \C^2$ be its fundamental spin-c representation.
Then there exists a canonical diffeomorphism $\ogr{V} \cong \mathbb{P}(S).$
\end{Lemma}

\begin{proof}
The Grassmanian $\ogr{V}$ is diffeomorphic to the unit sphere in $V^{\ast}$ via $$\mathrm{Sph}(V^{\ast} ,g ) \xrightarrow{\cong} \ogr{V}\, , \quad \alpha \mapsto \ker \alpha .$$ 
Recall that the fundamental spin-c representation provides an isomorphism $V^{\ast} \cong \mathfrak{su}(S)$ as $\mathrm{Spin}^{\C}(3)$ modules. Under this isomorphism, the Clifford multiplication,  a given $\alpha \in \mathrm{Sph}(V^{\ast} , g )$ acts on $S$ decomposing it into $\pm i$ eigenspaces $S = l^{+} \oplus l^{-}$.  Each $l^{\pm}$ is a complex line in $S$, and the assignment $$\mathrm{Sph}(V^{\ast} , g ) \rightarrow \mathbb{P}(S)\, , \quad \alpha \mapsto l^{+}$$ provides a diffeomophism, concluding the proof.
\end{proof}

To conclude, apply Lemma \ref{cliffordlemma} for each $t \in \mathcal{CM}(Y, \xi_0 )$ to the inner product spaces $(T_p Y , g_{\xi_t , \alpha_t , j_t} )$ and spin-c structures $(g_{\xi_t , \alpha_t , j_t } , S_{\xi_0 , \alpha_0 , j_0 } , \rho_{\xi_0 , \alpha_0 , j_0}\circ b^{\ast}_{g , g_{\xi_0 , \alpha_0 , j_0}}  )$. This provides an isomorphism between the (trivial) bundles 
\[
\widetilde{\mathrm{Gr}}_2 (T_p Y ) \times \mathcal{CM}(Y, \xi_0 ) \cong \mathbb{P}(S_p ) \times \mathcal{CM}(Y, \xi_0 )
\]
over $\mathcal{CM}(Y, \xi_0 )$. Under this isomorphism, the evaluation map $ev$ becomes identified with the map $e \circ f$ (both maps are regarded here as sections of the corresponding trivial bundles), and from this the result now follows.
 \end{proof}


Finally, we briefly comment on the action of $\gamma \in H_1 (Y , \Z )$ $$H_{\ast}(\mathcal{C}(Y , \xi_0 ); \Lambda_R) \xrightarrow{\gamma} H_{\ast-1}(\mathcal{C}(Y , \xi_0 ); \Lambda_R).$$ The geometric interpretation that we will need in the subsequent sections is already provided by Lemma \ref{lemmaholonomy}: upon choosing a reduction of the structure group of $\boldsymbol{\xi} \rightarrow Y \times \mathcal{C}(Y , \xi_0 )$ to $\mathrm{U}(1)$, and a family of unitary connections $\{ B_{\xi} \}$ over the $Y$-slices, one obtains a holonomy map $$\mathcal{C}(Y , \xi_0 ) \rightarrow \mathrm{U}(1) \, , \quad \xi \mapsto \mathrm{exp} \int_{\gamma} B_{\xi} $$whose regular fibres are Poincaré dual to $2 \overline{\mu} ( \gamma ) \in H^{1}( \mathcal{C}(Y , \xi_0 ) ; \Z )$. In particular, the canonical spin-c connections $\hat{B}_{\xi , \alpha , j}$ on $Y$ parametrised by $(\xi, \alpha , j ) \in \mathcal{CM}(Y , \xi_0 ) \simeq \mathcal{C}(Y , \xi_0 )$ provide such a family of connections.

\section{The neck-stretching argument}\label{proofmainsection}

In this section we establish Theorem \ref{mainthm} (B). This asserts that the families contact invariant $\fc : H_{\ast}( \mathcal{C}(Y , \xi_0 ); \Lambda_R  ) \rightarrow \widehat{HM}^{\ast}(Y , \s_{\xi_0} ; R)$ intertwines the module structures, which were introduced in \S \ref{modulesection}. We must show: for $T \in H_{\ast}(\mathcal{C}(Y , \xi_0 ) ; \Lambda_R )$ and a homology class $\gamma \in H_1 ( Y ; \Z )$
\begin{align}
U \cdot \fc(T) &= \fc(U \cdot T )\label{ident1} \\
\gamma \cdot \fc(T) & = \fc(\gamma \cdot T).\label{ident2}
\end{align}
We sketch now the main ideas in the case of $U$. The key is to consider, for a given simplex $\sigma:\Delta^n \rightarrow \mathcal{CM}(Y , \xi_0 )$, a moduli space $\mathcal{M}([\afr] , U  , \sigma ; \tau) \rightarrow \Delta^n \times \R \ni (t , s) $ of solutions to the Seiberg--Witten equations on $Z^+$ that meet certain evaluation constraint at the point $(s , p ) \in \R \times Y \cong Z^+$. Here $p \in Y$ is fixed, whereas $s \in \R$ is not, and hence the evaluation constraint is thought of as travelling through $Z^+$ from the cylindrical to the symplectic end. The evaluation constraint itself is that the spinor $\Phi$ satisfies $\tau_1 \Phi = 0$ at the point $(s,p)$, for a suitably chosen trivialisation $\tau = (\tau_1 , \tau_2 )$ of the bundle $S^+$ along the line $\R \times p \subset Z^+$. Such moduli spaces will be referred to as \textit{parametrised evaluation moduli spaces}. The main part of the argument is to analyse the ends of the (non-compact) moduli $\mathcal{M}([\afr] , U , \sigma ; \tau )$ as $s \rightarrow \pm \infty$. As $s\rightarrow - \infty$ we will see that the solutions to the equations degenerate into \textit{broken} configurations, which in the simplest case consist of pairs of configurations $(\gamma_1 , \gamma_0 )$, the first of which solves the Seiberg--Witten equations over an infinite cylinder $\R \times Y$ with an evaluation constraint, and the second is an unconstrained solution over $Z^+$. The interesting part of the moduli space, however, shows up as $s \rightarrow + \infty$. Here we will 
show that $\mathcal{M}([\afr] , U , \sigma ; \tau)$ looks like the product $\R \times M$ where $M$ is, in a sense, the intersection of the moduli $M([\afr] , \sigma )$ over the simplex (constructed in \S \ref{theinvariant}) with a fibre of the map $\mathcal{B}^{\ast}(Y , \s_{\xi_0 , \alpha_0 , j_0} ) \dashrightarrow \mathbb{P}^1$ from \S \ref{modulesection}.

This will allow us to construct compactifications of the parametrised evaluation moduli spaces, and the identities (\ref{ident1})-(\ref{ident2}) arise from counting the boundary points of the compactified $1$-dimensional parametrised evaluation moduli.

\subsection{Parametrised evaluation moduli spaces over $Z^+$}\label{universalUmoduli}

\subsubsection{A family of perturbations}\label{modulioverR}

As a starting point for the construction of the parametrised evaluation moduli spaces we introduce an intermediate moduli space 
\begin{align}
\mathcal{M}([\afr] , Z^+ ) \rightarrow \mathcal{CM}(Y, \xi_0 ) \times \mathcal{P} \times \R \label{pertmodulimain}
\end{align}
analogous to $\mathfrak{M}([\afr] , Z^+ )$ from \S \ref{moduliandpert}. The only new feature is that the $\R$ factor in the base will parametrise various perturbations of the Seiberg--Witten equations. The parametrised evaluation moduli space will result from imposing constraints on the configurations in $\mathcal{M}([\afr] , Z^+ )$. Following the same scheme as in \S \ref{moduliandpert}, we construct $\mathcal{M}([\afr] , Z^+)$ as a fibre product of moduli.

The first step is constructing a moduli space 
\begin{align}
\mathcal{M}_{k} \rightarrow \mathcal{CM}(Y, \xi_0 ) \times \mathcal{P} \times \R \label{pertmoduli}
\end{align}
 in the same flavour of $\mathfrak{M}_{k}(K^{\prime}) \rightarrow \mathcal{CM}(Y, \xi_0 ) \times \mathcal{P}$. The moduli space $\mathcal{M}_k$ consists of gauge-equivalence classes of quintuples $(A, \Phi , t , \mathfrak{p} , s )$, where the variables $(t , \mathfrak{p} , s ) \in \mathcal{CM}(Y, \xi_0 )\times \mathcal{P} \times \R$ provide the map in (\ref{pertmoduli}), and $(A, \Phi )$ are configurations over the region 
 \begin{align}
 K(s) = [ \mathrm{m}(s) , + \infty ) \times Y \subset Z^+ 
 \end{align}
Here $\mathrm{m}(s) $ stands for the function $\mathrm{min}(s-1 , 0 )$, or rather, a suitable smooth approximation of it. Such $(A, \Phi , t , \mathfrak{p} , s )$ must be asymptotic to canonical configurations as before, and satisfy the Seiberg--Witten equations $$\mathrm{sw}(A , \Phi , u ) + \lambda( A , \Phi ,  t , \mathfrak{p} , s ) = 0$$ perturbed by a certain quantity $\lambda ( A, \Phi , t , \mathfrak{p} ,  s ) \in \Upupsilon_{k-1}$ which we now describe. It is given by the section
 \begin{align}
& \lambda : \mathcal{C}_{k }(K(s))\times \mathcal{P} \times \R \rightarrow \Upupsilon_{k-1 } \nonumber \\
& ( A, \Phi , t , \mathfrak{p} , s)  \mapsto \varphi^{1}_{s} \hat{\mathfrak{q}}(A, \Phi) + \varphi^{2}_{s} \hat{\mathfrak{p}}(A, \Phi) + \eta_{s} \hat{\mathfrak{t}} (A, \Phi) + \varphi^{3} \hat{\mathfrak{p}}_{K,u}. 
\end{align}

Here, $\mathfrak{q}$, $\mathfrak{t}$ are fixed generic perturbations chosen as in \S \ref{moduliandpert} and \S \ref{modulefloer}.  Also, $\varphi_{s}^1$, $\varphi_{s}^2$, $\eta_{s}$ are fixed $\R$-families of non-negative functions on $\R$, and $\varphi^3$ is the function we chose in \S \ref{moduliandpert}. We require that they relate to the functions $\varphi^{1} , \varphi^2$ of \S \ref{moduliandpert}, and $\eta$ of \S \ref{modulefloer} as follows. Let $(\tau_s f)(t) := f(t+s)$. Then
\begin{enumerate}[(i)]
\item 
$\varphi_{s}^1= \tau_{-\mathrm{m}(s)} \varphi^1$ for all $s \in \R$

\item $\varphi_{s}^2 = \tau_{- \mathrm{m}(s)} \varphi^2$  for all $ s \in \R$


\item $\eta_s = \tau_{- s} \eta$ for $s < 0$ very negative, and identically vanishing for $s \geq 0$.


\end{enumerate} 
 
The choice of such perturbation data will ultimately ensure the behaviour of the parametrised evaluation moduli spaces that we have described at the beginning of \S \ref{proofmainsection}.

We want to make $\mathcal{M}_k$ into a Banach manifold. By applying the $\R$-family of translations $t \mapsto t - s$ we can view $\mathcal{M}_k$ as a subset of a suitable configuration space $\mathcal{B}_k = \mathcal{C}_k / \mathcal{G}_{k+1}$ over $[ 0 , + \infty ) \times Y$. As before, the latter is a $C^{l-k-2}$ Banach manifold and $\mathcal{M}_k$ is the transverse zero set of a section of a bundle over $\mathcal{B}_k$ given by the perturbed Seiberg--Witten map; hence a Banach manifold of the same regularity. The claimed transversality follows, once more, from the results in \S \ref{transversality}.

We then have restriction maps 
\begin{align*}
\mathcal{M}_k & \xrightarrow{\mathcal{R}_{-}} \mathcal{B}^{\sigma}_{k-1/2}(Y , \s_{\xi_0 , \alpha_0 , j_0} )\\
M_{k}([\afr] , (- \infty , 0 ] \times Y )& \xrightarrow{R_{+}} \mathcal{B}^{\sigma}_{k-1/2}(Y , \s_{\xi_0 , \alpha_0 , j_0} )
\end{align*}
onto the left-most and right-most end, respectively. From \S \ref{transversality} it will follow that the fibre product 
\begin{align*}
\mathcal{M}([\afr] , Z^+) = \mathrm{Fib}(R_{+} , \mathcal{R}_{-} )
\end{align*}
is transverse, and that the projection to $\mathcal{CM}(Y, \xi_0 ) \times \mathcal{P} \times \R$ is Fredholm. This completes the construction of (\ref{pertmodulimain}).

\subsubsection{The parametrised $U$-moduli space}\label{paramUmod}
We fix a point $p \in Y$ throughout. We denote by $\tau$ an arbitrary unitary splitting of the fibre of the spinor bundle $S_{\xi_0 , \alpha_0 , j_0} \rightarrow Y$ over the point $ p \in Y$, i.e. a unitary isomorphism $\tau = (\tau_1 , \tau_2 ): (S_{\xi_0 , \alpha_0 , j_0})_p \xrightarrow{\cong} \C^2$. Given such $\tau$, which we may view as an element in the unitary group $\mathrm{U}(2)$, we obtain an extension to a unitary splitting of the positive spinor bundle $S^+ \rightarrow Z^+$ as follows. First, over the cylindrical end $Z = (- \infty , 0] \times Y \subset Z^+$ by translation. For the symplectic end $K = [1 , + \infty ) \times Y \subset Z^+$ we proceed as follows. Recall that in \S \ref{translationinvariance} that we introduced a rescaling operator $\mathcal{R}_0$, which upon acting on canonical configurations yields them translation-invariant in some gauge. We have the translation-invariant bundle over $K$ given by $$\overline{S^+} = \mathcal{R}_{0}^{\ast} S^{+} = \C \oplus \Lambda^{0,2}_{\overline{J}_{0}}T^{\ast}K.$$ Usinf the identification $(S_{\xi_0 , \alpha_0 , j_0})_p = S^{+}_{(1,p)} \cong \mathcal{R}_{0}^{\ast} S^{+}_{(1,p)}$ we may simply translate the splitting $\tau$ along $K$. In the transition region $[0 , 1] \times Y \subset Z^+$ we extend $\tau$ in an arbitrary manner.

We have the \textit{canonical line bundle} $\mathcal{U} \rightarrow \mathcal{B}_k  \times \R$, which arises from the $\R$-family of representations of the group of gauge transformations $\mathcal{G}_{k+1} \rightarrow \mathrm{U}(1)$ given by $v \mapsto v(s , p )$, where $s$ varies within $\R$. We pullback this bundle over to $\mathcal{M}([\afr] , Z^+ )$, which can be identified naturally as a Banach submanifold of $\mathcal{B}(Z^+) \times \R$. We consider the section of this pullback bundle given by $$f_{\tau} ( A, \Phi , t , \mathfrak{p} , s ) = \tau_1 \overline{\Phi} (s , p )$$ where $\overline{\Phi} = \mathcal{R}_{0}^{\ast} \Phi$ is the rescaled version of $\Phi$. Note that we only defined $\mathcal{R}_0$ over the region $K$; we extend it here over the whole $Z^+$ as the identity over the cylindrical end $Z$.

The following will follow from \S \ref{transversality}:
\begin{Proposition}
The section $f_\tau$ is transverse to the zero section of $\mathcal{U} \rightarrow \mathcal{M}([\afr] , Z^+ )$.
\end{Proposition}

\begin{Definition}
The (universal) \textit{parametrised $U$-moduli space} is the Banach submanifold $$\mathcal{M}([\afr] , U , Z^+ ; \tau) \subset \mathcal{M}([\afr] , Z^+ )$$ given by the zero set of the section $f_{\tau}$.
\end{Definition}
\begin{Remark} Allowing for arbitrary splittings $\tau \in \mathrm{U}(2)$ might seem strange at this point. The main case to have in mind is the basic splitting $S_{\xi_0 , \alpha_0 , j_0} = \C \oplus \xi_0 $, which over the symplectic end corresponds to the splitting $S^{+} = \C \oplus \Lambda_{J_{0}}^{0,2} T^{\ast} K.$  The section of $\mathcal{U}$ that we would want to take in this case is simply given by projecting $\Phi (s , p )$ to the trivial $\C$ factor. However, it will soon become apparent that, in order for the ends of the relevant moduli spaces in the neck-stretching argument to have a nice structure, we have to pass to a generic splitting $\tau$.
\end{Remark}

\subsubsection{The parametrised $\gamma$-moduli space}

In a similar fashion, we fix a smooth oriented closed curve $\gamma \subset Y$ and consider the map 
\begin{align}
h_{\gamma} : \mathcal{M}([\afr] , Z^+ ) \rightarrow \mathrm{U}(1) \label{paramholonomy}
\end{align}
obtained by associating to $(A, \Phi , t , \mathfrak{p} , s )$ the half-holonomy of the induced connection $\hat{A}$ on $\Lambda^2 S^+$ around the loop $s \times \gamma \subset Z^+$ $$h_{\gamma} (A , \Phi , t , \mathfrak{p} , s ) = \mathrm{exp} \frac{1}{2} \int_{s \times \gamma } \hat{A} .$$

In \S \ref{transversality} we show:

\begin{Proposition}
The map (\ref{paramholonomy}) is a submersion.
\end{Proposition}

\begin{Definition}
The (universal) \textit{parametrised} $\gamma$-\textit{moduli space} is the Banach submanifold $$\mathcal{M}([\afr] , \gamma , Z^+ ; \kappa) \subset \mathcal{M}([\afr] , Z^+ )$$ given by the preimage of $\kappa \in \mathrm{U}(1)$ under (\ref{paramholonomy}). 
\end{Definition}

\subsection{Compactifications}\label{compactifications}

\subsubsection{The setup}\label{setup}

We first introduce the moduli spaces that will be the main players in the neck-stretching argument that will follow. These are associated to a singular chain $\sigma : \Delta^n \rightarrow \mathcal{C} := \mathcal{CM}(Y , \xi_0 ) \times \mathcal{P}$ equipped with a unitary splitting $\tau \in \mathrm{U}(2)$ and a value $\kappa \in \mathrm{U}(1)$.

First, by taking the fibre product of $\mathcal{M}_{z}([\afr] , Z^+ ) \xrightarrow{\pi} \mathcal{C} \times \R$ and $\sigma \times \mathrm{id}_{\R} : \Delta^n \times \R \rightarrow \mathcal{C} \times \R$ we obtain the space $\mathcal{M}_{z}([\afr] , \sigma)$ which is a $C^2$ manifold with corners provided the fibre product is transverse. Similarly, taking the fibre product of the each of the two maps $\mathcal{M}_z ([\afr] , U , Z^+ ) , \mathcal{M}_z ([\afr] , \gamma , Z^+ ) \xrightarrow{\pi} \mathcal{C} \times \R$ with $\sigma \times \mathrm{id}_{\R}$ we obtain $C^2$ manifolds with corners $\mathcal{M}_z ([\afr] , U , \sigma  ; \tau)$, $\mathcal{M}_z ([\afr] , \gamma , \sigma ; \kappa )$ if transversality holds. In both cases the required transversality can be achieved by a $C^2$ perturbation of $\sigma$ whenever the index of $\pi$ is $\leq 1- n$, by the Thom-Smale transversality theorem (see \S \ref{transversechains}).

The task that we take up for the remainder of this section is to analyze the ends of the $1$-dimensional non-compact moduli spaces $\mathcal{M}_z ([\afr] , U , \sigma  ; \tau)$, $\mathcal{M}_z ([\afr] , \gamma , \sigma ; \kappa )$ and construct suitable compactifications of them with a nice boundary structure.

\subsubsection{Exponential decay}


Consider a configuration $(A , \Phi , t , s) $ in the moduli $\mathcal{M}_{z}([\afr] , \sigma)$. Over the symplectic end $K \subset Z^+$ the positive spinor bundle $S^+ \rightarrow Z^+$ decomposes into the $\pm 2i$ eigenspaces of the Clifford action of $\omega_t$. The canonical spinor $\Phi_{t}$ provides a framing of the $-2i$ eigenspace, and we decompose $\Phi$ accordingly
\begin{align}
\Phi &= \alpha \Phi_t + \beta \label{spinordecompos}
\end{align}
where $\alpha$ is a function, and $\beta$ is a section of the $+2i$ eigenspace. Similarly, using the canonical connection $A_{t}$ we obtain a decomposition 
\begin{align}
A = A_{t} + a \label{connectiondecompos}
\end{align}
for an $i \R$-valued $1$-form $a$. We regard $a$ as a unitary connection $\nabla_a = d + a$ on the trivial line bundle, with curvature given by $F_a = da$. There is a also a unitary connection $\tilde{\nabla}_A$ on the $+ 2i$ eigenspace $E_+ (t)$ obtained from $A$ by orthogonal projection.

The main ingredient for the various compactness results needed in this article is the following \textit{exponential decay estimate}, which follows from the work of Kronheimer--Mrowka \cite{monocont} and Zhang \cite{boyu}.

\begin{Theorem}\label{expdecay1}
There exists constants $C , \epsilon > 0$ depending on $\sigma$, with the following significance: if $(A, \Phi , t , s) \in \mathcal{M}_{z}([\afr] , \sigma)$ for some $[\afr] , z$, we have the following estimate over $K \subset Z^+$
$$ |1- |\alpha|^2  + |\beta|^2 |^2 + |\beta|^2 + | \nabla_a \alpha |^2 + | \tilde{\nabla}_A \beta |^2 + |F_a |^2  \leq C e^{-\epsilon s }.$$
\end{Theorem}

\begin{Corollary}\label{expdecay2}
For any element in $\mathcal{M}_z ([\afr], \sigma )$ there is a gauge representative $(A , \Phi , t , s )$ of it such that $A - A_t$ and $\Phi - \Phi_t$ decay exponentially over $K$ with first derivatives (with constants $C, \epsilon >0$ only depending on $\sigma$). 
\end{Corollary}
\begin{proof}
The only part which doesn't follow directly from Theorem \ref{expdecay1} is that $|A - A_t|^2 \leq Ce^{-\epsilon s }$. This is proved exactly as in Corollary 3.16 of \cite{monocont}.
\end{proof}

\subsubsection{The boundary of $\mathcal{M}_{z}([\afr] , U , \sigma ; \tau)$ at $s = + \infty$}\label{compactificationU}

We now describe the behaviour of configurations $(A , \Phi , t , s) \in \mathcal{M}_{z}([\afr] , U , \sigma ; \tau )$ when $s$ approaches $+ \infty$.


Denote by $e_{\infty}: \Delta^n \rightarrow \mathbb{P} (S_{(1,p)}^{+}) $ the map that associates to $t \in \Delta^n$ the fibre over the point $(1,p) \in Z^+$ of the $-2i$-eigenspace for $\rho_t (\omega_t )$, namely the line $$\C \cdot \Phi_{t}(1, p ) \subset S^{+}_{(1,p)}.$$ We encountered this map in the proof of Proposition \ref{geominterpcont}. Recall that $\tau$ provides a translation-invariant unitary splitting $\tau = (\tau_1 , \tau_2 ) : \overline{S^+}|_{[1, + \infty) \times p} \cong \C^2$ as in \S\ref{paramUmod}. This provides us with a preferred line $l_\tau \in \mathbb{P} (S_{p}^+ )$, namely that line which corresponds with $0:1$ under the identification $\mathbb{P}(S^{+}_{p}) \cong \mathbb{P}( \mathbb{C}^2 )$ given by $\tau$. 

\begin{Definition}
The $U$-\textit{limiting locus} at $s = + \infty$ of $\sigma$ is the subset $Z_{\infty , \tau}(\sigma) := e_{\infty}^{-1}(l_\tau) \subset \Delta^n .$
\end{Definition}
The limiting set at infinity is a compact subset of $\Delta^n$. Later we will require that $l_\tau$ is a regular value (by varying $\tau$), so that $Z_{\infty , \tau}(\sigma)$ will be a submanifold (with corners) of $\Delta^n$. The terminology we chose is justified by the next observation:

\begin{Lemma} \label{locuslemma}
Suppose $(A_n , \Phi_n , t_n , s_n ) \in \mathcal{M}_{z}([\afr] , U , \sigma ; \tau)$ is a sequence of configurations such that $\lim_{n \rightarrow +\infty} s_n  = + \infty$ and $\lim_{n \rightarrow +\infty} t_n  = t^{\ast} $ for some $t^{\ast} \in \Delta^n$. Then $t^{\ast}$ lies in $Z_{\infty , \tau}(\sigma) \subset  \Delta^n$
\end{Lemma}
\begin{proof}
We choose a family of canonical configurations $(A_{t} , \Phi_{t})$ defined for $t \in \Delta^n$, since $\Delta^n$ is contractible. By Lemma \ref{spinorasymp} we may assume, after passing to a different gauge, that $\overline{\Phi_{t}}$ are translation-invariant spinors over the symplectic end $K$.

By Theorem \ref{expdecay1}, there exist constants $C > 0$ and $\epsilon > 0$ independent of $n$, such that for any $s \in \mathbb{R}$ and $y \in Y$ $$|\Phi_n (s,y) - \Phi_{t_n}(s , y) | \leq C e^{-\epsilon s}. $$ Thus
\begin{align*}
| \overline{\Phi}_n (s , y) - \overline{\Phi}_{0 , t_n}(s , y) |  \leq |\mathcal{R}_{0}^{\ast}(s)| |\Phi_n (s,y) - \Phi_{t_n}(s , y) |  \leq C  |\mathcal{R}_{0}^{\ast}(s)| e^{-\epsilon s} .
\end{align*}
where $|\mathcal{R}_{0}^{\ast}(s)| $ denotes the pointwise norm of the rescaling operator, which for $s \geq 1$ equals $1$. 

On the other hand, by the definition of $\mathcal{M}_z ([\afr] , U , \sigma ; \tau )$, at the point $(s_n , p)$ the evaluation constraint $\tau_1 \overline{\Phi}_n  (s_n , p ) = 0$ holds. By the above bound, $| \overline{\Phi}_n (s_n , p) - \overline{\Phi}_{ t_n}(s_n , p) | $ converges to zero, and hence
$\lim_{n \rightarrow \infty}\tau_1 \overline{\Phi}_{t_n}(s_n , p ) = 0.$ By translation-invariance we have $\tau_1 \overline{\Phi}_{t_n} (1, p ) = \tau_1 \overline{\Phi}_{t_n} (s_n, p ) = 0$. 
Hence we obtain $\tau_1 \overline{\Phi}_{ t^{\ast}} (1, p) = \lim_{n \rightarrow \infty}\tau_1 \overline{\Phi}_{t_n}(1, p) =0$, which means that $e_{\infty} (t^{\ast} ) = l_\tau $, as required.
\end{proof}

\begin{Definition}\label{Ulimiting}
The $U$-\textit{limiting moduli space} at $s = + \infty$ is the preimage of the $U$-limiting locus $Z_{\infty, \tau}(\sigma) = e_{\infty}^{-1}(l_\tau) \subset \Delta^n$ under the map $M_{z}([\afr] , \sigma) \rightarrow \Delta^n$. We denote it by $M_{z}([\afr] , Z_{\infty, \tau }(\sigma) )$.
\end{Definition}


The next is the main result of this section. It describes the shape of $\mathcal{M}_{z}([\afr] , U , \sigma ; \tau )$ as the evaluation constraint goes to $+ \infty$. 

\begin{Theorem}\label{structureU}
Let $\sigma$ be a $C^2$ singular chain in $\mathcal{C} =\mathcal{CM}(Y , \xi_0 ) \times \mathcal{P} $. After a $C^2$ pertubation of $\sigma$ and a residual choice of splitting $\tau \in \mathrm{U}(2)$, there exists a constant $s_0 > 0$ such that the following holds for all $[\afr] , z$ for which the moduli space $\mathcal{M}_z ([\afr] , U , \sigma ; \tau )$ has expected dimension $1$: 
\begin{itemize}
\item  the moduli spaces $\mathcal{M}_{z}([\afr] , U, \sigma ; \tau )$ are transversely cut out and the moduli spaces $M_{z}([\afr] , Z_{\infty , \tau}(\sigma ) )$ consist of a finite set of transversely cut out points
\item there is a homeomorphism of the open subset $\{ s > s_0\} \cap \mathcal{M}_{z}([\afr] , U , \sigma ; \tau )$ with the product $M_{z}([\afr], Z_{\infty , \tau}(\sigma )) \times ( s_0, + \infty )$, compatible with the projection to $(s_0 , + \infty )$.
\end{itemize}
\end{Theorem}

\begin{proof}
We start with some preliminary observations. First, note that the transversality assertion for the moduli spaces $\mathcal{M}_{z}([\afr] , U , \sigma ; \tau )$ of dimension $1$ follows by an application of the Thom-Smale transversality theorem, in the same way as for the moduli spaces $M_z ([\afr] , \sigma )$. In this case, again by standard finiteness results (see \S \ref{compactness}) we only have finitely many non-empty $\mathcal{M}_z ([\afr] , U , \sigma ; \tau )$ with dimension $1$. Also, the moduli space $M_{z}([\afr] , Z_{\infty , \tau}(\sigma ))$ is compact, since its expected dimension is $0$. Thus, if it is transversely cut out then it will consist of finitely-many points.

We choose a family of canonical configurations $(A_{t} , \Phi_{t} )$ parametrised by $t \in \Delta^n$ in translation-invariant form (see Proposition \ref{spinorasymp}, Definition \ref{transinvform}). The open subset $\{ s >0\} \cap \mathcal{M}_{z}([\afr] , \sigma)$ is canonically identified with the product $M_{z}([\afr] , \sigma ) \times ( 0 , + \infty) $, compatibly with the projection to $(0 , +\infty )$. For this product structure, the canonical line bundle $\mathcal{U} \rightarrow \mathcal{M}_{z}([\afr] , \sigma )$ is identified with a pullback to the first factor in the product. Next, we extend the section $f_\tau$ of $\mathcal{U}$ (whose zeros give $\mathcal{M}_z ([\afr] , U , \sigma ; \tau ) \subset \mathcal{M}_z ([\afr] , \sigma ; \tau )$) to a section $F$ defined over $s = + \infty$ as follows
\begin{align*}
  F : \,  &M_{z}([\afr] , \sigma) \times ( 0 , + \infty]   \rightarrow \mathcal{U} \\
& (  A , \Phi , t ,  s )  \mapsto   \tau_1 \overline{\Phi} ( s, p ), \text{  if  } s \neq + \infty \\
 &  ( A , \Phi , t ,  + \infty )  \mapsto  \tau_1 \overline{\Phi}_{t}(1, p) .
\end{align*}

We write $F_s$ for the smooth section given by restriction of $F$ to the slice $M_z ([\afr] , \sigma ) \times \{ s\}$. The required result is thus of Implicit Function Theorem type. Namely, from the version of this given in [\cite{KM}, Lemma 19.3.3] it will follow that the map $F^{-1}(0 )  \rightarrow (0 , + \infty ] $ (given by the projection $(A, \Phi , t , s) \mapsto s $) defines a \textit{topological} submersion over $(F_\infty )^{-1}(0) = M_z ([\afr] , Z_{\infty , \tau}(\sigma ) )$ (and therefore the required homeomorphism $\mathcal{M}_z ( [\afr] , U , \sigma ; \tau ) \cong M_z ([\afr] , Z_{\infty , \tau}(\sigma ) ) \times (s_0 , + \infty )$ for some large $s_0 > 0$) provided we can show
\begin{enumerate}
    \item[(i)] $F$ is continuous
    \item[(ii)] $F_s \rightarrow F_\infty$ in $C^{1}_{\mathrm{loc}}$
    \item[(iii)] $F_\infty$ is transverse to the zero section.
\end{enumerate}
Item (i) follows from the exponential decay estimates (Theorem \ref{expdecay1}). For item (ii) we proceed as follows. Recall that the configuration space $\mathcal{C}_k (K^\prime )$ over the symplectic end is a Banach manifold with tangent space
$$ T_{(A, \Phi  , t )} \mathcal{C}_k (K^\prime ) \cong \Big\{ (a , \phi , \dot{t} ) \, \, | \,\, a - \frac{\partial}{\partial \dot{t}}A_t \in L^{2}_{k}(K , g_t ) \, , \, \phi - \frac{\partial}{\partial \dot{t}}\Phi_t \in L^{2}_{k} (K , g_t ) \, , \, \dot{t} \in T_t \mathcal{CM}(Y, \xi_0 ) \Big\}$$
(see (\ref{tspace})) and a Banach space norm induced from a local chart is given by $$|| (a, \phi , \dot{t} ) || = ||a -\frac{\partial}{\partial \dot{t}}A_t ||_{L^{2}_{k}(K , g_t )} + ||\phi - \frac{\partial}{\partial \dot{t}}\Phi_t ||_{L^{2}_{k} (K , g_t )} + ||\dot{t}||.$$
The vertical components (taken with respect to the obvious connection on $\mathcal{U}$) of the derivatives of the sections $F_s$ and $F_\infty$ are
\begin{align*}
(\mathcal{D}f_s )_{(A, \Phi , t ) }(a, \phi , \dot{t} ) = \overline{\phi}(s,p) \\
(\mathcal{D}f_\infty )_{(A, \Phi , t ) }(a, \phi , \dot{t} ) = \frac{\partial}{\partial \dot{t}} \overline{\Phi}_t (1,p).
\end{align*}

We then use the continuous embedding $L^{2}_{k, \overline{g}_t , \overline{A}_t} (K , \overline{S}^+) \hookrightarrow C^{0} (K , \overline{S}^+) $ (recall that the \textit{cylindrical} metric is $\overline{g}_t = ds^2 + g_{\xi_t , \alpha_t , j_t }$ over $K$) together with the identity of Riemannian volume forms $d\mathrm{vol}_{g_t} = s^3  d \mathrm{vol}_{\overline{g}_t}$ to obtain the estimate
\begin{align*}
    s^{3/2} \cdot | (\mathcal{D} (F_s - F_\infty ))_{(A, \Phi , t )} (a , \phi , \dot{t} ) | & \leq  || s^{3/2} ( \overline{\phi} -\frac{\partial}{\partial \dot{t}} \overline{\Phi}_t ) ||_{C^{0} (K , \overline{S}^+ )} \\
    & \leq C \cdot || s^{3/2} (\overline{\phi} - \frac{\partial}{\partial \dot{t}} \overline{\Phi}_t ) ||_{L^{2}_{k, \overline{g}_t , \overline{A}_t} ( K , \overline{S}^+)}\\
    & \leq C ||\phi - \frac{\partial}{\partial \dot{t}} \Phi_t ||_{L^{2}_{k , g_t , A_t} ( K , S^+ )} .
    \end{align*}
From this we deduce that $|| F_s - F_\infty ||_{C^1 (M_z ([\afr] , \sigma) ) } \leq C / s^{3/2}$, and in particular we have $C^1_{\mathrm{loc}}$ convergence $F_s \rightarrow F_\infty$ as $s \rightarrow + \infty$ follows.

For (iii) recall that $Z_{\infty, \tau} = e_{\infty}^{-1}(l_\tau)$. The space of unitary splittings $\tau$ is identified with the unitary group $\mathrm{U}(2)$. There exists, by Sard's theorem, a residual subset of the space unitary splittings $\tau \in \mathrm{U}(2)$ for which $l_\tau \in \mathbb{P} (S_{(1,p)}^+ )$ is a regular value of $e_{\infty}$. It is straightforward to see then that 

\begin{Lemma}\label{surjectiveinfty1} If $l_\tau$ is a regular value of $e_\infty$, then the map from the configuration space
\begin{align*}
\mathcal{C}_{k}([0 , + \infty) \times Y)|_{\Delta^n} & \xrightarrow{F_{\infty}} \C\\
( A, \Phi , t ) &\mapsto \tau_1 \overline{\Phi}_{t}(1, p )
\end{align*}
has a regular value at $0$.
\end{Lemma}

In \S \ref{transversality} we establish a general transversality result for moduli spaces with evaluation constraints, which applies to certain evaluation maps that fall into an suitable class (Definition \ref{permiss}). Lemma \ref{surjectiveinfty1} shows that $F_{ \infty}$ falls into this class. This general result implies in this instance that $M_{z}([\afr] , \sigma) \xrightarrow{F_{\infty}} \mathcal{U}$ is transverse to the zero section of $\mathcal{U}$. This concludes the proof of Theorem \ref{structureU}.
\end{proof}

\subsubsection{The boundary of $\mathcal{M}_{z}([\afr] , \gamma , \sigma ; \kappa )$ at $s = + \infty$}

We now carry out an analogous study of the shape as $s$ approaches $+ \infty$ of the second kind of parametrised evaluation moduli spaces $\mathcal{M}_{z}([\afr] , \gamma , \sigma ; \kappa)$ where $\gamma \subset Y$ is a smooth oriented closed curve. First, have the analogue of Lemma \ref{locuslemma}. This time it involves the map $e_{\infty}^{\gamma}: \Delta^n \rightarrow \mathrm{U}(1)$ which associates to $t \in \Delta^n$ the half-holonomy $\mathrm{exp} \frac{1}{2} \int_{1 \times \gamma} \hat{A}_{t}$. 

\begin{Definition}
The $\gamma$-\textit{limiting locus} at $s = + \infty$ of $\sigma$ is the subset $Z_{\infty, \kappa}^{\gamma}(\sigma) =( e_{\infty}^{\gamma})^{-1}(\kappa ) \subset \Delta^n$. 
\end{Definition}

\begin{Lemma}
Suppose $(A_n , \Phi_n , t_n , s_n ) \in \mathcal{M}_{z}([\afr] , \gamma , \sigma ; \kappa)$ is a sequence of configurations such that $\lim_{n \rightarrow +\infty} s_n  = + \infty$ and $\lim_{n \rightarrow +\infty} t_n  = t^{\ast} $ for some $t^{\ast} \in \Delta^n$. Then $t^\ast$ lies in $Z_{\infty , \kappa}^{\gamma} \subset \Delta^n$.
\end{Lemma}

\begin{proof}
Let $a_n = A_n - A_{t_n}$. By Corollary \ref{expdecay2} we may assume $|a_n |^2 \leq C e^{- \epsilon s}$ over $K$. For convenience, regard $\mathrm{U}(1)$ as $i \R / 2 \pi i \Z$. There we have the identity
\begin{align*}
\frac{1}{2} \int_{s_n \times \gamma} \hat{A}_{n} - \frac{1}{2} \int_{1 \times \gamma} \hat{A}_{t_n} = \int_{s_n \times \gamma} a_n  +  \frac{1}{2} \int_{s_n \times \gamma - 1 \times \gamma} \hat{A}_{t_n}.  
\end{align*}
The second term on the right-hand side vanishes (mod $2 \pi i \mathbb{Z}$) by the translation-invariance property of the canonical connection $A_{t}$. From the exponential decay estimate on $|a_n |$ it follows that the first term goes to zero as $n\rightarrow \infty$. The result follows.
\end{proof}


\begin{Definition}
The $\gamma$-\textit{limiting moduli space} at $s = + \infty$ is the preimage of the $\gamma$-limiting locus $Z_{\infty , \kappa}^{\gamma} \subset \Delta^n$ under the map $M_{z}([\afr] , \sigma ) \rightarrow \Delta^n $. We denote it by $M_{z}([\afr] , Z_{\infty , \kappa}^{\gamma}(\sigma ))$.
\end{Definition}

We have the following analogue of Theorem \ref{structureU}, describing the shape of $\mathcal{M}_{z}([\afr] , \gamma , \sigma ; \kappa )$ as the evaluation constraint goes to $+ \infty$.

\begin{Theorem}\label{structuregamma}
Let $\sigma$ be a $C^2$ singular chain in $\mathcal{C} = \mathcal{CM}(Y , \xi_0 ) \times \mathcal{P}$. After a $C^2$ pertubation of $\sigma$ and a residual choice of $\kappa \in \mathrm{U}(1)$, there exists a constant $s_0 > 0$ such that the following holds for all $[\afr] , z$ such that $\mathcal{M}_z ([\afr] , \gamma , \sigma ; \kappa )$ has expected dimension $1$: 
\begin{itemize}
\item  all the moduli spaces $\mathcal{M}_{z}([\afr] , U, \sigma ; \kappa )$ are transversely cut out and the moduli spaces $M_{z}([\afr] , Z_{\infty , \kappa}^{\gamma}(\sigma ) )$ consist of a finite set of transversely cut out points
\item there is a homeomorphism of the open subset $\{ s > s_0\} \subset \mathcal{M}_{z}([\afr] , \gamma , \sigma ; \kappa )$ with the product $M_{z}([\afr], Z_{\infty , \kappa}^{\gamma}(\sigma )) \times ( s_0, + \infty )$, compatible with the projection to $(s_0 , + \infty )$.
\end{itemize}

\end{Theorem}

\begin{proof}
The strategy is the same as in the proof of Theorem \ref{structureU}. Rather than working with the half-holonomy map $h_{\gamma} (A , \Phi , t , s ) = \mathrm{exp} \frac{1}{2} \int_{t \times \gamma} \hat{A} $ we view $\mathrm{U}(1)$ as $i \R / 2 \pi i \Z$ and work with $f_\gamma ( A, \Phi , t , s ) = \frac{1}{2} \int_{\{s\} \times \gamma} \hat{A} $. We extend this to a map $F$ defined over $s = + \infty$ in a similar fashion as before:
\begin{align*}
F : \, &M_{z}([\afr] , \sigma ) \times ( 0 , + \infty ] \rightarrow i \R / 2 \pi i \Z \cong \mathrm{U}(1) \\
& (A , \Phi , t ,  s )\mapsto \frac{1}{2} \int_{\{s\} \times \gamma} \hat{A} \text{  if  } s \neq + \infty\\
& ( A , \Phi , t ,  + \infty ) \mapsto \frac{1}{2} \int_{\{1\} \times \gamma } \hat{A}_{t} .
\end{align*}

As in the proof of Theorem \ref{structureU} we need to show that the restrictions to the slices $F_{s}$ satisfy (i)-(iii). For (iii) we have the statement analogous to Lemma \ref{surjectiveinfty1}: for residual $\kappa \in \mathrm{U}(1)$ the map 
\begin{align*}
\mathcal{C}_k ([0 , + \infty) \times Y )|_{\Delta^n} \xrightarrow{f_{\infty}} \mathrm{U}(1)\\
(A , \Phi , t ) \mapsto \frac{1}{2} \int_{1 \times \gamma } \hat{A}_t
\end{align*}
has a regular value at $\kappa$. Indeed, $\kappa$ has this property whenever the map $e^{\gamma}_{\infty} : \Delta^n \rightarrow \mathrm{U}(1)$ has a regular value at $\kappa$. Then the general transversality results of \S \ref{transversality} imply (iii). This concludes the proof of Theorem \ref{structuregamma}.

\end{proof}

\subsubsection{The compactification of $\mathcal{M}_{z}([\afr] , U , \sigma ; \tau)$ and $\mathcal{M}_{z}([\afr] , [\afr] , \gamma , \sigma ; \kappa )$}

We are now set to describe the compactification of the parametrised evaluation moduli spaces over a simplex $\sigma$. This brings together the various moduli spaces that we have studied: $M_{z}([\afr] , \sigma )$ (\S \ref{theinvariant}), $\mathcal{M}_{z}([\afr], U , \sigma ; \tau )$ and $\mathcal{M}_{z}([\afr] , \gamma , \sigma ; \kappa )$ (\S \ref{setup}). In addition, we also have the usual moduli spaces of Floer trajectories $\breve{M}_{z}([\afr] , [\bfr] )$ (where we quotient by the reparametrization action of $\R$, as usual), and the $U$ and $\gamma$-moduli spaces over cylinders $M_{z}([\afr] , U , [\bfr] ; \tau )$, $M_{z}([\afr] , \gamma , [\bfr] ; \kappa )$ introduced in \S \ref{modulefloer}. 

First we introduce some definitions. Recall our notation $\mathcal{C} := \mathcal{CM}(Y , \xi_0 ) \times \mathcal{P}$.
\begin{Definition}
Given a unitary splitting $\tau \in U(2)$, a $C^2$ singular simplex $\sigma : \Delta^n \rightarrow \mathcal{C}$ will be called $\tau$-\textit{transverse} provided it satisfies the following transversality requirements:
\begin{enumerate}[(i)]
\item $\sigma$ is transverse to the Fredholm map $\mathfrak{M}(Z^+ ) \rightarrow \mathcal{C}$ (see \S \ref{moduliandpert}) along components with index $\leq 1- n$
\item the map $\Delta^n \times \R \ni (t , s) \mapsto (\sigma (t), s ) \in \mathcal{C} \times \R$ is transverse to the Fredholm map $\mathcal{M}(U , Z^+ ; \tau)   \rightarrow \mathcal{C} \times \R$ (see \S \ref{universalUmoduli}) along the components of index $\leq 1-n$.
\item the map $e_{\infty} : \Delta^n \rightarrow \mathbb{P}(S^{+}_{(1,p)})$ has a regular value at $l_\tau \in \mathbb{P}(S^{+}_{(1,p)})$.
\end{enumerate}

Likewise, given $\kappa \in U(1)$ a $C^2$ singular simplex $\sigma$ will be called $\kappa$-transverse if it satisfies (i) together with the following analogues of (ii) and (iii):

\begin{enumerate}
\item[(i)'] the map $\Delta^n \times \R \ni (t , s) \mapsto (\sigma (t), s ) \in \mathcal{C} \times \R$ is transverse to the Fredholm map $\mathcal{M}(\gamma , Z^+; \kappa )   \rightarrow \mathcal{C} \times \R$ along the components of index $\leq 1-n$.
\item[(ii)'] the map $e^{\gamma}_{ \infty}: \Delta^n  \rightarrow \mathrm{U}(1)$, defined in terms of $\sigma$, has a regular value at $\kappa \in U(1)$.
\end{enumerate}
\end{Definition}

By the Thom-Smale transversality Theorem, after a small perturbation one can ensure that any $\sigma$ becomes $\tau$-transverse for a residual choice of $\tau$ (likewise for $\kappa$-transverse).

The following result describes the structure of the compactifications of our moduli spaces:

\begin{Proposition}\label{strata}
Let $\sigma$ be $\tau$-transverse $C^2$ singular simplex $\sigma : \Delta^n \rightarrow \mathcal{C}$. If the moduli space $\mathcal{M}_{z}([\afr] , U , \sigma ; \tau)$ has expected dimension $0$, then it consists of finitely-many transversely cut-out points. If it has expected dimension $1$, then it is a $C^2$ manifold with boundary which admits a compactification $\mathcal{M}_{z}^{+}([\afr] , U , \sigma ; \tau)$ with the structure of a space stratified by manifolds. Its top stratum is given by $\mathcal{M}_{z}([\afr] , U , \sigma ; \tau)$ itself, and the boundary of the top stratum consists of configurations of the following types:
\begin{enumerate}[(a)]
\item $ M_{z}([\afr] , Z_{\infty, \tau}(\sigma))$
\item the moduli $\mathcal{M}_{z}([\afr] , U , \sigma|_{\Delta^{n-1}_{i}} ; \tau)$ over the codimension $1$ faces $\Delta^{n-1}_{i} \subset \Delta^n$ of $\sigma$
\item $\breve{M}_{z_1}([\afr] , [\bfr] ) \times \mathcal{M}_{z_0}([\bfr] , U , \sigma ; \tau )$
\item $\breve{M}_{z_2} ([\afr] , [\bfr] ) \times \breve{M}_{z_1}([\bfr] , [\cfr] ) \times \mathcal{M}_{z_0}([\cfr] , U , \sigma ; \tau)$
\item $M_{z_1} ([\afr] , U , [\bfr] ; \tau ) \times M_{z_0}([\bfr] , \sigma )$
\item $M_{z_2} ([\afr] , U , [\bfr]; \tau ) \times \breve{M}_{z_1}([\bfr] , [\cfr] ) \times M_{z_0}([\cfr] , \sigma  )$
\item $\breve{M}_{z_2} ([\afr] , [\bfr] ) \times M_{z_1}([\bfr] , U ,  [\cfr] ; \tau ) \times M_{z_0}([\cfr] , \sigma  )$.
\end{enumerate}
(Here, the middle factor in the triple products must be boundary-obstructed. 
The concatenation of the homotopy classes $z_i$ in every product must equal $z$.)

Furthermore, the structure near the boundary strata of type (a),(b),(c),(e) is that of a $C^0$ manifold with boundary, and the structure near (d), (f),(g) is that of a codimension $1$ $\delta$-structure (see \cite{KM}, Definition 19.5.3.)

The analogous result holds for the $\gamma$-moduli spaces.
\end{Proposition}

More generally, a compactification by broken trajectories of the moduli $\mathcal{M}_{z}([\afr] , U , \sigma ; \tau )$ of any dimension can be constructed, provided transversality holds. However, we will only use those of dimension $0$ or $1$. We refer to \S \ref{compactness} for an outline of the standard technical results that enable us to establish the compactness. We have carried out in this section the analysis of the structure of the compactification near the boundary stratum of type (a). This component of the boundary stratum is the most interesting, and will be the key to the proof of Theorem \ref{mainthm} (B). For the strata of type (c)-(g) the required gluing analysis follows similar techniques as those in \cite{KM}.

\subsection{The proof of Theorem \ref{mainthm} (B)}

We are now ready to complete the proof of Theorem \ref{mainthm} (B). This follows from chain level identities arising from enumeration of boundary points of $\mathcal{M}_{z}([\afr] , U , \sigma ; \tau )$ and $\mathcal{M}_{z}([\afr] , \gamma ; \sigma ; \kappa )$. 

\subsubsection{Orientations}\label{orientations2}

We explained in \S \ref{orientability1} how to orient the moduli $M_{z}([\afr] , \sigma )$. We now want to orient the parametrised moduli $\mathcal{M}_{z}([\afr] ,U ,  \sigma ; \tau )$ and $\mathcal{M}_{z}([\afr] , \gamma, \sigma ; \kappa )$, and for this we first orient the bigger moduli $\mathcal{M}_z ([\afr] , \sigma )$ that contains them. The latter moduli is defined as the fibre product $\mathrm{Fib}(\pi , \sigma )$ of the natural map $\pi: \mathcal{M}_{z}([\afr] , Z^+ ) \rightarrow \mathcal{C}:= \mathcal{CM}(Y, \xi_0 ) \times \mathcal{P}$ (as in (\ref{pertmodulimain}) but projecting the $\R$ factor away) and $\sigma$.

To orient $\mathcal{M}_{z}([\afr] , \sigma )$ we need to orient the determinant line $\mathrm{det} \pi$ of the Fredholm map $\pi: \mathcal{M}_{z}([\afr] , Z^+ ) \rightarrow \mathcal{C}$. Once that is done $\mathcal{M}_{z}([\afr] , \sigma )$ becomes oriented by Lemma \ref{signrule}. Since the moduli $\mathfrak{M}_{z}([\afr] , Z^+ )$ is the fibre over $s = 1$ of the natural map $$\mathcal{M}_{z}([\afr] , Z^+ )\rightarrow \R $$
then an orientation of $\mathrm{det} \pi$ is determined by an orientation of the determinant line of $\mathfrak{M}_{z}([\afr] , \sigma ) \rightarrow \mathcal{C}$ and the convention that the $\R$ factor \textit{goes first}. This is contrary to the usual fibre-first convention, but agrees with standard conventions in \cite{KM}.

The remaining moduli spaces are oriented as follows: 
\begin{itemize}
    \item $\mathcal{M}_{z}([\afr],U , \sigma ; \tau)$ is the zero set of a section of a complex line bundle over $\mathcal{M}_{z}([\afr] , \sigma )$, so we orient it as such.
\item $\mathcal{M}_{z}([\afr] , \gamma , \sigma ; \kappa )$ is the fibre of a map $\mathcal{M}_{z}([\afr] , \sigma ) \rightarrow \mathrm{U}(1)$, so we orient it using the fibre-first convention.
\item $M_{z}([\afr] , Z_{\infty , \tau}(\sigma ))$ is the fibre of a map $M_{z}([\afr] , \sigma ) \rightarrow \mathbb{P}^1$, so we orient it by the fibre-first convention.
\item $M_{z}([\afr] , U , [\bfr] ; \tau )$ and $M_{z}([\afr] , \gamma , [\bfr] ; \kappa)$ are analogous to the first two bullets.
\end{itemize}

We refer to the above as the \textit{canonical orientations} of the moduli spaces.

\subsubsection{Counting solutions to the Seiberg--Witten equations}\label{linearoperators}


If $\sigma$ is a transverse $C^2$ singular simplex in $\mathcal{C} = \mathcal{CM}(Y, \xi_0 ) \times \mathcal{P}$ (together with an element from the two-element set $\Lambda (\sigma (b ) )$, where $b \in \Delta^n$ is the barycenter, which we will omit from notation) then we have defined an associated monopole cochain in \S \ref{theinvariant}
\[
\psi (\sigma )  = \sum_{[\afr] , z} \# M_z ([\afr], \sigma ) \cdot [\afr] \in \widehat{C}^\ast (Y , \mathfrak{s}_{\xi_0 , \alpha_0 , j_0} ) .
\]
which induces the family contact invariant $\fc$. Given a $\tau$-transverse $C^2$ singular simplex $\sigma : \Delta^n \rightarrow \mathcal{C}$ for some $\tau \in U(2)$, we will now define new associated monopole cochains
\[
\theta_{ \tau}(\sigma ) \, , \psi_{\infty , \tau}(\sigma ) \in \widehat{C}^\ast (Y , \mathfrak{s}_{\xi_0 , \alpha_0 , j_0} )
\]
and if, in addition, $\sigma$ is $\kappa$-transverse for some $\kappa \in U(1)$, we will define new monopole cochains
\[ 
\theta_{\gamma , \kappa}(\sigma ) , \psi_{\infty, \gamma , \kappa }(\sigma ) \in \widehat{C}^\ast (Y , \mathfrak{s}_{\xi_0 , \alpha_0 , j_0} ).
\]
These are given by
\begin{align*}
&\theta_{\tau}( \sigma ) = \sum_{[\afr] ,z } \big( \# \mathcal{M}_{z}([\afr]  , U , \sigma ; \tau ) \big) \cdot [\afr] \\
& \theta_{\gamma , \kappa} ( \sigma  ) (\alpha )= \sum_{[\afr] , z} \big( \# \mathcal{M}_{z}([\afr] , \gamma , \sigma; \kappa ) \big) \cdot [\afr]\\
& \psi_{\infty } (\sigma ) = \sum_{[\afr] , z} \big( \# M_{z}([\afr] , Z_{\infty , \tau}(\sigma ) ) \big) \cdot [\afr]\\
& \psi_{\infty, \gamma , \kappa }^{\gamma}( \sigma ) = \sum_{[\afr] , z} \big( \# M_{z}([\afr] , Z_{\infty , \kappa}(\sigma ) ) \big) \cdot [\afr].
\end{align*}


That all these sums are indeed finite follows from standard compactness arguments as in \cite{KM} that we will review in \S \ref{compactness}. 

Similar arguments as for Proposition \ref{chainmapprop} show that $\psi_{\infty, \tau}(\sigma)$ and $\psi_{\infty, \gamma}(\sigma)$ satisfy chain map relations (up to signs): $\widehat{\partial}^{\ast} \psi_{\infty, \tau}(\sigma) =(-1)^{n} \psi_{\infty, \tau} (\partial \sigma )$ on simplices of dimension $n$, and similarly for $\psi_{\infty , \gamma , \kappa}(\sigma)$. In addition, given a $\sigma$ which is transverse for both $\tau$ and $\tau^\prime$, the difference $\psi_{\infty , \tau }(\sigma ) - \psi_{\infty , \tau^\prime} (\sigma )$ is easily seen to be \textit{exact} (likewise for two $\kappa$ and $\kappa^\prime$), so we obtain well-defined homomorphisms at the level of homology
\[ ( \psi_\infty )_\ast \, , \, (\psi_{\infty , \gamma} )_\ast : H_\ast (\mathcal{C}(Y , \xi_0 ) ; \Lambda_R ) \rightarrow \widehat{HM}^\ast (Y , \mathfrak{s}_{\xi_0 } )
\]
defined as follows: any homology class $[T]$ in $\mathcal{C}(Y , \xi_0 )$ with coefficients twisted by $\Lambda_R$ can be represented by a cycle $T$ whose simplices are $\tau$-transverse for a fixed $\tau$. Then one sets $(\psi_\infty )_\ast [T] := [\psi_{\infty , \tau}(T) ]$, and likewise for $(\psi_{\infty , \kappa})_\ast$. We have the following crucial identities:

\begin{Proposition}\label{infinitymoduli}
For any homology class $[T] \in H_{\ast} ( \mathcal{C}(Y , \xi_0 ) ; \Lambda_R )$ we have 
\begin{align*}
& (\psi_{\infty })_\ast ([T]) =  \fc( U \cdot [T] )\\ 
& (\psi_{\infty, \gamma})_\ast ([T]) =  \fc( [\gamma] \cdot [T] ).
\end{align*}
\end{Proposition}

\begin{proof}
We explain the first identity, and the second follows identically. Recall from Lemma \ref{dualU} that the cohomology class $\overline{\mu}(1) \in H^2 ( \mathcal{C}(Y , \xi_0 ) ; R )$ is Poincar\'{e} dual to the zero set of the section $f _{\tau}: \mathcal{B}^{\sigma}(Y , \s_{\xi_0 , \alpha_0 , j_0} ) \rightarrow \mathcal{U}$ restricted to $\mathcal{CM}(Y , \xi_0 ) \subset \mathcal{B}^{\sigma}(Y , \s_{\xi_0 , \alpha_0 , j_0 })$. 
For residual $\tau$ the section $f_{\tau}$ will be transverse to the zero section along $\mathcal{CM}(Y , \xi_0 )$. Any given homology class $[T]$ can be represented by a $C^2$ cycle $T$ in $\mathcal{C}(Y , \xi_0 )$ that intersects transversely the zero set $f_{\tau}^{-1}(0)$. This intersection is given by restricting $T$ to the union over the limiting loci $Z_{\infty , \tau}(\sigma)$, where $\sigma$ runs over subfaces $\sigma \subset T$ of all dimensions. This intersection can be given the structure of a cycle $T_{\infty}$, and it follows that in homology $[T_{\infty}] = \overline{\mu}(1) \cap [T] =: U \cdot [T]$. The result now follows from applying $\fc$ to both sides.
\end{proof}

\subsubsection{The chain homotopy relation}

Theorem \ref{mainthm} (B) now follows from combining Proposition \ref{infinitymoduli} and 
\begin{Proposition}\label{mainidentities}
Let $\sigma : \Delta^n \rightarrow \mathcal{C}$ be a singular $C^2$ chain which is $\tau$-transverse. Then the following chain homotopy relation (up to signs) holds
\begin{align*}
\widehat{m}_\tau (U  )^{\ast} \psi (\sigma ) - \psi_{\infty, \tau} ( \sigma )  = (\widehat{\partial})^{\ast} \theta_\tau ( \sigma ) + (-1)^{n}\theta_\tau ( \partial \sigma ).
\end{align*}
Likewise, if $\sigma$ is $\kappa$-transverse we have 
\begin{align*}
&\widehat{m}_\kappa(\gamma  )^{\ast} \psi (\alpha ) - \psi_{\infty, \gamma , \kappa} ( \sigma )  = (\widehat{\partial})^{\ast} \theta_{\gamma , \kappa} ( \sigma ) + (-1)^{n-1}\theta_{\gamma, \kappa} ( \partial \sigma ) .
\end{align*}
\end{Proposition}

\begin{proof}
We show how the first identity is obtained. For the second we proceed identically. 
We write down for reference the two operators involved (see [\cite{KM}, Definition 22.1.3, Definition 25.3.3] ), namely the differential $\widehat{\partial} : \widehat{C}_{\ast}(Y , \s_{\xi_0} )\rightarrow \widehat{C}_{\ast -1} (Y , \s_{\xi_0, \alpha_0 , j_0} )$ and the chain map $\widehat{m}(U) : \widehat{C}_{\ast}(Y , \s_{\xi_0, \alpha_0 , j_0} )\rightarrow \widehat{C}_{\ast -2} (Y , \s_{\xi_0, \alpha_0 , j_0} )$, \begin{align}
&\widehat{\partial} = \begin{pmatrix} \partial^{o}_{o} & \partial^{u}_{o} \\ - \overline{\partial}^{s}_{u} \partial^{o}_{s} & - \overline{\partial}^{u}_{u} - \overline{\partial}^{s}_{u} \partial^{u}_{s} \end{pmatrix} \label{floerdiff2}\\
&\widehat{m}_\tau (U) = \begin{pmatrix} m^{o}_{o}(U) & m^{u}_{o}(U) \\
\overline{m}^{s}_{u}(U) \partial^{o}_{s} - \overline{\partial}^{s}_{u} m^{o}_{s}(U) & \overline{m}^{u}_{u} (U)+ \overline{m}^{s}_{u} (U)\partial^{u}_{s} - \overline{\partial}^{s}_{u} m^{u}_{s} (U) \end{pmatrix}.\label{Uchainmap}
\end{align}
Recall that we are interested in the duals $(\hat{\partial})^{\ast}$ , $\hat{m}_\tau (U)^{\ast}$ of the above operators, acting on cochains.


We first let $[\afr]$ be an irreducible critical point, and $z$ a component for which $\dim \mathcal{M}_{z}([\afr] , U , \sigma; \tau )= 1$. By Proposition \ref{strata} its compactification $\mathcal{M}_{z}([\afr] , U , \sigma ; \tau )$ has a codimension $1$ $\delta$-structure near the boundary stratum. This has the desirable property that the total count of boundary points (with orientations) vanishes [\cite{KM} , Corollary 21.3.2]. Then, enumerating the points on the boundary strata yields the identity

\begin{align*}
0 =  + & \# M_{z}([\afr] , Z_{\infty, \tau}(\sigma) )\\
+ & (-1)^{n} \cdot \sum_{\mathrm{subfaces} \,\, \Delta_{i}^{n-1} \subset \Delta^n} (-1)^{i}\# \mathcal{M}_{z}([\afr] , U , \sigma |_{\Delta_{i}^{n-1}}  ; \tau)\\
+ &  \sum_{[\bfr] \in \mathfrak{C}^{o}, z_1 , z_0} \# \breve{M}_{z_1}([\afr] , [\bfr] ) \# \mathcal{M}_{z_0} ([\bfr] , U , \sigma ; \tau ) \\
 - &  \sum_{[\bfr] \in \mathfrak{C}^{s},  [\cfr] \in \mathfrak{C}^{u} ,  z_2 , z_1 , z_0} \# \breve{M}_{z_2}([\afr] , [\bfr] ) \# \breve{M}_{z_1}([\bfr] , [\cfr] ) \# \mathcal{M}_{z_0} ([\bfr] , U , \sigma ; \tau ) \\
 - & \sum_{[\bfr] \in \mathfrak{C}^{o} , z_1 , z_0} \# M_{z_1}([\afr] , U , [\bfr]  ; \tau) \# M_{z_0}([\bfr] , \sigma )  \\
 + & \sum_{[\bfr] \in \mathfrak{C}^{s} , [\cfr] \in \mathfrak{C}^{u} , z_2, z_1 , z_0 } \# M_{z_2}([\afr] , U , [\bfr] ; \tau) \# \breve{M}_{z_1}([\bfr] , [\cfr] ) \# M_{z_0}([\cfr] , \sigma )  \\
 - & \sum_{[\bfr] \in \mathfrak{C}^{s} , [\cfr] \in \mathfrak{C}^{u} , z_2, z_1 , z_0 } \# \breve{M}_{z_2}([\afr], [\bfr] ) \# M_{z_1}([\bfr] , U , [\cfr] ; \tau ) \# M_{z_0}([\cfr] , \sigma ) 
   \end{align*}

Let $[\afr]$ be boundary-unstable now. The corresponding enumeration yields the identity
\begin{align}
0 = + & \# M_{z}([\afr] , Z_{\infty, \tau}(\sigma) ) .\nonumber \\
+ & (-1)^{n}\cdot\sum_{\mathrm{subfaces} \,\, \Delta^{n-1}_{i} \subset \Delta^n} (-1)^{i} \# \mathcal{M}_{z}([\afr] , U , \sigma |_{\Delta_{i}^{n-1}} ; \tau)\\
+  &  \sum_{[\bfr] \in \mathfrak{C}^{o}, z_1 , z_0} \# \breve{M}_{z_1}([\afr] , [\bfr] ) \# \mathcal{M}_{z_0} ([\bfr] , U , \sigma ; \tau )  \nonumber   \\
 - &  \sum_{[\bfr] \in \mathfrak{C}^{u} , z_1 , z_0}\# \breve{M}_{z_1}([\afr] , [\bfr] ) \# \mathcal{M}_{z_0} ([\bfr] , U , \sigma ; \tau ) \label{line1} \\
 - &  \sum_{[\bfr] \in \mathfrak{C}^{s},  [\cfr] \in \mathfrak{C}^{u} ,  z_2 , z_1 , z_0} \# \breve{M}_{z_2}([\afr] , [\bfr] ) \# \breve{M}_{z_1}([\bfr] , [\cfr] ) \# \mathcal{M}_{z_0} ([\bfr] , U , \sigma ; \tau )  \nonumber \\
 - & \sum_{[\bfr] \in \mathfrak{C}^{o} , z_1 , z_0} \# M_{z_1}([\afr] , U , [\bfr] ; \tau ) \# M_{z_0}([\bfr] , \sigma ) \nonumber \\
 - & \sum_{[\bfr] \in \mathfrak{C}^{u} , z_1 , z_0} \# M_{z_1}([\afr] , U , [\bfr] ; \tau) \# M_{z_0}([\bfr] , \sigma )\label{line2} \\
 + & \sum_{[\bfr] \in \mathfrak{C}^{s} , [\cfr] \in \mathfrak{C}^{u} , z_2, z_1 , z_0 } \# M_{z_2}([\afr] , U , [\bfr] ; \tau) \# \breve{M}_{z_1}([\bfr] , [\cfr] ) \# M_{z_0}([\cfr] , \sigma )  \nonumber \\
 - & \sum_{[\bfr] \in \mathfrak{C}^{s} , [\cfr] \in \mathfrak{C}^{u} , z_2, z_1 , z_0 } \# \breve{M}_{z_2}([\afr], [\bfr] ) \# M_{z_1}([\bfr] , U , [\cfr] ; \tau) \# M_{z_0}([\cfr] , \sigma ) \nonumber 
\end{align}


For the origin of the signs above we refer to Lemma \ref{signs2}\footnote{Again, we encounter the technical point that we must change some signs if we follow the \textit{reducible convention} for orienting the moduli $M_{z_1}([\afr] , [\bfr] )$ or $M_{z_1}([\afr] , U , [\bfr] )$ when both $[\afr], [\bfr]$ are boundary-unstable (see \S 20.6 \cite{KM}). This reducible convention is meant when writing the term $- \overline{\partial}^{u}_{u}$ in the Floer differential (\ref{floerdiff2}) and the term $\overline{m}^{u}_{u} (U)$ in (\ref{Uchainmap}). The signs listed in Lemma \ref{signs2} follow the usual convention. The only sign that one must add is $(-1)^{\mathrm{dim}_{z_1}M([\afr] , [\bfr] )} = -1 $ for line (\ref{line1}). In line (\ref{line2}) the sign is correct, since the difference between the two conventions is given by the sign $(-1)^{\mathrm{dim} M_{z_1}([\afr] , U , [\bfr] ; \tau )} = +1$. }.

For each of the two cases considered above, the corresponding identity can be written in terms of the natural pairing $\inner{ \cdot}{\cdot} : \widehat{C}^{\ast}(Y) \otimes_{R} \widehat{C}_{\ast}(Y) \rightarrow R$ as 
\begin{align*}
\inner{ \theta(U)(\sigma ) }{\widehat{\partial}[\afr]} + \inner{\theta (U) ( \partial \sigma )}{[\afr]}- \inner{ \psi(\sigma)}{\widehat{m}(U)[\afr]} + \inner{\psi_{\infty}(\sigma ) }{[\afr]} = 0 .
\end{align*}
This concludes the proof of the desired identity.
\end{proof}

\section{Exact triangles}\label{trianglessection}

For the whole of this section we assume that $\Lambda$ is a trivial double cover of $\mathcal{C}(Y , \xi_0 )$ and fix a trivialization. See Corollary \ref{triviallocalsystem} for a criterion that ensures this and which applies in particular if $\xi_0$ is strongly fillable. We work throughout with homology and cohomology with coefficients in a ring $R$. 

We recall that $\mathcal{C}(Y , \xi_0 , B ) \subset \mathcal{C}(Y , \xi_0 )$ denotes the subspace of contact structures $\xi$ which agree with $\xi_0$ over a Darboux ball $B$ (for $\xi_0 $) around the point $p\in Y$. The goal of this section is to establish Theorem \ref{thmtriangles1}. We rewrite this result in cohomological terms:

\begin{Theorem} \label{thmtriangles2}
Associated to any closed contact $3$-manifold $(Y , \xi_0 )$ for which the local system $\Lambda$ is trivial, there is a natural diagram which is commutative up to signs
\begin{center}
\begin{tikzpicture}[baseline= (a).base , trim left=-7cm]
\node[scale=.9] (a) at (0,0){
\begin{tikzcd}
\cdots \arrow{r}  & \widehat{HM}^{\ast} ( Y , \s_{\xi_0 }) \arrow{r}{U} &  \widehat{HM}^{\ast+ 2} ( Y , \s_{\xi_0} ) \arrow{r} & \widetilde{HM}^{\ast}(Y , \s_{\xi_0 } ) \arrow{r} & \widehat{HM}^{\ast+1}(Y , \s_{\xi_0 } ) \arrow{r}{U} & \cdots \\
\cdots \arrow{r} & H_{\ast}(\mathcal{C}(Y , \xi_0 )) \arrow{r} \arrow{u}{\fc} & H_{\ast-2}(\mathcal{C}(Y , \xi_0 , B )) \arrow{r} \arrow{u}{\fc} &  H_{\ast-1}(\mathcal{C}(Y , \xi_0 , B )) \arrow{r} \arrow{u}{\fctilde} &  H_{\ast-1}(\mathcal{C}(Y , \xi_0 )) \arrow{r} \arrow{u}{\fc}&  \cdots \\
 \end{tikzcd}
};
\end{tikzpicture}
\end{center}

where the top row is the long exact sequence of the mapping cone of $U$ in Floer cohomology, the bottom row is Wang's long exact sequence associated to the Serre fibration $\mathcal{C}(Y , \xi_0 , B ) \rightarrow \mathcal{C}(Y , \xi_0 ) \xrightarrow{ev} S^2$, the vertical arrows $\fc$ denote the families contact invariant, and $\fctilde$ is another families invariant which is to be defined.
\end{Theorem}

We recall that the "tilde" group $\HMtilde^\ast (Y, \mathfrak{s} )$ is defined as the cohomology of the algebraic mapping cone of the chain map $\widehat{m}_\tau (U)^\ast : \widehat{C}^\ast (Y, \mathfrak{s}) \rightarrow \widehat{C}^{\ast+2}(Y, \mathfrak{s} )$, for any choice of $\tau \in U(2)$. We denote $\widetilde{C}^\ast (Y, \mathfrak{s} ) := \mathrm{cone}(\widehat{m}_\tau (U)^\ast )$, and thus $\HMtilde (Y, \s ) = H^\ast (  \widetilde{C}^\ast (Y, \s ) )$. 

\begin{Remark}
The conventions we use for the algebraic mapping cone $\mathrm{cone}(f) $ of a chain map $f : A_{\ast} \rightarrow B_{\ast}$ (and similarly for a cochain map) are the following: as a module $\mathrm{cone}(f)_{\ast} = A_{\ast - 1} \oplus B_{\ast}$, and the differential in the cone is given by 
$$ d_{\mathrm{cone}} = \begin{pmatrix} - d_{A} & 0 \\ - f & d_{B} \end{pmatrix}.$$
Associated to a chain map $f: A_{\ast} \rightarrow B_{\ast}$ there is a sequence of chain maps $$ A_{\ast} \xrightarrow{f} B_{\ast} \xrightarrow{i} \mathrm{cone}(f) \xrightarrow{\delta} A_{\ast-1} $$ where $i (\beta ) = (0 , \beta )$ and $\delta (\alpha , \beta ) =  \alpha $. 
We recall that the sequence above becomes exact upon taking homology.
\end{Remark}

\subsection{Achieving transversality with constant families of perturbations} \label{betterchain}

To establish the above result, it is convenient to work with transverse singular chains in $\mathcal{C} = \mathcal{CM}(Y , \xi_0 ) \times \mathcal{P}$ a stronger property: that the perturbation term $\mathfrak{p} \in \mathcal{P}$ is constant for each simplex. The main result is the following



\begin{Proposition}\label{fixpert}
Let $\sigma : \Delta^n \rightarrow \mathcal{CM}(Y , \xi_0 )$ be any $C^2$ singular simplex. Then there exists a residual subset of perturbations $\mathfrak{p} \in \mathcal{P}$ and unitary splittings $\tau \in U(2)$ for which the singular chain $\sigma_{\mathfrak{p}} : \Delta^{n} \rightarrow \mathcal{C} = \mathcal{CM}(Y , \xi_0 ) \times \mathcal{P}$ defined by $\sigma_{\mathfrak{p}}(u ) = (\sigma(u ) , \mathfrak{p} )$ is $\tau$-transverse.
\end{Proposition}

We explain now how to establish this result. Let $\mathcal{M}$ stand for either of the moduli spaces $\mathfrak{M}(Z^+ ), \mathcal{M}(Z^+ )$ or $\mathcal{M}(U , Z^+ ; \tau )$. Recall that there is a natural Fredholm map $\mathcal{M} \xrightarrow{\pi} A \times \mathcal{P}$ where $A = \mathcal{CM}(Y , \xi_0 )$ in the first case, and $A = \mathcal{CM}(Y , \xi_0 ) \times \R$ in the other two. We write $\mathrm{pr}: A \times \mathcal{P} \rightarrow \mathcal{CM}(Y , \xi_0 ) \times \mathcal{P}$ for the natural projection in all cases above. In \S \ref{transversality} we deal with establishing the various transversality statements used in this paper. From the arguments there, we can deduce a finer transversality property than those stated thus far: that in order to achieve transversality for $\mathcal{M}$ one does not need to consider variations along the $A$ direction. Essentially, this is a consequence of the fact that the fibre product construction of the moduli space involved a restriction map to the slice $ 0 \times Y$, over which the family of spin-c structures was constant, independent of $A$. We have the following result, which will follow from \S \ref{transversality}:
\begin{Lemma}\label{transv2}
The map $\mathcal{M} \xrightarrow{\pi} A \times \mathcal{P} \xrightarrow{\mathrm{pr}_1} A$ is a submersion.
\end{Lemma}




\begin{proof}[Proof of Proposition \ref{fixpert}]

Fix a residual $\tau$ (such that Lemma \ref{surjectiveinfty1} holds). Then we need to show that $\sigma_{\mathfrak{p}}$ is transverse to $\pi^{\prime} := \mathrm{pr} \circ \pi: \mathcal{M} \xrightarrow{\pi} A \times \mathcal{P} \xrightarrow{\mathrm{pr}} \mathcal{CM}(Y , \xi_0 ) \times \mathcal{P}$ along components of $\mathcal{M}$ with $\mathrm{ind} \pi^{\prime} \leq 1-n$. 

By Lemma \ref{transv2}, the product map $\sigma \times \mathrm{id}_{\mathcal{P}} : \Delta^n \times \mathcal{P} \rightarrow \mathcal{CM}(Y , \xi_0 ) \times \mathcal{P}$ is transverse to $\pi^{\prime}$, and so their fibre product is transverse:\\

\begin{tikzcd}
\mathcal{M}(\sigma ):= \mathrm{Fib}(\sigma \times \mathrm{id}_{\mathcal{P}} , \pi^{\prime} ) \arrow{r} \arrow{d}{\pi_{\sigma}^{\prime}} &\mathcal{M} \arrow{d}{ \pi^{\prime}} \\
\Delta^n \times \mathcal{P} \arrow{r}{\sigma \times \mathrm{id}_{\mathcal{P}}} &  = \mathcal{CM}(Y , \xi_0 ) \times \mathcal{P} .\\
\end{tikzcd}

Now, the $C^{2}$ map $\mathrm{pr}_2 \circ \pi_{\sigma}^{\prime}: \mathcal{M}(\sigma) \rightarrow \mathcal{P}$ is Fredholm and has index $\mathrm{ind}(\pi_{\sigma}^{\prime}) = \mathrm{ind} \pi^{\prime} + n$, where $\mathrm{ind} \pi^{\prime}$ depends on the component of $\mathcal{M}$. The Sard-Smale theorem \cite{sardsmale} gives us a residual subset of perturbations $\mathfrak{p} \in \mathcal{P}$ which are regular values for the map $\mathrm{pr}_2 \circ \pi_{\sigma}^{\prime}$, provided that $\mathrm{ind} (\mathrm{pr}_2 \circ \pi_{\sigma}^{\prime} ) \leq 1$ (because $\mathrm{pr}_2 \circ \pi_{\sigma}^{\prime}$ is $C^2$). For those $\mathfrak{p}$, the map $\iota_{\mathfrak{p}} :\Delta^n \rightarrow \Delta^n \times \mathcal{P}$ given by $u \mapsto (u , \mathfrak{p} )$ is transverse to $\pi_{\sigma}^{\prime}$, and we obtain a transverse fibre product:\\

\begin{tikzcd}
M(\sigma ) : = \mathrm{Fib}(\iota_{\mathfrak{p}} , \pi_{\sigma}^{\prime} ) \arrow{r} \arrow{d} & \mathcal{M}(\sigma ) \arrow{d}{\pi_{\sigma}^{\prime}} \\
\Delta^n \arrow{r}{\iota_{\mathfrak{p}}} & \Delta^n \times \mathcal{P} \\
\end{tikzcd}

A simple diagram chasing argument involving the two diagrams above shows now that $\sigma_{\mathfrak{p}} = (\sigma \times \mathrm{id}_{\mathcal{P}} ) \circ \iota_{\mathfrak{p}}$ is transverse to $\pi^{\prime}$.
\end{proof}

\subsection{The map between triangles}


\begin{Remark}
The assignments $\sigma \mapsto \psi (\sigma )\, , \, \psi_{\infty, \tau} (\sigma ) \, , \, \theta_{\tau}(U) (\sigma )$ satisfied the chain map or chain homotopy up to signs (see Proposition \ref{chainmap}, Proposition \ref{mainidentities}). To avoid dealing with this plethora of signs, we find it convenient to resolve this issue now, by redefining $\psi_{\tau}, \psi_{\infty, \tau}, \theta_{\tau}(U)$ simply by placing the sign $(-1)^{\frac{n(n+1)}{2}}$ whenever they act on a singular simplex $\sigma$ of dimension $n$. It is straightforward to verify that these redefined operations now satisfy the chain map or chain homotopy relations strictly:
\begin{align*}
& \psi \partial = \widehat{\partial}^{\ast} \psi\\
& \psi_{\infty, \tau }\partial = \widehat{\partial}^{\ast} \psi_{\infty, \tau}\\
& \widehat{m}_{\tau}(U)^{\ast} \psi - \psi_{\infty, \tau} = \widehat{\partial}^{\ast} \theta_{\tau} (U) + \theta_{\tau}(U) \partial.
\end{align*}
\end{Remark}

Consider the subspace $\mathcal{CM}(Y , \xi_0 , B) \subset \mathcal{CM}(Y , \xi_0 )$ of triples $(\xi, \alpha , j )$ which over $B$ agree with the fixed triple $( \xi_0 , \alpha_0 , j_0 )$. Of course, the forgetful map gives a homotopy equivalence $\mathcal{CM}(Y , \xi_0 , B ) \xrightarrow{\simeq} \mathcal{C}(Y, \xi_0 , B )$. The next result uses the $\tau$-transverse singular simplices $\sigma_{\mathfrak{p}} : \Delta^n \rightarrow \mathcal{C} = \mathcal{CM}(Y, \xi_0 , B ) \times \mathcal{P}$ from Proposition \ref{fixpert}, and it shows that the $U$ action annihilates the image of $H_{\ast}(\mathcal{C}(Y , \xi_0 , B ))$ in $H_{\ast}(\mathcal{C}(Y , \xi_0 ) )$:

\begin{Lemma}\label{vanishingfixed}
Let $\sigma : \Delta^n \rightarrow \mathcal{CM}(Y, \xi_0 , B )$ be a singular simplex such that $\sigma_{\mathfrak{p}}$ is $\tau$-transverse for the perturbation $\mathfrak{p} \in \mathcal{P}$ (see Proposition \ref{fixpert}). Then $\psi_{\infty, \tau} (\sigma_\mathfrak{p} ) = 0$. Thus, by Proposition \ref{mainidentities} we have 
\[\widehat{m}_{\tau}(U)^{\ast} \psi( \sigma_{\mathfrak{p}}) = \widehat{\partial}^\ast \theta_{\tau} (U) (\sigma_{\mathfrak{p}}) + \theta_{\tau} ( U) (\partial \sigma_{\mathfrak{p}} ).\]
\end{Lemma}


\begin{proof}
Observe that the moduli spaces $M_{z}([\afr] , Z_{\infty, \tau} (\sigma_\mathfrak{p} ))$ of dimension $0$ are empty. The point is that $\sigma : \Delta^n \rightarrow \mathcal{CM}(Y , \xi_0 , B )$ parametrises triples that agree with $( \xi_0 , \alpha_0 , j_0 )$ on a neighbourhood of $p \in Y$. Thus the limiting set $Z_{\infty , \tau} (\sigma_\mathfrak{p} ) \subset \Delta^n$ must be either empty, or equal to $\Delta^n$. But $Z_{\infty, \tau} (\sigma_\mathfrak{p} ) \subset \Delta^n$ is a codimension $2$ submanifold with corners that is cut out transversely by the $\tau$-transversality assumption. Thus $Z_{\infty , \tau}(\sigma_\mathfrak{p})$ must be empty, and hence the moduli $M_{z}([\afr] , Z_{\infty , \tau}(\sigma_\mathfrak{p} ))$ of dimension $0$ are empty.
 \end{proof}

With this in place we now start building the map between the exact triangles from Theorem \ref{thmtriangles2}. First, we define a homomorphism
\[
\Psi_{\infty} : H_\ast \big( \mathcal{C}(Y , \xi_0 ) , \mathcal{C}(Y, \xi_0 , B ) \big) \rightarrow \widehat{HM}^{\ast+2}(Y, \mathfrak{s}_{\xi_0} )
\]
as follows. Choose a relative homology class in $H_\ast \big( \mathcal{C}(Y , \xi_0 ) , \mathcal{C}(Y, \xi_0 , B ) \big)$ together with a chain representing it, and lift this to a chain $T$ in $\mathcal{CM}(Y, \xi_0 )$ with boundary in $\mathcal{CM}(Y, \xi_0 , B)$. Choosing a generic $\tau$ and $\mathfrak{p}$, we can ensure that for each simplex $\sigma$ of $T$ the simplex $\sigma_\mathfrak{p}$ is $\tau$-transverse. Let $T_\mathfrak{p}$ stand for the cycle in $\mathcal{CM}(Y, \xi_0 , B ) \times \mathcal{P}$ obtained by using the perturbation $\mathfrak{p}$ on each simplex of $T$. We then have that $\psi_{\infty, \tau}(T_\mathfrak{p}) = \psi_{\infty, \tau}( \partial T_\mathfrak{p} ) = 0$, where the vanishing occurs due to Lemma \ref{vanishingfixed}. We can thus define $\Psi_{\infty} ([T]) := [ \psi_{\infty, \tau} (T_\mathfrak{p} ) ]$.

Next, we define a homomorphism
\[
\fctilde : H_{\ast-1} \big( \mathcal{C}(Y, \xi_0 , B ) \big) \rightarrow \widetilde{HM}^\ast (Y, \mathfrak{s}_{\xi_0} )
\]
as follows. We choose a class in $H_\ast \big( \mathcal{C}(Y, \xi_0 , B ) \big)[1]$ and a cycle in $\mathcal{C}(Y, \xi_0 , B )$ representing it. Then lift this to a cycle $T$ in $\mathcal{CM}(Y, \xi_0 , B )$. By choosing $\tau$ and $\mathfrak{p}$ generic we can again ensure that $T$ is made up of simplices $\sigma$ such that $\sigma_\mathfrak{p}$ is $\tau$-transverse. It follows from Lemma \ref{vanishingfixed} that the cochain $(\psi (T_\mathfrak{p} ) , - \theta_\tau (U)(T_\mathfrak{p} ) ) \in \widetilde{C}^\ast (Y, \s )$ is closed for the mapping cone differential. Thus, we can define
\[
\fctilde([T]) := [ ( \psi (T_\mathfrak{p}) ,  \theta_\tau (U) (T_\mathfrak{p} ) )].
\]

\begin{Lemma}\label{triangleLem}
The maps $\fc$, $\Psi_\infty$ and $\fctilde$ fit into a commutative diagram (up to signs)
\begin{center}
\begin{tikzpicture}[baseline= (a).base , trim left=-7cm]
\node[scale=.9] (a) at (0,0){
\begin{tikzcd}
\cdots \arrow{r}  & \widehat{HM}^{\ast} ( Y , \s_{\xi_0 }) \arrow{r}{U} &  \widehat{HM}^{\ast+ 2} ( Y , \s_{\xi_0} ) \arrow{r} & \widetilde{HM}^{\ast}(Y , \s_{\xi_0 } ) \arrow{r} & \widehat{HM}^{\ast+1}(Y , \s_{\xi_0 } ) \arrow{r}{U} & \cdots \\
\cdots \arrow{r} & H_{\ast}(\mathcal{C}(Y , \xi_0 )) \arrow{r} \arrow{u}{\fc} & H_{\ast}\big( \mathcal{C}(Y , \xi_0 ), \mathcal{C}(Y, \xi_0 , B )\big) \arrow{r} \arrow{u}{\Psi_\infty } &  H_{\ast-1}(\mathcal{C}(Y , \xi_0 , B )) \arrow{r} \arrow{u}{\fctilde} &  H_{\ast-1}(\mathcal{C}(Y , \xi_0 )) \arrow{r} \arrow{u}{\fc}&  \cdots \\
 \end{tikzcd}
 };
 \end{tikzpicture}
 \end{center}
where the top row is the Gysin long exact sequence associated to the mapping cone of $\widehat{m}_\tau (U)$, and the bottom row is the long exact sequence in homology of the pair $\big( \mathcal{C}(Y, \xi_0 ) , \mathcal{C}(Y, \xi_0 , B )\big) $.
\end{Lemma}

\begin{proof}
The third square is immediate. The first square commutes because of the identity $\widehat{m}_\tau (U) \psi - \psi_{\infty, \tau } = \widehat{\partial}^\ast \theta_\tau (U) + \theta_\tau (U) \partial $, and the second square anti-commutes because of the identity
\[
\begin{pmatrix}
0 \\
- \psi_{\infty, \tau} 
\end{pmatrix}
 = 
 \begin{pmatrix}
\psi \partial\\
\theta_\tau (U) \partial
 \end{pmatrix}
 + 
 \begin{pmatrix}
 - \widehat{\partial}^\ast & 0 \\
 - \widehat{m}_\tau (U) & \widehat{\partial}^\ast
 \end{pmatrix}
 \cdot 
 \begin{pmatrix}
\psi \\
\theta_\tau (U)
 \end{pmatrix}.
\]
\end{proof}

The final step is to identify the bottom row in Lemma \ref{triangleLem} as Wang's long exact sequence, and this is done by using a standard excision isomorphism 
\begin{align}
 H_{n-2}(\mathcal{C}(Y , \xi_0 , B ) )  \xrightarrow{\cong} H_n (\mathcal{C}(Y , \xi_0 ) , \mathcal{C}(Y , \xi_0 , B ) ).\label{excisioniso}
 \end{align}
 which follows from the fibration property of $ \mathcal{C}(Y , \xi_0 , B ) \rightarrow \mathcal{C}(Y , \xi_0 ) \xrightarrow{ev} S^2 $ (recall that $ev (\xi) = \xi (p) $). Let us recall how this isomorphism is constructed since we will need it. Let $x_0 \in S^2$ be the point corresponding to the plane $\xi_0 (p )$, and $-x_0 \in S^2$ its antipodal. We take the standard CW structure of $S^2$, where $x_0$ is the $0$-cell, and the $2$-cell $D^2$ is centered at $-x_0$. The map $1 \times \mathrm{pr} : ev^{-1}(-x_0 ) \times D^2 \rightarrow S^2 $ which collapses $\partial D^2$ to the point $x_0 \in S^2$ can be lifted through the fibration $ev : \mathcal{C}(Y , \xi_0 ) \rightarrow S^2$ to produce a map of pairs $$f : ( ev^{-1} (-x_0 ) \times D^2 , ev^{-1} (-x_0 ) \times \partial D^2 ) \rightarrow ( \mathcal{C}(Y , \xi_0 ) , ev^{-1} (x_0 ) ) $$ which at the center $-x_0 \in D^2$ agrees with the fibre inclusion $ev^{-1} (-x_0 ) \hookrightarrow \mathcal{C}(Y , \xi_0 )$. The map $f$ is a homotopy equivalence of pairs. The pair $( \mathcal{C}(Y , \xi_0 ) , ev^{-1} (x_0 ) )$ is weakly homotopy equivalent to the pair $( \mathcal{C}(Y ,\xi_0 ) , \mathcal{C}(Y , \xi_0 , B ) )$, so their homology is identified. The fibre transport along a path joining $x_0$ to $-x_0$ combined with the Künneth isomorphism yields an isomorphism
 $$ t_\ast : H_{n-2}(\mathcal{C}(Y , \xi_0 , B ) ) \xrightarrow{\cong} H_{n-2}(ev^{-1}(-x_0 ) )\xrightarrow{\cong} H_n( ev^{-1} (-x_0 ) \times D^2  , ev^{-1} (-x_0 ) \times \partial D^2 ).$$
The map $t_{\ast}$ is independent of the chosen path joining $x_0 , -x_0$ by $\pi_1 S^2 = 0$. Then the excision isomorphism (\ref{excisioniso}) is concretely described as the map $f_{\ast} \circ t_{\ast}$. Equipped with this description, we re-identify the map $\Psi_{\ast}$:

\begin{Lemma}\label{ident2lem}
Under the excision isomorphism (\ref{excisioniso}), the map $\Psi_{\infty}$ is identified with the restriction to $H_{\ast-2} \big( \mathcal{C}(Y, \xi_0 , B ) \big)$ of the families contact invariant $\fc : H_{\ast -2} \big( \mathcal{C}(Y, \xi_0 , B ) \big) \rightarrow \widehat{HM}^{\ast +2}(Y, \mathfrak{s}_{\xi_0} )$. 
\end{Lemma}
 
 \begin{proof}
 Choose an $(n-2)$-cycle $T$ in $\mathcal{C}(Y , \xi_0 , B )$. The homology class $[T]$ corresponds on the right-hand side of (\ref{excisioniso}) to the class of the chain $\widetilde{T} = f ( T^{\prime} \times D^2) $, where $T^{\prime}$ is the cycle in $ev^{-1}(-x_0 )$ obtained by transporting $T$ along a path from $x_0$ to $-x_0$. Thus, we need to compute the class $\Psi_{\infty} ([\widetilde{T}])$. 
 
The chain $\widetilde{T}$ can be lifted to a chain in $\mathcal{CM}(Y, \xi_0 )$ with boundary in $\mathcal{CM}(Y, \xi_0 , B )$, which we denote by $\widetilde{T}$ as well. As usual, we can assume that all simplices that make up $\widetilde{T}$ are $\tau$-transverse when equipped with a fixed perturbation $\mathfrak{p}$, for a generic choice of $\tau$ and $\mathfrak{p}$. We have $\Psi_{\infty}([\widetilde{T}]) = [ \psi_{\infty, \tau} ( \widetilde{T} )]$, and by construction the chain $\psi_{\infty, \tau} ( \widetilde{T} )$ agrees with the chain $\psi ( T_{\infty} )$ where $T_{\infty}$ is obtained by intersecting $\widetilde{T}$ with the union over the limiting loci $Z_{\infty , \tau}(\sigma)$ with $\sigma$ running over the subfaces of $\tilde{T}$ of all dimensions. By transversality of these intersections, $T_{\infty}$ can be given the structure of a chain. Now, as in the proof of Proposition \ref{geominterpcont}, we know that $T_{\infty}$ agrees with the intersection of $\widetilde{T}$ and a fibre of $ev : \mathcal{C}(Y , \xi_0 ) \rightarrow S^2$. By the description of $\widetilde{T}$ as $f (T^\prime \times D^2 )$ we can now see that $T_{\infty}$ is, in fact, a cycle in $\mathcal{C}(Y, \xi_0 )$ which is homologous to $T$. Thus, we have $\Psi_\infty ([\widetilde{T}]) = \psi_{\ast}( [T] )$, as required. 
 \end{proof}

Combining Lemma \ref{triangleLem} and Lemma \ref{ident2lem}, the proof of Theorem \ref{thmtriangles2} is now complete.

\section{Transversality, compactness and orientations}\label{appendix}

\subsection{Transversality}\label{transversality}

We now take up the task of establishing the transversality results claimed in the previous sections. The arguments used follow quite closely those of \cite{KM} and \cite{marianothesis}, and we will focus on describing the differences. This section has the nature of an appendix.

We recall that we have chosen integers $k \geq 4$ and $l$ with $l-k-2 \geq 1$.

\subsubsection{Main results}

We consider the following setup, in the spirit of the one considering thus far. We consider a $P$-family of Riemannian metrics $\{g_p\}$ on $Z^+ = \R \times Y$. As before, we consider metrics of regularity $C^{l-1}$. The parametrising space $P$ is a Banach manifold, possibly just finite-dimensional. The cases we have in mind are mainly $P = \mathcal{CM}(Y, \xi_0 )$ and $P = \mathcal{CM}(Y, \xi_0 ) \times \R$. We assume that the metrics $g_p$ coincide with a fixed cylindrical metric $g_0 = dt^2 + g_{0 , Y}$ over the region $(- \infty , 1/2] \times Y$. We assume that $K = [1, + \infty) \times Y$ is equipped with a family of almost-K\"{a}hler structures $\{(\omega_p , J_p , g_p )\}$ such that $g_p = \omega_p ( \cdot , J_p \cdot )$. We assume that each $( \omega_p , J_p , g_p)$ makes $K$ an asymptotically flat almost-K\"{a}hler end (for the definition see \cite{monocont} , \S 3(i)). We also assume that the differences $g_{p} g_{0}^{-1}$ are bounded over $Z^+$ (though not necessarily uniformly in $P$). There is then a $P$-family of spin-c structures on $Z^+$ constructed as in \S \ref{confscon} using the triple $(\omega_0 , J_0 , g_0 )$ and the $P$-family of metrics.

We remark at this point that if we consider a compact end or a cylindrical end, rather than an asymptotically flat almost-K\"{a}hler end, then the results of this section will still apply. 

The corresponding space of configurations $(A, \Phi , p)$ over $K^{\prime} = [0, + \infty) \times Y$, and its quotient by the group $\mathcal{G}_{k+1}(K^{\prime})$ of $L^{2}_{k+1}$ gauge transformations asymptotic to the identity, are denoted by $\mathcal{C}_k (K^{\prime} )$ and $\mathcal{B}_{k}(K^{\prime} )$ respectively. $\mathcal{B}_{k}(K^{\prime})$ is a $C^{l-k-2}$ Banach manifold away from the reducible locus. The moduli space $\mathfrak{M} (K^{\prime} ) \subset \mathcal{B}_{k}(K^{\prime}) \times \mathcal{P}$ is defined as the zero set of the perturbed Seiberg--Witten map $\mathrm{sw}_{\eta}$, which is naturallly a section of a Hilbert bundle $\Upupsilon_{k-1}$ over $\mathcal{B}_k (K^{\prime} ) \times \mathcal{P}$. The perturbation $\eta$ is taken of the form
$$\eta ( A, \Phi , p,  \mathfrak{p}  ) = \varphi_{p}^1 \hat{\mathfrak{q}}(A, \Phi ) + \varphi_{p}^2 \hat{\mathfrak{p}} (A, \Phi )   + \varphi_{p}^{3} \hat{\mathfrak{p}}_{K, p} .$$ 

Here $\mathfrak{q}$ is a fixed admissible perturbation, and $\mathfrak{p}_{K,p}$ is the Taubes perturbation used earlier. We consider $P$-families of functions satisfying similar constraints as before. Namely, (i) cutoff functions $\varphi_{p}^1$ which are identically $1$ on a neighbourhood of $(- \infty ,  0]$ and vanishing on a neighbourhood of $[1/2, + \infty )$; (ii) bump functions $\varphi_{p}^2$ with compact support inside $(0, 1/2)$; (iii) and cutoff functions $\varphi_{p}^{3}$ which are identically $1$ over $[1, + \infty)$ and vanish on a neighbourhood of $(- \infty , 1/2]$. One can include more perturbation terms in $\eta$ to adjust to each particular situation, and as long as they don't depend on $\mathcal{P}$ the results of this section apply.

Below we introduce a suitable class of maps $\mathrm{ev} : \mathcal{C}_k (K^{\prime} ) \rightarrow V$ that we call \textit{good} (Definition \ref{permiss} below). These are equivariant sections of a $\mathcal{G}_{k+1}(K^{\prime})$-equivariant fibre bundle $V$ over $\mathcal{C}_k (K^{\prime} )$, and we wish to impose the constraint that $\mathrm{ev}(A , \Phi , p , \mathfrak{p} ) = \sigma (A , \Phi , p , \mathfrak{p} )$ on the moduli $\mathfrak{M}(K^{\prime} )$. Here $\sigma $ is a fixed equivariant section of $V$ (Definition \ref{permiss}). We prove:

\begin{Proposition}\label{Transv1}
The Seiberg--Witten map $\mathrm{sw}_{\eta} :  \mathcal{B}_{k}(K^{\prime} ) \times \mathcal{P} \rightarrow \Upupsilon_{k-1}$  is transverse to the zero section. If $\mathrm{ev}$ is a good evaluation map, then $\mathrm{ev}$ and $\sigma$ are transverse sections of $V \rightarrow \mathfrak{M}(K^{\prime} )$.
\end{Proposition}

Thus, the space $\mathfrak{M}_k (K^{\prime} , \mathrm{ev} ) := \mathfrak{M}(K^{\prime}) \cap \mathrm{ev}^{-1}(\mathrm{Im} \sigma )$ is a Banach manifold, of class $C^{l-k-2}$. We have two restriction maps to the configuration space of the boundary \begin{align}
& R_{+} : M^{\ast}([\afr] , (- \infty , 0] \times Y  ) \rightarrow \mathcal{B}_{k-1/2}^{\ast}(Y) \label{R+}\\
& \mathfrak{R}_{-} : \mathfrak{M} (K^{\prime} , \mathrm{ev}  ) \rightarrow \mathcal{B}_{k-1/2}^{\ast}(Y).\label{R-}
\end{align}
and their fibre product $\mathfrak{M} ( [\afr] , Z^+, \mathrm{ev} ) = \mathrm{Fib} ( R_{+} , \mathfrak{R}_{-} )$ is the moduli space we are interested in. The main result of this section is:

\begin{Proposition}\label{Transv2}
For a good evaluation map, the maps $R_{+}$ and $\mathfrak{R}_{-}$ are transverse. Thus, the moduli space $\mathfrak{M}([\afr] , Z^+, \mathrm{ev} )_P$ is a $C^{l-k-2}$ Banach manifold. The map $\mathfrak{M}([\afr] , Z^+, \mathrm{ev} ) \rightarrow P \times \mathcal{P}$ is $C^{l-k-2}$ and Fredholm.
\end{Proposition}

We now describe the class of evaluation maps for which our transversality results apply.


\begin{Definition}\label{permiss}
Fix a smooth $\mathcal{G}_{k+1}(K^{\prime})$-equivariant fibre bundle $V \rightarrow \mathcal{C}_k (K^{\prime})$ with finite-dimensional fibre, together with a preferred equivariant section $\sigma$ and a connection on $V$ along $\sigma$ (that is, a connection on the pullback fibre bundle $\sigma^{\ast} V$). A \textit{good} evaluation map compatible with such data is a section $\mathrm{ev} : \mathcal{C}_{k}(K^{\prime}) \rightarrow V$ subject to the following conditions:
\begin{enumerate}[(i)]
\item $\mathrm{ev}$ is a $\mathcal{G}_{k+1}(K^{\prime})$-equivariant section
\item $\mathrm{ev}$ is transverse to $\sigma$ as sections of $V \rightarrow \mathcal{C}_k (K^\prime )$
\item There exists a compact set $E \subset (1/2 , 1) \times Y$ with $([0 , + \infty) \times Y) \setminus E$ \textit{connected} such that, for any $( A, \Phi , p) \in \mathrm{ev}^{-1}(\mathrm{Im} \sigma)$, all the smooth configurations tangent to $( A , \Phi , p)$ of the form $$( a, \phi , 0) \in T_{( A , \Phi , p)} \mathcal{C}_{k}([0 , + \infty) \times Y)$$ which are compactly supported away from $E$ are contained in $$T_{(A , \Phi , p )} \big( \mathrm{ev}^{-1}(\mathrm{Im} v) \big) = \ker (\mathcal{D} \mathrm{ev} - \mathcal{D} \sigma)_{(A, \Phi , p)}.$$
\end{enumerate}

In (iii), $\mathcal{D}$ denotes the differential of a section projected onto the the vertical direction using the connection $V$ defined along $\sigma$.
\end{Definition}

The evaluation constraints we have considered thus far in the article fall into the above category. These are: 
\begin{Example}\label{exev1}
One of the main examples (see \S \ref{universalUmoduli}) is the evaluation map $\mathrm{ev}: (A, \Phi , p) \mapsto \tau_1 \overline{\Phi}(x_0)$ induced by an unitary splitting $\tau$ of the fibre of the spinor bundle $S^+$ at a point $x_0 \in (1/2 , 1 ) \times Y$. Here $V$ is the trivial bundle with fibre $\C$ carrying the $\mathcal{G}_{k+1}(K^{\prime})$-action $v \cdot \lambda = v(x_0) \lambda$, and $\sigma$ is the zero section. The subset $E$ can be taken to be the point $x_0$. 

In \S \ref{universalUmoduli} we considered the moduli $\mathcal{M}([\afr] , Z^+ , U ; \tau )$ with an evaluation constraint that travelled along the $\R$ direction: $\tau_1 \overline{\Phi}(s,y_0) = 0$, $s \in \R$. By applying an $\R$-family of diffeomorphisms taking the point $(s , y_0 )$ to $(1/2 , 1) \times \{y_0 \}$ we see that the situation considered in \S \ref{universalUmoduli} fall into our current setup.
\end{Example}

\begin{Example}\label{exev2}
Another example (see \S \ref{universalUmoduli}) is the half-holonomy evaluation map, associated to a smooth oriented closed curve $\gamma \subset  Y$, given by $h_{\gamma}(A , \Phi ,p ) = \mathrm{exp} \frac{1}{2} \int_{s_{0} \times \gamma} \hat{A}$, where $s_0 $ is a fixed number in $(1/2 ,1 )$. Here, the fibre bundle $V$ is equivariantly trivial, with fibre $\mathrm{U}(1)$, and we take $\sigma$ constant. $E$ can be taken to be $s_0 \times \gamma \subset (1/2 , 1) \times Y$.
\end{Example}

\begin{Example}\label{exev3}
In \S \ref{compactifications} we considered a map that evaluates the canonical spinors at a point $x_0$. The zero set of this map are the limiting moduli space $M([\afr] , Z_{\infty , \tau}(\sigma ) )$ (Definition \ref{Ulimiting}). This was defined, after choosing an unitary splitting, by $\mathrm{ev} : ( A , \Phi ,p ) \mapsto \tau_1 \overline{\Phi}_{ p} (x_0)$. The bundle $V$ of which $\mathrm{ev}$ is a section is a vector bundle with trivial $\mathcal{G}_{k+1}(K^{\prime})$ action, and $\sigma$ is the zero section. This defines a good evaluation map for generic unitary splittings (Lemma \ref{surjectiveinfty1}). An analogous map was considered for the half-holonomy of the canonical connections. 
\end{Example}

Recall from \S \ref{confspaces} that $\mathcal{C}(K^{\prime} ) \rightarrow P$ is a bundle of affine Hilbert spaces equipped with a preferred connection on $\mathcal{C}(K^{\prime} ) \rightarrow P$ i.e. a complementary (horizontal) subbundle for the vertical subbundle of $T \mathcal{C}_{k}(K^{\prime} )$.

\begin{Definition}
An admissible evaluation map $\mathrm{ev}: \mathcal{C}_{k}(K^\prime ) \rightarrow V$ is \textit{very good} for the data $V , \sigma$ if the transversality condition (ii) from Definition \ref{permiss} can be achieved without variation along the horizontal direction of $T \mathcal{C}(K^\prime)$.
\end{Definition}
\begin{Remark}
Examples \ref{exev1} and \ref{exev2} are very good, while \ref{exev3} is not.
\end{Remark}
Finally, we will show:

\begin{Proposition}\label{Transv3}
For a very good evaluation map $\mathrm{ev}$, the map $\mathfrak{M}([\afr] , Z^+, \mathrm{ev} ) \rightarrow P$ is a submersion.
\end{Proposition}

\subsubsection{Proof of Proposition \ref{Transv1}}

 
Let $\gamma = (A , \Phi , p , \mathfrak{p})$ be a configuration in $\mathcal{C}_{k}([0, + \infty) \times Y)$ solving the equations $\mathrm{sw}_{\eta}(\gamma ) = 0$, and denote by $\mathbf{d}_{\gamma} : L^{2}_{k+1}([0 , + \infty) \times Y, i \R ) \rightarrow T_{p} P \times L_{k}^{2}([0, + \infty) \times Y , i \Lambda^1 \oplus S^+ )$ the linearisation of the gauge action at $\gamma$. 
To establish the first statement in Proposition \ref{Transv1} it suffices to show the stronger result that the operator $Q_{\gamma} = (\mathcal{D} \mathrm{sw}_\eta )_{\gamma} + \mathbf{d}^{\ast}_{\gamma}$ is surjective. This operator takes the form
\begin{align*}
   L_{k}^{2}(K^{\prime} , i \Lambda^1 \oplus S^{+}) \times T_p P  & \times \mathcal{P}  \rightarrow  \\
   &  L_{k-1}^{2}(K^{\prime} , i \mathfrak{su}(S^+) \oplus S^- \oplus i \R )\end{align*}

The desired surjectivity is established in [\cite{marianothesis}, p.51] using similar ideas to \cite{KM}. We explain how to adapt these ideas to the case of the moduli space with evaluation constraint $\mathfrak{M}(K^{\prime} , \mathrm{ev} )$. For this suppose that $\gamma = ( A, \Phi , p , \mathfrak{p})$ satisfies, in addition, the constraint $\mathrm{ev} = \mathrm{\sigma}$. We have the vertical derivative at $\gamma $ of $\mathrm{ev} : \mathcal{C}_{k}([0 , + \infty ) \times Y)_P \rightarrow V$, which is a linear map of the form $$\mathcal{D} \mathrm{ev}_{\gamma } :   L_{k}^{2}(K^{\prime} , i \Lambda^1 \oplus S^{+}) \times T_p P \times \mathcal{P} \rightarrow V_0 $$
where $V_0$ denotes the fibre of $V$ at $\gamma$. We wish to establish the surjectivity of the operator $Q_{\gamma} + (\mathcal{D} \mathrm{ev}-\mathcal{D}\sigma )_{\gamma}$, which takes the form
\begin{align*}
 L_{k}^{2}(K^{\prime} , i \Lambda^1 \oplus S^{+}) & \times T_p P \times \mathcal{P}    \rightarrow \\
  &  L_{k-1}^{2}(K^{\prime} , i \mathfrak{su}(S^+) \oplus S^- \oplus i \R )\oplus V_0. 
\end{align*}

Equivalently, because $(\mathcal{D} \mathrm{ev}-\mathcal{D}\sigma )_{\gamma}$ is surjective (condition (ii) of Definition \ref{permiss}), it suffices to prove the surjectivity of 
\begin{align}
 \mathrm{ker}(\mathcal{D} \mathrm{ev}-\mathcal{D}\sigma )_{\gamma} \xrightarrow{Q_{\gamma }} L_{k-1}^{2}(K^{\prime} , i \mathfrak{su}(S^+) \oplus S^- \oplus i \R ). \label{defop}
 \end{align}

\begin{Lemma} \label{lembound}
For all $\gamma$ with $\mathrm{sw}_{\eta}(\gamma) = 0$ and $\mathrm{ev}(\gamma) = \sigma( \gamma )$, the operator (\ref{defop}) is surjective. 
\end{Lemma}

\begin{proof}
We follow the argument in the proof of [\cite{KM}, Proposition 24.3.1], and explain the necessary modifications. We first show that the image of $\mathrm{ker} (\mathcal{D} \mathrm{ev}-\mathcal{D}\sigma )_{\gamma} $ under $Q_{\gamma}$ is dense in the $L^2$ topology on the target. Suppose not, for a contradiction, and hence choose a non-zero element $V  \in L^2$ which annihilates the image of $\mathrm{ker} (\mathcal{D} \mathrm{ev}-\mathcal{D}\sigma )_{\gamma}$. 

In particular, it annihilates the image under $\mathcal{D} : = Q_{\gamma }(- , - , 0 , 0 )$ of the subspace consisting of all smooth configurations $(a , \phi , 0 , 0 )$ compactly supported away from $E$ (see (iii) of Definition \ref{permiss}). Now, $\mathcal{D}$ is an elliptic differential operator, so by elliptic regularity we obtain that $V$ is in $ L_{1, \mathrm{loc}}^{2}$ on $([0 , + \infty) \times Y ) \setminus E$. In particular $V$ is in $L^{2}_{1}$ on the collar neighbourhood $[0 , 1/2) \times Y$, where it satisfies the formal adjoint equation $\mathcal{D}^{\ast} V = 0$. By the unique continuation principle (see \cite{KM}: Lemma 7.1.3 for the cylindrical case, and the argument in Lemma 7.1.4 for arbitrary manifolds), because $V$ is not identically zero over $K^{\prime} = [0 , + \infty) \times Y$ and $E$ does not disconnect this set (see (iii) of Definition \ref{permiss}), we know that $V$ does not vanish identically on the collar $[0 , 1/2) \times Y$. The fact that $\mathcal{D}^{\ast}$ satisfies the unique continuation property follows from [\cite{KM}, eq. (24.15)]. We thus obtain that the restriction of $V$ to the boundary $0 \times Y$ is non-zero, again by the unique continuation principle [\cite{KM} , Lemma 7.1.3]. 

However, using the argument in the proof of [\cite{KM} ,Corollary 17.1.5] we can show that the restriction must be zero, by orthogonality of $V$ to the image. This gives the desired contradiction.

Thus, the image under $Q_{\gamma }$ of $\mathrm{ker} \mathcal{D} ( \mathrm{ev}- \sigma)_\gamma $ is dense in the $L^2$ topology. As we mentioned in the previous section, $\mathcal{D} : L^{2}_{k} \rightarrow L^{2}_{k-1}$ is surjective [\cite{marianothesis}, p.51], for which the argument is similar but simpler than this. Hence, the image of $\ker (\mathcal{D} \mathrm{ev}-\mathcal{D}\sigma )_{\gamma}$ under $Q_{\gamma }$ is of finite-codimension and dense in $L^2$, so (\ref{defop}) is surjective.

\end{proof}

The previous lemma, together with the fact that $Q_{\gamma}$ is also surjective on the bigger domain, and the surjectivity of $(\mathcal{D} \mathrm{ev}-\mathcal{D}\sigma )_{\gamma}$, complete the proof of Proposition \ref{Transv1} with the aid of the following observation from linear algebra:
\begin{Lemma}
Suppose $X \xrightarrow{q} Y$ and $X \xrightarrow{e} V$ are linear maps of vector spaces. Assume that the following maps are surjective: $q$ , $e$ and the restriction $\ker e \xrightarrow{q} Y $. Then $\ker q \xrightarrow{e} V$ is also surjective.
\end{Lemma}

\begin{proof}
The cokernel of $\mathrm{ker} e \xrightarrow{q} Y$ is 
\begin{align*}
\frac{V}{e(\ker q)}&  = \frac{e(X)}{e(\ker q)} \cong \frac{ X}{e^{-1} \big( e( \ker q ) \big)} = \frac{X}{\ker q + \ker e }\\
 &\cong \frac{ X / \ker q}{\ker e / \ker q \cap \ker e} \cong q(X) / q(\ker e ) = 0 .
\end{align*}
\end{proof}
\begin{Remark}
The proof of Lemma \ref{lembound} shows that the surjectivity of (\ref{defop}) is already achieved by the tangent configurations $\{ (a , \phi , 0 , 0 ) \} \subset \ker (\mathcal{D} \mathrm{ev}-\mathcal{D}\sigma )_{\gamma}$. In particular, the map $\mathfrak{M}(K^{\prime} , \mathrm{ev} ) \rightarrow \mathcal{P}$ is a submersion. If, in addition, $(\mathcal{D} \mathrm{ev}-\mathcal{D}\sigma )_{\gamma}$ achieves surjectivity without varying in $T_p P$ (the "very good" condition), then the map $\mathfrak{M}(K^{\prime} , \mathrm{ev} ) \rightarrow P \times \mathcal{P}$ will also be a submersion.
\end{Remark}

\subsubsection{Proof of Propositions \ref{Transv2} and \ref{Transv3}}


For Proposition \ref{Transv2} we need to establish the transversality of the fibre product. In other words, we need to check that the sum of the derivatives $$(d R_{+})_a + (d \mathfrak{R}_{-})_b : T_a M^{\ast}([\afr] , (- \infty , 0 ] \times Y ) \oplus T_b \mathfrak{M}(K^{\prime}   , \mathrm{ev}) \rightarrow T_{[\mathfrak{c}]} \mathcal{B}_{k - 1/2}(Y)$$ is surjective for each $(a,b)$ in the fibre product $\mathrm{Fib}(R_{+} , \mathfrak{R}_{-} )$, and $[\cfr]$ the restriction to the boundary. 

The sum $( (d R_{+})_a + (d \mathfrak{R}_{-})_b) (- ,-, 0, 0)$, i.e. acting on tangent directions which vanish on the $P$ and $\mathcal{P}$ directions, is a Fredholm operator. This can be extracted from \cite{marianothesis}, Lemma 26 (see assertions (3),(4),(7),(8)). Thus, $(dR_{+})_a + (d \mathfrak{R}_{-})_b$ has finite dimensional cokernel. This, together with Lemma \ref{density}, coming up next, shows that $(dR_{+})_a + (d \mathfrak{R}_{-})_b$ is surjective.

\begin{Lemma} \label{density}
Let $\gamma = ( A, \Phi , p , \mathfrak{p}) \in \mathfrak{M}([0 , + \infty) \times Y , \mathrm{ev} )_P$, and let $[\cfr] = \mathfrak{R}_{-} ( \gamma  )$. Then $$(d \mathfrak{R}_{-} )_{\gamma} : T_{ \gamma} \mathfrak{M}(K^{\prime} , \mathrm{ev} ) \rightarrow T_{[\cfr]} \mathcal{B}_{k-1/2}(Y)$$ has dense image in the $L_{1/2}^2$ topology. 
\end{Lemma}

\begin{proof}
We follow the proof of [\cite{KM} ,Lemma 24.4.8]. The result will follow if we show that the following operator has dense image in the $L^2 \times L_{1/2}^2$ topology. It is the operator given by the restriction of $( \mathcal{D} \mathrm{sw}_{\eta})_{\gamma}$ to $\ker (\mathcal{D} \mathrm{ev}-\mathcal{D}\sigma )_{\gamma} \oplus \mathcal{P}$, coupled with the derivative of the restriction $\tilde{\mathfrak{R}}_{-}$ to the configuration space of the boundary $\mathcal{C}_{k-1/2}(Y)$, which has the form 
\begin{align}
 \mathrm{ker} (\mathcal{D} \mathrm{ev}-\mathcal{D}\sigma )_{\gamma} \oplus \mathcal{P} \rightarrow    & L_{k-1}^{2}(K^{\prime} , i \mathfrak{su}(S^+) \oplus S^-  ) \label{operator6}\\
 & \oplus L_{k -1/2}^2(Y ; i \Lambda^1 \oplus S_{Y}) .\nonumber
\end{align}
Here $S_{Y}$ is the restriction of $S^+$ to $0 \times Y$, and note that $$\mathrm{ker} (\mathcal{D} \mathrm{ev}-\mathcal{D}\sigma )_{\gamma} \subset  L_{k}^2(K^{\prime} ; i \Lambda^1 \oplus S^+ ) \oplus T_p P .$$

We suppose for a contradiction that the image of this operator is not dense in $L^2 \times L_{1/2}^2$, and pick a non-zero $(V , v) \in L^2 \times L_{-1/2}^2$ which annihilates the image. By considering directions in $\mathcal{C}_{k}(K^{\prime})$ which are tangent to the gauge-orbit of $\gamma$ (these are contained in $\mathrm{ker} (\mathcal{D} \mathrm{ev}-\mathcal{D}\sigma )_{\gamma}$ by (i) of Definition \ref{permiss}) we see that $v$ is orthogonal to the directions tangent to the gauge-orbit through $\gamma_{|Y}$. \\

Consider the restriction map $r$ to the $dt$ component of the connection form at the boundary. We couple the previous operator with $r$ and the operator $\mathbf{d}^{\ast}_{\gamma}$ to obtain an operator on $\mathrm{ker} (\mathcal{D} \mathrm{ev}-\mathcal{D}\sigma )_{\gamma}$ by restriction of
\begin{align*}
 L_{k}^2 (K^{\prime};  & i \Lambda^1 \oplus  S^+ ) \oplus T_p P  \oplus \mathcal{P}  \xrightarrow{\mathfrak{Q}\oplus r}\\
& L_{k-1}^2 (K^{\prime} ; i \mathfrak{su}(S^+) \oplus S^-) \oplus L^{2}_{k-1}(K^{\prime} ;  i \R ) \\
\oplus &  L_{k-1/2}^{2}(Y ; i \Lambda^1 \oplus S_{Y} ) \oplus L_{k-1/2}^{2}(Y ; i \R ) .\end{align*}

The image of $\mathrm{ker} (\mathcal{D} \mathrm{ev}-\mathcal{D}\sigma )_{\gamma}$ under this operator is orthogonal to $(V, 0 , v , 0)$. The operator $Q := \mathfrak{Q}(   - , 0 , 0 )$ is elliptic. As in the proof of Lemma \ref{lembound}, $(V , 0 , v, 0)$ is orthogonal to the image of the smooth configurations $(a , \phi , 0 , 0 )$ that are compactly supported away $E$, which are contained inside $\ker (\mathcal{D} \mathrm{ev}-\mathcal{D}\sigma )_{\gamma}$ by (iii) of Definition \ref{permiss}. Elliptic regularity then implies that $V$ is in $L_{1, \mathrm{loc}}^{2}$ away from $E$; so $V$ is in $L^{2}_{1}$ on the collar $[0,1/2) \times Y$ since $E \subset (1/2 , 1)$.  Thus, $V$ satisfies the formal adjoint equation $Q^{\ast} V = 0$ over the collar $[0 , 1/2) \times Y$, and so $V$ does not vanish identically over this region by the unique continuation principle (similar argument as in the proof of Lemma \ref{lembound}). 

From this point on, the argument of [\cite{KM} ,Lemma 24.4.8] carries through without modification. Namely, by integrating by parts we see that $V_{|Y} = -v$ (under standard identifications of the corresponding bundles), and combining this with the fact that $v$ was orthogonal to the gauge orbit, an argument as in [\cite{KM} ,Lemma 15.1.4] shows that $V$ is orthogonal to the gauge-orbit on every slice $t \times Y$. Finally, the argument of [\cite{KM} , Proposition 15.1.3] produces, because $V$ does not vanish identically on the collar, a perturbation $\mathfrak{t} \in T_{\mathfrak{p}}\mathcal{P} = \mathcal{P}$ for which the derivative of (\ref{operator6}) in the direction of $(0 , 0,0, \mathfrak{t} )$ is not orthogonal to $(V,v)$, a contradiction. 
\end{proof}

If $\mathrm{ev}$ is very good, then no variation in the $P$ direction (horizontal) will be needed to achieve transversality in the previous Lemma. This, together with the Remark at the end of the previous subsection, gives us the stronger result of Proposition \ref{Transv3}.



\subsection{Compactness } \label{compactness}

Here we briefly describe some of the compactness results that lead to the construction of the compactified moduli spaces by broken configurations. Large part of the material presented here is a straightforward adaptation of results found in \cite{monocont}, \cite{boyu} and \cite{KM}.


The main moduli spaces that will concern us in this section are $\mathfrak{M}(Z^+ )$ and $ \mathcal{M}(Z^+)$, the latter because it contains the parametrised evaluation moduli $\mathcal{M}(U , Z^+ ; \tau)$ and $\mathcal{M}(\gamma , Z^+ ; \kappa )$. For simplicity in notation we state all results for $\mathcal{M}(Z^+ )$ or $\mathcal{M}(U , Z^+ ; \tau )$, but we note that the corresponding results for $\mathfrak{M}(Z^+)$ are analogous and simpler. As before, we use the notation $\mathcal{M}(\sigma )$ for the fibre product of $\mathcal{M}(Z^+) \rightarrow \mathcal{CM}(Y , \xi_0 ) \times \mathcal{P}$ and a $C^2$ singular simplex $\sigma : \Delta^n \rightarrow \mathcal{CM}(Y , \xi_0 ) \times \mathcal{P}$. $M(\sigma)$ is a $C^2$ manifold with corners provided transversality holds, and its points consist of gauge-equivalence classes $[(A, \Phi , t , s)]$ of Seiberg--Witten monopoles with $t \in \Delta^n$ , $s \in \R$. The projection to the $s \in \R$ coordinate is denoted $\pi_{\R} : \mathcal{M}(\sigma) \rightarrow \R$. The simplex $\sigma : \Delta^n \rightarrow \mathcal{CM}(Y , \xi_0 ) \times \mathcal{P}$ will be kept fixed throughout this section.



It is also convenient to introduce the moduli space $\mathcal{M}_{\mathrm{loc}}(Z^+ )$ of gauge-equivalence classes of solutions to the same equations as $\mathcal{M}(Z^+)$, also approaching the canonical configurations in $L^2_{k}$ on the conical end, but with \textit{no} asymptotics to critical points on the cylindrical end. The relevant gauge group involved in the quotient is now the topological group $\mathcal{G}_{\mathrm{loc}}$ of locally $L^{2}_{k+1}$ gauge transformations which along the conical end approach the identity in $L^{2}_{k}$. The moduli space $\mathcal{M}_{\mathrm{loc}}(Z^+)$ is not a Banach manifold, but carries a natural topology -- that of convergence in $L^2_{k}$ away from infinite cylindrical regions $(- \infty, l ) \times Y$. We use the notation $\mathcal{M}_{\mathrm{loc}}(\sigma)$ for the corresponding space obtained by a fibre product as above.


\subsubsection{A local compactness result}

The exponential decay estimates of Theorem \ref{expdecay1} can be interpreted as telling us that certain energy along the conical end for configurations $(A, \Phi , t , s ) \in \mathcal{M}(\sigma)$ is uniformly bounded, by a constant depending on $\sigma$, and hence that the conical end $K$ behaves like a compact end for the purpose of the compactness analysis of $\mathcal{M}(\sigma)$. We now introduce the relevant notion of energy along the cylindrical energy, and describe the main local compactness result.

We fix $r \geq 1$. Later we will require that $r$ is large enough, depending on $\sigma$ only, so that for all configurations $(A, \Phi, t , s) \in \mathcal{M}(\sigma)$ we have $|\alpha| \geq 1/2$ (using the notation of (\ref{connectiondecompos}-\ref{spinordecompos})) over the portion $[ r , + \infty ) \times Y$ of the conical end $K$. That this can be done follows from the exponential decay estimate of Theorem \ref{expdecay1}. Let $Z_{r}^{+} = Z^+ \setminus ( r , + \infty ) \times Y$. 

Throughout the article we have been fixing an admissible perturbation $\mathfrak{q}$ of the Chern-Simons-Dirac functional on $(Y , \s_{\xi_0 , \alpha_0 , j_0} )$. By the construction in [\cite{KM}, \S 10.1], the admissible perturbation $\mathfrak{q}$ is the formal $L^2$-gradient of some gauge-invariant function $f$ on the configuration space $\mathcal{C}(Y)$.

\begin{Definition}
The \textit{cylindrical energy} of a configuration $\gamma = (A , \Phi , t , s ) \in \mathcal{M}_{z}([\afr])$ is 
\begin{align}
\mathcal{E}_{ r} (\gamma)  =  & \frac{1}{4} \int_{Z_{r}^{+}} F_{\hat{A}} \wedge F_{\hat{A}} - \int_{r \times Y} \inner{ \Phi|_{r \times Y}}{D_{B} \Phi|_{ r \times Y} } \nonumber \\
& +   \int_{r \times Y } (H/ 2) | \Phi |^2  + 2 f ( [ \afr ] ). \end{align}
Above, $B$ denotes the restriction of $A$ to the boundary $\partial Z_{r}^{+} = r \times Y $. By $H$ we denote the mean curvature vector field of the boundary $\partial Z_{r}^{+} = r \times Y$. 
\end{Definition}

That one should just integrate over $Z^{+}_{r} \subset Z^+$ was proposed by B.Zhang (see p.54 \cite{boyu}). The point of cutting off at $r$ is that $\mathcal{E}_{r}(\gamma )$ approches $+ \infty$ as $r$ grows. This can be deduced from Lemma \ref{topan} below. In \cite{KM} this type of energy is called \textit{topological}: the analogous integral over a compact manifold with a cylindrical end attached only depends on the critical point $[\afr]$, the homotopy class $z$ and the chosen perturbation $\mathfrak{q}$. This interpretation is lost in our case, due to the cutting off that is forced upon us, but we do have the identity
\begin{align}
\mathcal{E}_{r}(\gamma ) = 2 \mathrm{CSD}_{\mathfrak{q}}(\afr ) - 2 \mathrm{CSD}(\gamma|_{r} ) + \frac{1}{4} \int_{Z^{+}_{r}} F_{\hat{A}_{0}} \wedge F_{\hat{A}_{0}} \label{csdidentity}
\end{align}
whose terms we describe now.  First recall that for a closed oriented $3$-manifold $Y$ with a spin-c structure $\s$, the Chern-Simons-Dirac functional (see \cite{KM}, \S 4.1) is defined on the configuration space of pairs $(B, \Psi)$ by
\begin{align}
\mathrm{CSD}(B , \Psi ) = - \frac{1}{8} \int_{r \times Y} (\hat{B} - \hat{B}_0 ) \wedge (F_{\hat{B}} + F_{\hat{B}_0}) + \frac{1}{2} \int_{r \times Y} \inner{ D_{B} \Psi}{\Psi}. \label{csd}
\end{align}
The above formula needs the choice of a base spin-c connection $B_0$. Then in formula (\ref{csdidentity}), $\mathrm{CSD}_{\mathfrak{q}} = \mathrm{CSD} + f$ is the $\mathfrak{q}$-perturbed Chern-Simons-Dirac functional for $(Y , \s_{\xi_0 , \alpha_0 , j_0})$ and some choice of base connection $B_0$. The term $\gamma|_r $ is the restriction of $\gamma$ onto the slice $r \times Y$. We have chosen a spin-c connection $A_0$ over $Z^{+}_{r}$ which becomes translation invariant over the cylindrical end with the form $A_0 = d/dt + B_0$, and we use the restriction of $A_0$ onto the slice $r \times Y$ as base connection for the function $\mathrm{CSD}$ on configurations on the slice $r \times Y$. The identity (\ref{csdidentity}) is obtained by integrating by parts as in [\cite{KM}, \S 4.1].

That the cylindrical energy provides a good notion of energy along the cylindrical end is provided by the fact that it controls the $L^2$ norms of $F_{\hat{A}}$, $\Phi$ and $\nabla_A \Phi $ over compact sets:

\begin{Lemma} \label{topan}There exists a constant $C > 0$ depending on $\sigma$, such that for any configuration $\gamma =  (A , \Phi , t , s) \in \mathcal{M}(\sigma)$ we have the following estimate: for any $l \leq 0$ $$\mathcal{E}_{r } (\gamma) \geq \frac{1}{16}\int_{[l , r ] \times Y}( |F_{\hat{A}}|^2 + (|\Phi|^2 - C )^2 + |\nabla_{A} \Phi |^2)  - C ( r  - l  +1) .$$
\end{Lemma}

The proof of the above is analogous to that of Lemma 24.5.1 in \cite{KM}. By an argument as in [\cite{monocont}, pp. 26-27], we can combine Theorem \ref{expdecay1} and Lemma \ref{topan} and obtain, following the standard compactness argument (based on the proof of Theorem 5.1.1 in \cite{KM}), the following local compactness result:
\begin{Proposition}\label{localcompactness-prop}
For any sequence $\gamma_n \in \mathcal{M}(\sigma)$ with uniform bounds $\mathcal{E}_{r} (\gamma_n ) \leq C$ and $-C \leq\pi_{\R}(\gamma_n ) \leq C$, there exist a subsequence which converges in $\mathcal{M}_{\mathrm{loc}}(\sigma )$.
\end{Proposition}

At this point the compactness of the moduli spaces of broken configurations $\mathcal{M}_{z}^{+}([\afr], U , \sigma ; \tau )$, $\mathcal{M}_{z}^{+}([\afr] , \gamma , \sigma ; \kappa )$ or $M_{z}^{+}([\afr] , \sigma )$ follows. We state the result for the first. The broken configurations that can appear in the $1$-dimensional case were listed in Proposition \ref{strata}, and in the general one may see further breaking on the cylindrical end. The statement that we obtain is the following, and its proof follows the arguments of \S16.1 and \S 24.6 of \cite{KM}:

\begin{Corollary}\label{compbroken}
For a fixed $[\afr]$ and $C > 0$, the space of broken configurations $\gamma \in \bigcup_{z} \mathcal{M}_{z}^{+}([\afr] , U , \sigma ; \tau )$ with $\mathcal{E}_{r}(\gamma) \leq C$ is compact. In particular, each $\mathcal{M}_{z}^{+}([\afr] , U , \sigma ; \tau )$ is compact.
\end{Corollary}

Above, the cylindrical energy $\mathcal{E}_{r}$ has been extended to broken configurations $\gamma$ as in \cite{KM}: by adding up the energies of each component of $\gamma$. We recall that the energy of a configuration $\gamma$ in the cylinder moduli space $M_{z}([\afr] , [\bfr] )$ is $2 \cdot ( \mathrm{CSD}_{\mathfrak{q}}(\afr) - \mathrm{CSD}_{\mathfrak{q}}(\bfr ) )$ provided $\gamma$ approaches $\afr$ and $\bfr$ on the corresponding ends. The second assertion in Corollary \ref{compbroken} uses that $\mathcal{E}_r$ is bounded on $\mathcal{M}_{z}^{+}([\afr] , U , \sigma ; \tau )$, which can be seen from (\ref{csdidentity}) combined with Lemma \ref{csdbounds} below.

\begin{Lemma} \label{csdbounds}
There is a constant $C > 0$ depending on $\sigma$ such that for any $\gamma = ( A, \Phi , t , s ) \in \mathcal{M}(U, \sigma ; \tau )$ one has $| \mathrm{CSD}(\gamma|_{r}) - \mathrm{CSD}((A_{t}, \Phi_{t} )|_r ) | \leq C $.
\end{Lemma}

\begin{proof}
This follows from the exponential decay estimates in Theorem \ref{expdecay1} and Corollary \ref{expdecay2}.
\end{proof}

\subsubsection{Finiteness results}

We now outline how to deduce the finiteness result below. This result is the input needed to conclude that the counts of zero dimensional moduli in this paper are indeed finite. We state our results for $\mathcal{M}_{z}([\afr] , U , \sigma ; \tau )$ but the same holds for $\mathcal{M}_{z}^{+}([\afr] , \gamma , \sigma ; \kappa )$ or $M_{z}^{+}([\afr] , \sigma )$.

\begin{Proposition}\label{fin} Suppose that the moduli spaces $\mathcal{M}_{z}([\afr] , U , \sigma ; \tau)$ of expected dimension at most $1$ are transversely cut out. Then there exist only finitely many pairs $([\afr] , z)$ such that the compactified moduli spaces $\mathcal{M}^{+}_{z}([\afr], U , \sigma ; \tau )$ are non-empty and of dimension $\leq 1$.
\end{Proposition}
\begin{Remark}
The reason why the dimension is cut to at most $1$ has to do with the fact that we have been working with $C^l$ contact structures and $C^2$ simplices. This poses a problem if we want that all moduli spaces of all dimensions are transversely cut out after perturbing $\sigma$, due to the assumptions of the Thom-Smale transversality theorem. Raising the differentiability of our data would allow us to conclude the above result for moduli of higher dimensions.
\end{Remark}

The main estimate one needs to prove the above is

\begin{Lemma}[Bounds on energy by dimension] \label{endim}
There exists constants $r \geq 1$  $C > 0$ depending on $\sigma$ such that the following holds. For any $\gamma \in \mathcal{M}_{z}([\afr] , U , \sigma ; \tau)$ we have $$e - C \leq \mathcal{E}_{r}(\gamma) + 4 \pi^2 (\mathrm{gr}_z([\afr] , U , \sigma ; \tau ) - 2\iota([\afr]))    \leq e + C $$
where $\mathrm{gr}_{z}([\afr] , U , \sigma ; \tau )$ denotes the expected dimension of $\mathcal{M}_{z}([\afr] , U , \sigma ; \tau )$, $\iota ([\afr] )$ is defined in [\cite{KM}, p.286] and $e \in \R$ is a constant only depending on $\sigma$ and the image of the critical point $[\afr]$ under the blow-down map $\mathcal{B}_{k-1/2}^{\sigma}(Y)\rightarrow \mathcal{B}_{k-1/2}(Y)$.
\end{Lemma}

\begin{proof}
The corresponding result for the topological energy over a compact manifold with boundary would state that the quantity in the middle, denote it $Q(\gamma )$, only depends on the blow-down of $[\afr]$ (see Proposition 24.6.6 in \cite{KM} and its proof). This time, given two configurations $\gamma \in \mathcal{M}_{z}([\afr] , U , \sigma ; \tau )$ ,  $\tilde{\gamma} \in \mathcal{M}_{\tilde{z}}([\tilde{\afr}] , U , \sigma ; \tau )$ with $[\afr]$ and $[\tilde{\afr}]$ having the same blow-down, their difference in $Q$ can be computed using (\ref{csdidentity}) and we see
\begin{align}
Q(\gamma ) - Q(\tilde{\gamma} )  = -2 \mathrm{CSD}(\gamma|_{r} ) + 2 \mathrm{CSD}( \tilde{\gamma}|_{r} )\label{csddiff}
\end{align} 
We want to establish that $|Q(\gamma ) - Q (\tilde{\gamma})| \leq C$ for a constant $C$ only depending on $\sigma$, and this follows from Lemma \ref{csdbounds}.
\end{proof}

\begin{Lemma}\label{fin1}
Suppose the moduli spaces $\mathcal{M}_{z}([\afr] , U , \sigma ; \tau )$ of dimension $\leq 1$ are transversely cut out. Then for fixed $[\afr]$ there are only finitely many $z$ for which the compactification $\mathcal{M}_{z}^{+}([\afr], U , \sigma ; \tau)$ is non-empty and of dimension $\leq 1$.
\end{Lemma}

\begin{proof}
We note that Lemma \ref{endim} also holds for broken trajectories, with identical proof. For $[\afr]$ and $z$ with $\mathcal{M}^{+}_{z} ([\afr] , U , \sigma ; \tau )$ non-empty, and transversely cut out, we obtain from Lemma \ref{endim} that any broken trajectory $\gamma$ in the moduli space has
\begin{align*}
\mathcal{E}_{r} (\gamma ) \leq C - \mathrm{gr}_{z}([\afr], U , \sigma ; \tau ) + 8 \pi^2 \iota (\afr ) \leq C + 8 \pi^2 \iota (\afr )
\end{align*}
where the second inequality follows from $\mathrm{gr}_{z}([\afr] , U , \sigma ; \tau) \geq 0$ because the moduli is non-empty and transverse. Since $\mathfrak{q}$ is an admissible perturbation there are finitely many critical points in the blow-down, and the quantity $\iota ([\afr] )$ depends on the blow-down of $[\afr]$ only. So we obtain a uniform bound $\mathcal{E}_{r}(\gamma ) \leq C$. Then Corollary \ref{compbroken} yields finitely many such $z$.
\end{proof}

\begin{proof}[Proof of Proposition \ref{fin}]
If the first Chern class of contact structure $\xi_0$, or equivalently that of the spin-c structure $\s_{\xi_0}$, is non-torsion, then there are only finitely-many critical points $\afr$ and the result follows from Lemma \ref{fin1}. 

In the torsion case, we can still argue that there is a bound, independent of $[\afr]$ or $z$, on the cylindrical energy of all broken configurations. Indeed, consider just the case of an unbroken configuration $[\gamma] \in \mathcal{M}_{z}([\afr] , U , \sigma ; \tau )$ and the identity (\ref{csdidentity}) for $\mathcal{E}_{r}(\gamma )$. Since the Chern-Simons-Dirac function is fully gauge-invariant in the torsion case, then there is a bound $|\mathrm{CSD}_{\mathfrak{q}}(\afr )| \leq C$, since $\mathrm{CSD}_{\mathfrak{q}}$ only depends on the blow-down of the critical point $\afr$, for which there are only finitely-many possibilities. Also there is a bound $|\mathrm{CSD}(\gamma|_r )| \leq C$ from applying Lemma \ref{csdbounds}. The remaining term in (\ref{csdidentity}) can also be bounded, so this shows that $\mathcal{E}_r ( \gamma )$ is bounded. The case of a broken configuration is no different.

Now, Lemma \ref{endim} provides upper and lower bounds on $$\mathcal{E}_r (\gamma ) + 4\pi^2 ( \mathrm{dim} \mathcal{M}_{z}([\afr] , U , \sigma ; \tau ) - 2 \iota ([\afr] )).$$ Since we have upper and lower bounds on both the energy and dimension, we obtain $|\iota ( [\afr] )| \leq C$. This gives finitely-many choices for $[\afr]$ again.
\end{proof}

\subsection{Orientations} \label{orientations}

We described in \S \ref{orientability1} and \S \ref{orientations2} the rule for orienting all the moduli spaces in this article, which we called the \textit{canonical orientations}. Whenever these moduli are $0$-dimensional and we use them to make counts of points, each point is counted with a sign corresponding to its canonical orientation (relative to the natural orientation of a point).
The compactifications $M_{z}^{+}([\afr], \sigma)$ and $\mathcal{M}_{z}^{+}([\afr] , U , \sigma ; \tau )$ of $1$-dimensional moduli are $1$-dimensional stratified spaces with a codimension-$1$ $\delta$-structure near its boundary -- a more general form than a manifold with boundary structure (see \cite{KM}, Definition 19.5.3). In this situation each boundary point inherits a \textit{boundary orientation} (see \cite{KM}, Definition 20.5.1) generalising the usual outward-normal first convention for orienting the boundary of a manifold. The total enumeration of the boundary points of the compactified $1$-dimensional moduli equals zero, provided the boundary points are counted with their boundary orientation. 

The next two results compare the canonical and boundary orientations for the relevant moduli spaces. These provide the final touch to the proofs of Proposition \ref{chainmapprop} and Proposition \ref{mainidentities}.

\begin{Lemma} \label{signs1}
Let $M_{z}([\afr] , \sigma )$ be a $1$-dimensional moduli. For each of its codimension-$1$ stratum components listed in Proposition \ref{strata1}, the difference between the canonical and boundary orientation is given by the sign
\begin{enumerate}[(a)]
\item $+1$
\item $(-1)^{\dim M_{z_1} ([\bfr] , [\cfr] )} = -1$.
\item $(-1)^{n-1}(-1)^{i}$ for the moduli over the face $\Delta^{n-1}_i \subset \Delta^n$.
\end{enumerate}
\end{Lemma}
\begin{proof}
(a) and (b) are analogous to cases (i) and (iii) in Proposition 25.2.2 of \cite{KM}. 

For (c) we sketch the main idea. The key result is the following (see \cite{KM}, p.379, formula (20.3)): if $P_1$ and $P_2$ are two Fredholm linear maps of Banach spaces, and the determinant lines $\mathrm{det} P_1$ and $\mathrm{det}P_2$ are oriented, then both $\mathrm{det}( P_1 \oplus P_2)$ and $\mathrm{det}(P_2 \oplus P_1 )$ inherit orientations in a natural way, which under the obvious isomorphism $\mathrm{det}( P_1 \oplus P_2)= \mathrm{det}(P_2 \oplus P_1 )$ differ by the sign $$(-1)^{\mathrm{ind} P_1 \times \mathrm{dim}\, \mathrm{coker} P_2 + \mathrm{dim}\, \mathrm{coker} P_1 \times \mathrm{ind} P_2 }. $$
Suppose now $P_2 = 0_{N} : N \rightarrow 0$ is the zero map out of a finite-dimensional oriented vector space $N$. We also assume $N = \R \times B$ is a product of oriented vector spaces, and that the orientation on $N$ is the product orientation. Then we write $0_{N} = 0_{\R} \oplus 0_{B}$, and the previous result now gives us that the orientations of $\mathrm{det}(P_1 \oplus 0_{\R} \oplus 0_{B})$ and $\mathrm{det}(\oplus_{\R} \oplus P_1 \oplus 0_{B} )$ differ by the sign $(-1)^{\mathrm{dim} \,\mathrm{coker}P_1}$. 

Going back to our case of interest, what we want is to compute the boundary orientation (relative to the canonical orientation) of the boundary stratum component $M_{z}([\afr], \partial \sigma )$ of $ M_{z}([\afr] , \sigma )$, where $\sigma :\Delta^n \rightarrow \mathcal{CM}(Y , \xi_0 ) \times \mathcal{P}$ is a singular simplex of dimension $n$ and $M_z ([\afr] , \sigma )$ is $1$-dimensional. The deformation operator for a configuration $\gamma$ in $M_{z}([\afr] , \sigma )$ (say lying over the point $u \in \Delta^n $) can be transformed by a homotopy to an operator $P \oplus 0_{N}$ where $P$ is the deformation operator for the configuration $\gamma$ in the moduli over a point $M_{z}([\afr] , \sigma (u) )$ and $N := T_u \Delta^n$ is the tangent space to the simplex. When $u $ lies on the interior of the face $\Delta^{n-1}_{i}$, whose boundary orientation is given by the sign $(-1)^{i}$, we decompose $N = \R \oplus B$ where $B= T_u \Delta^{n-1}_{i}$ (with boundary orientation) and $\R$ is the outward-normal direction. The number $\mathrm{dim} \,\mathrm{coker} P$ is $n-1$, which by the above arguments gives the sign (c).
\end{proof}

\begin{Lemma} \label{signs2} Let $\mathcal{M}_{z}([\afr] , U , \sigma ; \tau )$ be a $1$-dimensional moduli. For each of its codimension-$1$ stratum components listed in Proposition \ref{strata}, the difference between the canonical and boundary orientation  is given by the sign
\begin{enumerate}[(a)]
\item $+1$
\item $(-1)^{n}(-1)^{i}$ for the moduli over the face $\Delta^{n-1}_{i} \subset \Delta^n$
\item $+1$
\item $(-1)^{\mathrm{dim} M_{z_1}([\bfr] , [\cfr])} = -1$
\item $-1$
\item $(-1)^{\mathrm{dim} M_{z_1}([\bfr] , [\cfr]) +1} = +1$
\item $(-1)^{\mathrm{dim}M_{z_2}([\afr] , [\bfr] ) + \mathrm{dim}M_{z_1}([\bfr] , U , [\cfr] )} = -1$
\end{enumerate}
\end{Lemma}
\begin{proof}
(a) is clear, and (b) is analogous to (c) of the previous Lemma. (c) and (d) are analogous to cases (i) and (iii) in Proposition 25.2.2 of \cite{KM}, whereas (e), (g) and (f) are analogous to cases (i),(ii) and (iii) in Proposition 26.1.7 of \cite{KM}.
\end{proof}

\printbibliography

\end{document}